\documentclass[11pt]{amsart}

\usepackage{amsmath}
\usepackage{amssymb}

\newtheorem{Theorem}{Theorem}[section]
\newtheorem{Claim}[Theorem]{Claim}
\newtheorem{Fact}[Theorem]{Fact}
\newtheorem{Lemma}[Theorem]{Lemma}
\newtheorem{Proposition}[Theorem]{Proposition}

\newtheorem{Observation}[Theorem]{Observation}

\theoremstyle{definition}
\newtheorem{Definition}[Theorem]{Definition}
\newtheorem{Thesis}[Theorem]{Thesis}
\newtheorem{Problem}[Theorem]{Problem}

\theoremstyle{remark}
\newtheorem{Remark}[Theorem]{Remark}
\newtheorem{Discussion}[Theorem]{Discussion}
\newtheorem{Conclusion}[Theorem]{Conclusion}
\newtheorem{Hypothesis}[Theorem]{Hypothesis}
\newtheorem{notation}[Theorem]{Notation}
\newtheorem{nConclusion}[Theorem]{nice? Conclusion}

\newcount\skewfactor
\def\mathunderaccent#1#2 {\let\theaccent#1\skewfactor#2
\mathpalette\putaccentunder}
\def\putaccentunder#1#2{\oalign{$#1#2$\crcr\hidewidth
\vbox to.2ex{\hbox{$#1\skew\skewfactor\theaccent{}$}\vss}\hidewidth}}

\def\smallbox#1{\leavevmode\thinspace\hbox{\vrule\vtop{\vbox
   {\hrule\kern1pt\hbox{\vphantom{\tt/}\thinspace{\tt#1}\thinspace}}
   \kern1pt\hrule}\vrule}\thinspace}

\newcommand{\Cal}{\mathcal}
\newcommand{\Dom}{{\rm Dom}}
\newcommand{\Rang}{{\rm Rang}}
\newcommand{\pcf}{{\rm pcf}}
\newcommand{\bt}{{\bf t}}
\newcommand{\lev}{{\rm lev}}
\newcommand{\otp}{{\rm otp}}
\newcommand{\inv}{{\rm inv}}
\newcommand{\INV}{{\rm INV}}
\newcommand{\INv}{{\rm INv}}
\newcommand{\Inv}{{\rm Inv}}
\newcommand{\id}{{\rm id}}
\newcommand{\bP}{{\bf P}}
\newcommand{\bU}{{\bf U}}
\newcommand{\nacc}{{\rm nacc}}
\newcommand{\acc}{{\rm acc}}
\newcommand{\Reg}{{\rm Reg}}
\newcommand{\OP}{{\rm OP}}
\newcommand{\IND}{{\rm IND}}
\newcommand{\cf}{{\rm cf}}
\newcommand{\tcf}{{\rm tcf}}
\newcommand{\bh}{{\bf h}}
\newcommand{\cov}{{\rm cov}}
\newcommand{\pp}{{\rm pp}}
\newcommand{\Univ}{{\rm Univ}}

\usepackage[hidelinks]{hyperref}

\begin{document}
\makeatletter\def\shfiuwefootnote{\gdef\@thefnmark{}\@footnotetext}\makeatother\shfiuwefootnote{Version 2002-04-03\_11. See \url{https://shelah.logic.at/papers/552/} for possible updates.}

\title[Non-existence of Universals]{Non-existence of Universals for Classes
like reduced torsion free Abelian groups under embeddings which are not
necessarily pure}
\author{Saharon Shelah}
\address{Institute of Mathematics
 The Hebrew University of Jerusalem
 Jerusalem 91904, Israel
 and  Department of Mathematics
 Rutgers University
 New Brunswick, NJ 08854, USA}
\email{shelah@math.huji.ac.il}
\urladdr{http://www.math.rutgers.edu/\char`\~shelah}

\thanks{I would like to thank Alice Leonhardt for the beautiful
typing of original version.\newline
The research was supported by the German-Israeli Foundation for
Scientific Research\ \&\ Development, Grant No. G-294.081.06/93 .\newline
Publication No 552. \newline
Revised after publication 2001-2}

\date{September 1995}

\begin{abstract}
We consider a class $K$ of structures, e.g.\ trees with $\omega+1$ levels,
\index{$\omega+1$-trees}
metric spaces
\index{metric space}
 and mainly, classes of Abelian groups
\index{abelian group}
 like the one mentioned
in the title
\index{torsion free abelian group}
 and the class of
\index{reduced separable (abelian) $p$-groups}
reduced separable (Abelian) $p$-groups. We say
$M\in K$ is universal for $K$ if any member $N$ of $K$ of cardinality not
bigger than the cardinality of $M$ can be embedded into $M$.  This is a
natural, often raised, problem. We try to draw consequences of
\index{cardinal arithmetic }
cardinal arithmetic
\index{pcf}
 to non--existence of universal
\index{universal}
 members for such natural classes.
\index{forcing}
\end{abstract}

\maketitle

\newpage

\section*{Table of Contents}

\noindent \S0 - INTRODUCTION
\medskip

\noindent \S1 - THEIR PROTOTYPE IS ${\mathfrak K}^{tr}_{\langle\lambda_n:n< \omega
\rangle}$, NOT ${\mathfrak K}^{tr}_\lambda$!
\medskip

\noindent \S2 - ON STRUCTURES LIKE $(\prod\limits_n
\lambda_n,E_m)_{m<\omega}$, $\eta E_m\nu =:\eta(m) = \nu(m)$
\medskip

\noindent \S3 - REDUCED TORSION FREE GROUPS; NON-EXISTENCE OF UNIVERSALS
\medskip

\noindent \S4 - BELOW THE CONTINUUM THERE MAY BE UNIVERSAL STRUCTURES
\medskip

\noindent \S5 - BACK TO ${\mathfrak K}^{rs(p)}$, REAL NON-EXISTENCE RESULTS
\medskip

\noindent \S6 - IMPLICATIONS BETWEEN THE EXISTENCE OF UNIVERSALS
\medskip

\noindent \S7 - NON-EXISTENCE OF UNIVERSALS FOR TREES WITH SMALL DENSITY
\medskip

\noindent \S8 - UNIVERSALS IN SINGULAR CARDINALS
\medskip

\noindent \S9 - METRIC SPACES AND IMPLICATIONS
\medskip

\noindent \S10 - ON MODULES
\medskip

\noindent \S11 - OPEN PROBLEMS
\medskip

\noindent REFERENCES

\eject

\section{Introduction}
{\bf Context.}\hspace{0.15in}  In this paper, model theoretic notions (like
superstable, elementary classes) appear in the introduction but not in the
paper itself (so the reader does not need to know them). Only naive set
theory and basic facts on Abelian groups (all in \cite{Fu}) are necessary
for understanding the paper. The basic definitions are reviewed at the end
of the introduction. On the history of the problem of the existence of
universal members, see Kojman, Shelah \cite{Sh:409}; for more direct
predecessors see Kojman, Shelah \cite{Sh:447}, \cite{Sh:455} and
\cite{Sh:456}, but we do not rely on them. For other advances see
\cite{Sh:457}, \cite{Sh:500} and D\v{z}amonja, Shelah \cite{Sh:614}.
Lately \cite{Sh:622} continues this paper.
\medskip

A class ${\mathfrak K}$ is a class of structures with an embeddability notion.
If not said otherwise, an embedding, is a one to one function preserving
atomic relations and their negations.  If ${\mathfrak K}$ is a class and
$\lambda$ is a cardinal, then ${\mathfrak K}_\lambda$ stands for the collection
of all members of ${\mathfrak K}$ of cardinality $\lambda$.\\
We similarly define ${\mathfrak K}_{\leq\lambda}$.

A member $M$ of ${\mathfrak K}_\lambda$ is universal, if every $N\in {\mathfrak
K}_{\le \lambda}$, embeds into $M$. An example is $M=:\bigoplus\limits_\lambda
{\mathbb Q}$, which is universal in ${\mathfrak K}_\lambda$ if ${\mathfrak K}$ is
the class of all torsion-free Abelian groups, under usual embeddings.

We give some motivation to the present paper by a short review of the above
references. The general thesis in these papers, as well as the present one
is:

\begin{Thesis}
\label{0.1}
General Abelian groups and trees with $\omega+1$ levels behave in universality
theorems like stable non-superstable theories.
\end{Thesis}

The simplest example of such a class is the class ${\mathfrak K}^{tr} =:$ trees
$T$ with $(\omega+1)$-levels, i.e. $T\subseteq {}^{\omega\ge}\alpha$ for some
$\alpha$, with the relations $\eta E^0_n\nu =: \eta\restriction n=\nu
\restriction n\ \&\ \lg(\eta)\geq n$. For ${\mathfrak K}^{tr}$ we know
that $\mu^+<\lambda=\cf(\lambda)
<\mu^{\aleph_0}$ implies there is no universal for ${\mathfrak K}^{tr}_\lambda$
(by \cite{Sh:447}). Classes as ${\mathfrak K}^{rtf}$ (defined in the title),
or ${\mathfrak K}^{rs(p)}$ (reduced separable Abelian $p$-groups) are similar
(though they are not elementary classes) when we consider pure embeddings
(by \cite{Sh:455}). But it is not less natural to consider usual embeddings
(remembering they, the (Abelian) groups under consideration, are reduced). The
problem is that the invariant has been defined using divisibility, and so
under non-pure embedding those seemed to be erased.

Then in \cite{Sh:456} the non-existence of universals is proved restricting
ourselves to $\lambda>2^{\aleph_0}$ and $(< \lambda)$-stable groups
(see there). These restrictions hurt the generality of the theorem; because
of the first requirement we lose some cardinals. The second requirement
changes the class to one which is not established among Abelian group
theorists (though to me it looks natural).

Our aim was to eliminate those requirements, or show that they are necessary.
Note that the present paper is mainly concerned essentially with results in
ZFC, but they have roots in ``difficulties" in extending independence results
thus providing a case for the

\begin{Thesis}
\label{0.2}
Even if you do not like independence results you better look at them, as you
will not even consider your desirable ZFC results when they are camouflaged
by a litany of many independence results you can prove things about.
\end{Thesis}

Of course, independence has interest {\em per se}; still for a given problem in
general a solution in ZFC is for me preferable to an independence result. But
if it gives a method of forcing (so relevant to a series of problems) the
independence result is preferable (of course, I assume there are no other
major differences; the depth of the proof would be of first importance to me).

As occurs often in my papers lately, quotations of {\bf pcf} theory appear.

This paper is also a case of, mainly

\begin{Thesis}
\label{0.3}
Assumption of cases of not GCH, mainly 
at singular (more generally pp$(\lambda)>
\lambda^+$) are ``good", ``helpful" assumptions; i.e. traditionally uses of
GCH proliferate mainly not from conviction but as you can prove many theorems
assuming $2^{\aleph_0}=\aleph_1$ but very few from $2^{\aleph_0}>\aleph_1$,
but assuming $2^{\beth_\omega}>\beth^+_\omega$ is helpful in proving.
\end{Thesis}

Unfortunately, most results \,are only almost in ZFC as they use extremely weak
assumptions from {\bf pcf}, assumptions whose independence is not known. So
practically it is not tempting to try to remove them as they may be
true, and it is unreasonable to try to prove independence results before
independence results on {\bf pcf} will advance.

In \S1 we give an explanation of the earlier difficulties: the problem
(of the existence of universals for
${\mathfrak K}^{rs(p)}$) is not like looking for ${\mathfrak K}^{tr}$ (trees with
$\omega+1$ levels) but for ${\mathfrak K}^{tr}_{\langle\lambda_\alpha:\alpha<
\omega\rangle}$ where
\begin{description}
\item[($\oplus$)]  $\lambda^{\aleph_0}_n<\lambda_{n+1}<\mu$, $\lambda_n$
are regular and $\mu^+<\lambda=\lambda_\omega=\cf(\lambda)<\mu^{\aleph_0}$
and ${\mathfrak K}^{tr}_{\langle\lambda_n:n<\omega\rangle}$ is
\[\{T:T\mbox{ a tree with $\omega+1$ levels, in level $n< \omega$ there are
$\lambda_n$ elements}\}.\]
\end{description}
We also consider ${\mathfrak K}^{tr}_{\langle\lambda_\alpha:\alpha\leq\omega
\rangle}$, which is defined similarly but the level $\omega$ of $T$ is
required to have $\lambda_\omega$ elements.

\noindent For ${\mathfrak K}^{rs(p)}$ this is proved fully, for ${\mathfrak K}^{rtf}$
this is proved for the natural examples.
\medskip

In \S2 we define two such basic examples: one is ${\mathfrak K}^{tr}_{\langle
\lambda_\alpha:\alpha \le \omega \rangle}$, \smallbox{??} and the second is
${\mathfrak K}^{fc}_{\langle \lambda_\alpha:\alpha\leq\omega\rangle}$.
The first is a tree with $\omega+1$ levels;
in the second we have slightly less restrictions. We have $\omega$ kinds of
elements and a function from the $\omega$-th-kind to the $n$th kind. We can
interpret a tree $T$ as a member of the second example: $P^T_\alpha = \{x:x
\mbox{ is of level }\alpha\}$ and
\[F_n(x) = y \quad\Leftrightarrow \quad x \in P^T_\omega\ \&\ y \in
P^T_n\  \&\ y <_T x.\]
For the second we recapture the non-existence theorems.

But this is not one of the classes we considered originally.

In \S3 we return to ${\mathfrak K}^{rtf}$ (reduced torsion free Abelian groups)
and prove the non-existence of universal ones in $\lambda$ if $2^{\aleph_0}
< \mu^+<\lambda=\cf(\lambda)<\mu^{\aleph_0}$ and an additional very weak set
theoretic assumption (the consistency of its failure is not known).
\medskip

\noindent Note that (it will be proved in \cite{Sh:622}):
\begin{description}
\item[($\otimes$)]  if $\lambda<2^{\aleph_0}$ then ${\mathfrak K}^{rtf}_\lambda$
has no universal members.
\end{description}

\noindent Note: if $\lambda=\lambda^{\aleph_0}$ then ${\mathfrak K}^{tr}_\lambda$
has universal member also ${\mathfrak K}^{rs(p)}_\lambda$ (see \cite{Fu}) but not
${\mathfrak K}^{rtf}_\lambda$ (see \cite[Ch IV, VI]{Sh:e}).

\noindent We have noted above that for ${\mathfrak K}^{rtf}_\lambda$ requiring
$\lambda\geq 2^{\aleph_0}$ is reasonable: we can prove (i.e. in ZFC) that
there is no universal member. What about ${\mathfrak K}^{rs(p)}_\lambda$?  By \S1
we should look at ${\mathfrak K}^{tr}_{\langle\lambda_i:i\le\omega\rangle}$,
$\lambda_\omega=\lambda<2^{\aleph_0}$, $\lambda_n<\aleph_0$.

In \S4 we prove the consistency of the existence of universals for
${\mathfrak K}^{tr}_{\langle\lambda_i:i \le\omega\rangle}$
when $\lambda_n\leq \omega$, $\lambda_\omega=\lambda< 2^{\aleph_0}$ but of
cardinality
$\lambda^+$; this is not the original problem but it seems to be a reasonable
variant, and more seriously, it shoots down the hope to use the present
methods of proving non-existence of universals. Anyhow this is
${\mathfrak K}^{tr}_{\langle \lambda_i:i \le\omega\rangle}$ not
${\mathfrak K}^{rs(p)}_{\lambda_\omega}$, so we proceed to reduce this problem to
the previous one under a mild variant of MA. The intentions are to deal with
``there is universal of cardinality $\lambda$" in D\v{z}amonja Shelah
\cite{Sh:614}.

The reader should remember that the consistency of e.g.
\begin{quotation}
{\em
\noindent $2^{\aleph_0}>\lambda>\aleph_0$ and there is no $M$ such that $M\in
{\mathfrak K}^{rs(p)}$ is of cardinality $<2^{\aleph_0}$ and universal for
${\mathfrak K}^{rs(p)}_\lambda$
}
\end{quotation}
is much easier to obtain, even in a wider context (just add many Cohen reals).
\medskip

As in \S4 the problem for ${\mathfrak K}^{rs(p)}_\lambda$ was reasonably resolved
for $\lambda<2^{\aleph_0}$ (and for $\lambda = \lambda^{\aleph_0}$, see
\cite{Sh:455}), we now, in \S5 turn to $\lambda>2^{\aleph_0}$ (and
$\mu,\lambda_n$)
as in $(\oplus)$ above.
As in an earlier proof we use $\langle C_\delta: \delta\in S\rangle$ guessing
clubs for $\lambda$ (see references or later here), so $C_\delta$ is a subset
of $\delta$ (so the invariant depends on the representation of $G$ but this
disappears when we divide by suitable ideal on $\lambda$).
What we do is: rather than trying to code a subset of
$C_\delta$ (for $\bar G=\langle G_i:i<\lambda\rangle$ a representation or
filtration of the structure $G$ as the union of an increasing continuous
sequence of structures of smaller cardinality)
by an element of $G$, we do it, say, by some set $\bar x=\langle x_t:t\in
\Dom(I)\rangle$, $I$ an ideal on $\Dom(I)$ (really by $\bar x/I$). At first
glance if $\Dom(I)$ is infinite we cannot list {\em a priori} all possible
such sequences for a candidate $H$ for being a universal member, as their
number is $\ge\lambda^{\aleph_0}=\mu^{\aleph_0}$. But we can find a family
\[{\Cal F}\subseteq\{\langle x_t:t\in A\rangle:\ A\subseteq\Dom(I),\ A\notin
I,\ x_t\in\lambda\}\]
of cardinality $<\mu^{\aleph_0}$ such that for any $\bar{x}=\langle x_t:t\in
\Dom(I)\rangle$, for some $\bar y\in {\Cal F}$ we have $\bar y=\bar
x\restriction \Dom(\bar y)$.
\medskip

As in \S3 there is such ${\Cal F}$ except when some set theoretic statement
related to {\bf pcf} holds. This statement is extremely strong, also in the
sense that we do not know how to prove its consistency at present. But
again, it seems unreasonable to try to prove its consistency before the
{\bf pcf} problem was dealt with. Of course, we may try to improve the
combinatorics to avoid the use of this statement, but are naturally
discouraged by the possibility that the {\bf pcf} statement can be proved in
ZFC; thus we would retroactively get the non-existence of universals in ZFC.
\medskip

In \S6, under weak {\bf pcf} assumptions, we prove: if there is a universal
member in ${\mathfrak K}^{fc}_\lambda$ then there is one in
${\mathfrak K}^{rs(p)}_\lambda$; so making the connection between the
combinatorial structures and the algebraic ones closer.
\medskip

In \S7 we give other weak {\bf pcf} assumptions which suffice to prove
non-existence of universals in ${\mathfrak K}^x_{\langle\lambda_\alpha:\alpha
\le\omega\rangle}$ (with $x$ one of the ``legal'' values):
$\max\pcf\{\lambda_n:n<\omega\}=\lambda$ and ${\Cal P}(
\{\lambda_n:n<\omega\})/J_{<\lambda}\{\lambda_n:n<\omega\}$ is an infinite
Boolean Algebra (and $(\oplus)$ holds, of course).
\medskip

In \cite{Sh:409}, for singular $\lambda$ results on non-existence of
universals (there on orders) can be gotten from these weak {\bf pcf}
assumptions.
\medskip

In \S8 we get parallel results from, in general, more complicated assumptions.
\medskip

In \S9 we turn to a closely related class: the class of metric spaces with
(one to one) continuous embeddings, similar results hold for it. We also
phrase a natural criterion for deducing the non-existence of universals
from one class to another.
\medskip

In \S10 we deal with modules and in \S11 we discuss the open problems of
various degrees of seriousness.

We thank Mirna D\v{z}amonja and Noam Greenberg for many correction.

The sections are written in the order the research was done.

\begin{notation}
\label{0.4}
Note that we deal with trees with $\omega+1$, levels rather than, say, with
$\kappa+1$, and related situations, as those cases are quite popular. But
inherently the proofs of \S1-\S3, \S5-\S9 work for $\kappa+1$ as well (in
fact, {\bf pcf} theory is stronger).

\noindent For a structure $M$, $\|M\|$ is its cardinality.

\noindent For a model, i.e. a structure, $M$ of cardinality $\lambda$,
where $\lambda$ is regular uncountable, we say that $\bar M$ is a
representation (or filtration) of $M$ if $\bar M=\langle M_i:i<\lambda\rangle$
is an increasing continuous sequence of submodels of cardinality $<\lambda$
with union $M$.

\noindent For a set $A$, we let $[A]^\kappa = \{B:B \subseteq A \mbox{ and }
|B|=\kappa\}$.

\noindent For a set $C$ of ordinals,
$$\acc(C)=\{\alpha\in C: \alpha=\sup(\alpha \cap C)\}, \mbox{(set of
accumulation points)}$$
$$
\nacc(C)=C\setminus \acc(C) \ (=\mbox{ the set of non-accumulation points}).
$$
We usually use $\eta$, $\nu$, $\rho$ for sequences of ordinals; let
$\eta\vartriangleleft\nu$ means $\eta$ is an initial segment of $\nu$.

Let $\cov(\lambda, \mu, \theta, \sigma)= \min\{|{\Cal P}|: {\Cal P}\subseteq
[\lambda]^{<\mu}$, and for every $A\in [\lambda]^{<\theta}$ for some $\alpha<
\sigma$ and $B_i\in {\Cal P}$ for $i< \alpha$ we have $A\subseteq
\bigcup\limits_{i< \alpha} B_i\}$.
\noindent Remember that for an ordinal $\alpha$, e.g. a natural number,
$\alpha=\{\beta:\beta<\alpha\}$.
\end{notation}

\begin{notation}
\noindent ${\mathfrak K}^{rs(p)}$ is the class of (Abelian) groups which are
$p$-groups (i.e. $(\forall x\in G)(\exists n)[p^nx = 0]$) reduced (i.e. have
no divisible non-zero subgroups) and separable (i.e. every cyclic pure
subgroup is a direct summand).  See \cite{Fu}.

\noindent For $G\in{\mathfrak K}^{rs(p)}$ define a norm $\|x\|=\inf\{2^{-n}:
p^n \mbox{ divides } x\}$.  Now every $G\in {\mathfrak K}^{rs(p)}$ has a
basic subgroup $B=\bigoplus\limits_{\scriptstyle n<\omega\atop\scriptstyle
i<\lambda_n} {\mathbb Z} x^n_i$, where $x^n_i$ has order $p^{n+1}$, and every
$x\in G$ can be represented as $\sum\limits_{\scriptstyle n<\omega\atop
\scriptstyle i<\lambda_n} a^n_ix^n_i$, where for each $n$, $w_n(x)=\{i<
\lambda_n:a^n_ix^n_i\ne 0\}$ is finite.

\noindent ${\mathfrak K}^{rtf}$ is the class of Abelian groups which are reduced
and torsion free (i.e. $G \models nx = 0$, $n>0$\qquad
$\Rightarrow\qquad x = 0$).

\noindent For a group $G$ and $A\subseteq G$ let $\langle A\rangle_G$ be the
subgroup of $G$ generated by $A$, we may omit the subscript $G$ if
clear from the context.

\noindent Group will mean an Abelian group, even if not stated
explicitly.

\noindent Let $H\subseteq_{pr} G$ means  $H$ is a pure subgroup of
$G$.

\noindent Let $nG=\{nx: x\in G\}$ and let $G[n]=\{x\in G: nx=0\}$.
\end{notation}

\begin{notation}
${\mathfrak K}$ will denote a class of structures with the same
vocabulary, with a notion of embeddability, equivalently a notion
$\leq_{{\mathfrak K}}$ of submodel.
\end{notation}

\section{Their prototype is ${\mathfrak K}^{tr}_{\langle \lambda_n:n<\omega
\rangle}$ not ${\mathfrak K}^{tr}$!}

If we look for universal member in ${\mathfrak K}^{rs(p)}_\lambda$, thesis
\ref{0.1} suggests to us to think it is basically ${\mathfrak K}^{tr}_\lambda$
(trees with $\omega+1$ levels, i.e. ${\mathfrak K}^{tr}_{\lambda}$ is our
prototype), a way followed in \cite{Sh:455}, \cite{Sh:456}. But, as
explained in the introduction, this does not give an answer for the case
of usual embedding for the family of all such groups. Here we show
that for this case the thesis should be corrected. More concretely,
the choice of the prototype means the choice of what we expect is the
division of the possible classes. That is for a family of classes a
choice of a prototype assert that we believe that they all behave in
the same way.

We show that looking for a universal member $G$ in ${\mathfrak K}^{rs(p)}_\lambda$
is like looking for it among the $G$'s with density $\le\mu$
($\lambda,\mu$, as usual, as in $(\oplus)$ from \S0). For
${\mathfrak K}^{rtf}_\lambda$ we get weaker results which still cover the examples
usually constructed, so showing that the restrictions in \cite{Sh:455} (to
pure embeddings) and \cite{Sh:456} (to $(<\lambda)$-stable groups) were
natural.

\begin{Proposition}
\label{1.1}
Assume that $\mu=\sum\limits_{n<\omega}\lambda_n=\lim\sup\limits_n\lambda_n$,
$\mu\le\lambda\le\mu^{\aleph_0}$, and $G$ is a reduced separable
$p$-group such that
\[|G|=\lambda\quad\mbox{ and }\quad\lambda_n(G)=:\dim((p^n G)[p]/
(p^{n+1}G)[p])\le\mu\]
(this is a vector space over ${\mathbb Z}/p {\mathbb Z}$, hence the dimension is
well defined). \\
{\em Then} there is a reduced separable $p$-group $H$ such that
$|H|=\lambda$, $H$ extends $G$ and $(p^nH)[p]/(p^{n+1}H)[p]$ is a
group of dimension $\lambda_n$ (so if $\lambda_n\geq \aleph_0$, this
means cardinality $\lambda_n$).
\end{Proposition}

\begin{Remark}
\label{1.1A}
So for $H$ the invariants from \cite{Sh:455} are trivial.
\end{Remark}

\proof (See Fuchs \cite{Fu}). We can find $z^n_i$ (for
$n<\omega$, $i<\lambda_n(G)\le\mu$) such that:
\begin{description}
\item[(a)]  $z^n_i$ has order $p^n$,
\item[(b)]  $B=\sum\limits_{n,i}\langle z^n_i \rangle_G$ is a direct sum,
\item[(c)]  $B$ is dense in $G$ in the topology induced by the norm
\[\|x\|=\min\{2^{-n}:p^n \mbox{ divides } x \mbox{ in } G\}.\]
\end{description}
For each $n<\omega$ and $i<\lambda_n(G)$ ($\le\mu$) choose $\eta^n_i\in
\prod\limits_{m<\omega}\lambda_m$, pairwise distinct such that for $(n^1,i^1)
\neq (n^2,i^2)$ for some $n(*)$ we have:
\[\lambda_n \ge \lambda_{n(*)}\qquad \Rightarrow\qquad \eta^{n^1}_{i^1}(n)
\neq \eta^{n^2}_{i^2}(n).\]
Let $H$ be generated by $G$, $x^m_i$ ($i<\lambda_m$, $m<\omega$),
$y^{n,k}_i$ ($i<\lambda_n$, $n<\omega$, $n\le k<\omega)$ freely except for:
\begin{description}
\item[($\alpha$)]  the equations of $G$,
\item[($\beta$)]   $y^{n,n}_i = z^n_i$,
\item[($\gamma$)]  $py^{n,k+1}_i - y^{n,k}_i = x^k_{\eta^n_i(k)}$,
\item[($\delta$)]  $p^{n+1}x^n_i = 0$,
\item[($\varepsilon$)]  $p^{k+1}y^{n,k}_i = 0$.
\end{description}
Now check, (also see the proof of \ref{1.3}).
\hfill$\square_{\ref{1.1}}$

\begin{Definition}
\label{1.2}
\begin{enumerate}
\item $\bt$ denotes a sequence $\langle t_i:i<\omega\rangle$, $t_i$ a
natural number $>1$.
\item For a group $G$ we define
\[G^{[\bt]}=\{x\in G:\bigwedge_{j<\omega}[x\in (\prod_{i<j} t_i)
G]\}.\]
\item We can define a semi-norm $\|-\|_{\bt}$ on $G$
\[\|x\|_{\bt}=\min\{2^{-i}:x\in (\prod_{j<i} t_j)G\}\]
and so the semi-metric
\[d_{\bt}(x,y)=\|x-y\|_{\bt}.\]
\end{enumerate}
\end{Definition}

\begin{Remark}
\label{1.2A}
So, if $\|-\|_{\bt}$ is a norm, $G$ has a completion under $\|-\|_{\bt}$,
which we call $\|-\|_{\bt}$-completion; if $\bt=\langle i!:i<\omega
\rangle$ we refer to $\|-\|_{\bt}$ as $\mathbb Z$-adic norm, and this induces
$\mathbb Z$-adic topology, so we can speak of $\mathbb Z$-adic completion.
\end{Remark}

\begin{Proposition}
\label{1.3}
Suppose that
\begin{description}
\item[($\otimes_0$)]   $\mu=\sum\limits_n\lambda_n$ and $\mu\le\lambda\le
\mu^{\aleph_0}$ for simplicity, $2<2\cdot\lambda_n\le\lambda_{n+1}$ (maybe
$\lambda_n$ is finite!),
\item[($\otimes_1$)]   $G$ is a torsion free group, $|G|=\lambda$, and
$G^{[\bt]}=\{0\}$,
\item[($\otimes_2$)]   $G_0\subseteq G$, $G_0$ is free and $G_0$ is
$\bt$-dense in $G$ (i.e. in the topology induced by the metric $d_{\bt}$),
where $\bt$ is a sequence of primes.
\end{description}
{\em Then} there is a torsion free group $H$, $G\subseteq H$, $H^{[\bt]}
=\{0\}$, $|H|=\lambda$ and, under $d_{\bt}$, $H$ has density $\mu$.
\end{Proposition}

\proof  Let $\{x_i:i<\lambda\}$ be a basis of $G_0$. Let $\eta_i\in
\prod\limits_{n<\omega} \lambda_n$ for $i<\mu$ be distinct such that
$\eta_i(n+1)\geq \lambda_n$ and
\[i\ne j\qquad \Rightarrow\qquad (\exists m)(\forall n)[m \le n \quad
\Rightarrow\quad \eta_i(n) \ne \eta_j(n)].\]
Let $H$ be generated by
\[G,\ \ x^m_i \mbox{ (for $i<\lambda_m$, $m<\omega$), }\ y^n_i \mbox{
(for $i<\mu$, $n<\omega$)}\]
freely except for
\begin{description}
\item[(a)]  the equations of $G$,
\item[(b)]  $y^0_i = x_i$,
\item[(c)]  $t_n\, y^{n+1}_i + y^n_i = x^n_{\eta_i(n)}$.
\end{description}
\medskip

\noindent{\bf Fact A}\hspace{0.15in}  $H$ extends $G$ and is torsion free.

\noindent [Why?  As $H$ can be embedded into the divisible hull of $G$.]
\medskip

\noindent{\bf Fact B}\hspace{0.15in}  $H^{[\bt]}= \{0\}$.

\proof  Let $K$ be a countable pure subgroup of $H$ such that $K^{[\bt]}\ne
\{0\}$. Now without loss of generality $K$ is generated by
\begin{description}
\item[(i)]  $K_1\subseteq G\cap\mbox{ [the $d_{\bt}$--closure of $\langle x_i:
i\in I\rangle_G]$]}$, where $I$ is a countable infinite subset of $\lambda$
and $K_1\supseteq\langle x_i:i\in I\rangle_G$,
\item[(ii)] $y^m_i$, $x^n_j$ for $i\in I$, $m<\omega$ and $(n,j)\in
J$, where $J\subseteq \omega\times \lambda$ is countable and
\[i\in I,\ n<\omega\qquad\Rightarrow\qquad (n,\eta_i(n))\in J.\]
\end{description}
Moreover, the equations holding among those elements are deducible from the
equations of the form
\begin{description}
\item[(a)$^-$]    equations of $K_1$,
\item[(b)$^-$]   $y^0_i=x_i$ for $i \in I$,
\item[(c)$^-$]  $t_n\,y^{n+1}_i+y^n_i=x^n_{\eta_i(n)}$ for $i\in I,n<\omega$.
\end{description}

\noindent We can find $\langle k_i:i<\omega\rangle$ such that
\[[n\ge k_i\  \&\ n\ge k_j\ \&\ i \ne j\qquad \Rightarrow\qquad \eta_i(n)\ne
\eta_j(n)].\]
Let $y \in K[t] \setminus\{0\}$. Then for some $j$, $y\notin
(\prod\limits_{i<j}
t_i)G$, so for some finite $I_0\subseteq I$ and finite $J_0\subseteq J$ and
\[y^* \in\langle\{x_i:i\in I_0\}\cup\{x^n_\alpha:(n,\alpha)\in J_0\}
\rangle_K\]
we have $y-y^*\in(\prod\limits_{i<j} t_i) G$. Without loss of generality $J_0
\cap\{(n,\eta_i(n)):i\in I,\ n\ge k_i\}=\emptyset$. Now there is a
homomorphism $\varphi$ from $K$ into the divisible hull $K^{**}$ of
\[K^* = \langle\{x_i:i\in I_0\}\cup\{x^n_j:(n,j)\in J_0\}\rangle_G\]
such that $\Rang(\varphi)/K^*$ is finite.  This is enough.
\medskip

\noindent{\bf Fact C}\hspace{0.15in} $H_0=:\langle x^n_i:n<\omega,i<\lambda_n
\rangle_H$ is dense in $H$ by $d_{\bt}$.

\proof Straight as each $x_i$ is in the $d_{\bt}$-closure of $H_0$ inside $H$.
\medskip

Noting then that we can increase the dimension easily, we are done.
\hfill$\square_{\ref{1.3}}$

\section{On structures like $(\prod\limits_n \lambda_n,E_m)_{m<\omega}$,
$\eta E_m \nu =: \eta(m)=\nu(m)$}

\begin{Discussion}
\label{2.1}
We discuss the existence of universal members in cardinality $\lambda$,
$\mu^+<\lambda<\mu^{\aleph_0}$, for certain classes of groups. The claims
in \S1 indicate that the problem is similar not to the problem of the
existence of a universal member in ${\mathfrak K}^{tr}_\lambda$ (the class of
trees with $\lambda$ nodes, $\omega+1$ levels) but to the one where the first
$\omega$ levels, are each with $<\mu$ elements. We look more carefully and
see that some variants are quite different.

The major concepts and Lemma (\ref{2.3}) are similar to those of \S3, but
easier. Since detailed proofs are given in \S3, here we give somewhat
shorter proofs.
\end{Discussion}

\begin{Definition}
\label{2.2}
For a sequence $\bar\lambda=\langle\lambda_i:i\le\delta\rangle$ of cardinals
we define:
\begin{description}
\item[(A)] ${\mathfrak K}^{tr}_{\bar \lambda}=\{T:\,T$ is a tree with $\delta
+1$ levels (i.e. a partial order such that

\qquad\quad for $x\in T$, $\lev_T(x)=:\otp(\{y:y<x\})$ is an ordinal
$\le\delta$) such

\qquad\quad that:\quad $\lev_i(T)=:\{x\in T:\lev_T(x)=i\}$  has cardinality
$\le\lambda_i\}$,
\item[(B)] ${\mathfrak K}^{fc}_{\bar\lambda}=\{M:\,M=(|M|,P_i,F_i)_{i\le\delta}$,
$|M|$ is the disjoint union of

\qquad\quad $\langle P_i:i\le\delta\rangle$, $F_i$ is a function from
$P_\delta$ to $P_i$, $\|P_i\|\le\lambda_i$,

\qquad\quad $F_\delta$ is the identity (so can be omitted)$\}$,
\item[(C)] If $[i\le\delta\quad \Rightarrow\quad \lambda_i=\lambda]$ then we
write $\lambda$, $\delta+1$ instead of
$\langle\lambda_i:i\le\delta\rangle$.
\end{description}
\end{Definition}

\begin{Definition}
\label{2.2A}
Embeddings for ${\mathfrak K}^{tr}_{\bar\lambda}$, ${\mathfrak K}^{fc}_{\bar\lambda}$
are defined naturally: for ${\mathfrak K}^{tr}_{\bar\lambda}$ embeddings preserve
$x<y$, $\neg x<y$, $\lev_T(x)=\alpha$; for ${\mathfrak K}^{fc}_{\bar\lambda}$
embeddings are defined just as for models.

If $\delta^1=\delta^2=\delta$ and $[i<\delta\quad\Rightarrow\quad\lambda^1_i
\le\lambda^2_i]$ and $M^\ell\in{\mathfrak K}^{fc}_{\bar\lambda^\ell}$, (or $T^\ell
\in{\mathfrak K}^{tr}_{\bar\lambda^\ell}$) for $\ell=1,2$, then an embedding of
$M^1$ into $M^2$ ($T^1$ into $T^2$) is defined naturally.
\end{Definition}

\begin{Lemma}
\label{2.3}
Assume  $\bar\lambda=\langle\lambda_i:i\le\delta\rangle$ and $\theta$, $\chi$
satisfy (for some $\bar C$):
\begin{description}
\item[(a)]  $\lambda_\delta$, $\theta$ are regular, $\bar C=\langle C_\alpha:
\alpha\in S\rangle$, $\emptyset \not= S\subseteq\lambda=:\lambda_\delta$,
$C_\alpha\subseteq
\alpha$, for every club $E$ of $\lambda$ for some $\alpha$ we have $C_\alpha
\subseteq E$, $\lambda_\delta<\chi\le |C_\alpha|^\theta$ and $\otp(C_\alpha)
\ge\theta$,
\item[(b)]  $\lambda_i\le\lambda_\delta$,
\item[(c)]  there are $\theta$ pairwise disjoint sets $A\subseteq\delta$
such that $\prod\limits_{i\in A}\lambda_i\ge\lambda_\delta$.
\end{description}
{\em Then}
\begin{description}
\item[($\alpha$)] there is no universal member in ${\mathfrak K}^{fc}_{\bar
\lambda}$;\quad moreover
\item[($\beta$)]  if $M^*_\alpha\in {\mathfrak K}^{fc}_{\bar\lambda}$ or even
$M^*_\alpha\in {\mathfrak K}^{fc}_{\lambda_\delta}$ for $\alpha<\alpha^*<\chi$
{\em then} some $M\in {\mathfrak K}^{fc}_{\bar\lambda}$ cannot be embedded into
any $M^*_\alpha$.
\end{description}
\end{Lemma}

\begin{Remark}
\label{2.3A}
Note that clause $(\beta)$ is relevant to our discussion in \S1: the
non-universality is preserved even if we increase the density and,
also, it is witnessed even by non-embeddability in many models.
\end{Remark}

\proof Let $\langle A_\varepsilon:\varepsilon<\theta\rangle$ be as in clause
(c) and let $\eta^\varepsilon_\alpha\in\prod\limits_{i\in A_\varepsilon}
\lambda_i$ for $\alpha<\lambda_\delta$ be pairwise distinct. We fix $M^*_\alpha
\in {\mathfrak K}^{fc}_{\lambda_\delta}$ for $\alpha<\alpha^*<\chi$.

\noindent For $M\in {\mathfrak K}^{fc}_{\bar\lambda}$, let $\bar M=(|M|,P^M_i,
F^M_i)_{i\le\delta}$ and let $\langle M_\alpha: \alpha<
\lambda_\delta\rangle$ be a representation (=filtration) of $M$; for
$\alpha\in S$, $x\in P^M_\delta$, let

\[\begin{aligned}
\inv(x,C_\alpha;\bar M)=\big\{\beta\in C_\alpha:&\mbox{for some }\varepsilon
<\theta\mbox{ and } y\in M_{\min(C_\alpha\setminus (\beta+1))}\\
  &\mbox{we have }\ \bigwedge\limits_{i\in A_\varepsilon} F^M_i(x)=F^M_i(y)\\
  &\mbox{\underbar{but} there is no such } y\in M_\beta\big\}.
\end{aligned}\]
\[\Inv(C_\alpha,\bar M)=\{\inv(x,C_\alpha,\bar M):x\in P^M_\delta\}.\]
\[\INv(\bar M,\bar C)=\langle\Inv(C_\alpha,\bar M):\alpha\in S\rangle.\]
\[\INV(\bar M,\bar C)=\INv(\bar M,\bar C)/\id^a(\bar C).\]
Recall that
\[\id^a(\bar C)=\{T\subseteq\lambda:\mbox{ for some club $E$ of
$\lambda$ for no $\alpha\in T$ is $C_\alpha\subseteq E$}\}.\]
The rest should be clear (for more details see proofs in \S3), noticing
\begin{Fact}
\label{2.3B}
\begin{enumerate}
\item $\INV(\bar M,\bar C)$ is well defined, i.e. if $\bar M^1$, $\bar M^2$
are representations (=filtrations) of $M$ then $\INV(\bar M^1,\bar
C)=\INV(\bar M^2,\bar C)$.
\item $\Inv(C_\alpha,\bar M)$ has cardinality $\le\lambda$.
\item $\inv(x,C_\alpha;\bar M)$ is a subset of $C_\alpha$ of cardinality
$\le \theta$.
\end{enumerate}
\end{Fact}
\hfill$\square_{\ref{2.3}}$

\begin{Conclusion}
\label{2.4}
If $\mu=\sum\limits_{n<\omega}\lambda_n$ and $\lambda^{\aleph_0}_n<
\lambda_{n+1}$ and $\mu^+<\lambda_\omega=\cf(\lambda_\omega)<\mu^{\aleph_0}$,
{\em then} in ${\mathfrak K}^{fc}_{\langle\lambda_\alpha:\alpha\le\omega\rangle}$
there is no universal member and even in ${\mathfrak K}^{fc}_{\langle
\lambda_\omega:\alpha\le\omega\rangle}$ we cannot find a member universal
for it.
\end{Conclusion}

\proof Should be clear or see the proof in \S3.
\hfill$\square_{\ref{2.4}}$

\begin{Claim}
\label{2.4A}
Suppose $T$ is a first order complete theory such that
\begin{description}
\item[$\boxtimes$] there are formulas $\varphi_i(\bar{x},\bar{y}_i)$ for
$i<\delta$, such that 
\begin{enumerate}
\item[(a)] $(\forall \bar{y}',\bar{y}^{''}) (
\bar{y}' \not= \bar{y}^{''} \rightarrow
\neg (\exists \bar{x}) (\varphi_i(\bar{x},\bar{y}') \wedge \varphi_i (\bar{x},
\bar{y}^{''})]$ and
\item[(b)] $\Gamma=\{\varphi_i(x_\eta,\bar{y}_{i,\alpha})^{if(\eta(i)=\alpha)}
:\eta \in {}^\delta\omega,i<\delta, \alpha<\omega\}$ is
consistent with $T$. 
\end{enumerate}
\end{description}
 If $T$ has a universal model in $\lambda$ and
$\lambda_i=\lambda$ for $i<\delta$ \underline{then} 
${\mathfrak K}^k_{\bar{\lambda}}$ has a universal member.
\end{Claim}

\begin{proof}
Should be clear.
\end{proof}

\section{Reduced torsion free groups: Non-existence of universals}
We try to choose torsion free reduced groups and define invariants so that
in an extension to another such group $H$ something survives. To this end
it is natural to stretch ``reduced" near to its limit.

\begin{Definition}
\label{3.1}
\begin{enumerate}
\item ${\mathfrak K}^{tf}$ is the class of torsion free (abelian) groups.
\item ${\mathfrak K}^{rtf}=\{G\in {\mathfrak K}^{tf}:{\mathbb Q}$ is not embeddable into
$G$ (i.e. $G$ is reduced)$\}$.
\item $\bP^*$ denotes the set of primes.
\item For $x\in G$, $\bP(x,G)=:\{p\in\bP^*: \bigwedge\limits_n x\in
p^n G\}$.
\item ${\mathfrak K}^x_\lambda=\{G\in{\mathfrak K}^x:\|G\|=\lambda\}$.
\item If $H\in {\mathfrak K}^{rtf}_\lambda$, we say $\bar H$ is a representation
or filtration of $H$ if $\bar H=\langle
H_\alpha:\alpha<\lambda\rangle$ is increasing continuous and
$H=\bigcup\limits_{\alpha<\lambda} H_\alpha$,
$H\in {\mathfrak K}^{rtf}$ and each $H_\alpha$ has cardinality $<\lambda$.
\end{enumerate}
\end{Definition}

\begin{Proposition}
\label{3.2}
\begin{enumerate}
\item If $G\in {\mathfrak K}^{rtf}$, $x\in G\setminus\{0\}$, $Q\cup\bP(x,G)
\subsetneqq\bP^*$, $G^+$ is the group generated by $G,y,y_{p,\ell}$ ($\ell
<\omega$, $p\in Q$) freely, except for the equations of $G$ and
\[y_{p,0}=y,\quad py_{p,\ell+1}=y_{p,\ell}\quad \mbox{ and }\quad
y_{p,\ell}=z\mbox{ when } z\in G,p^\ell z=x\]
{\em then} $G^+\in {\mathfrak K}^{rtf}$, $G\subseteq_{pr}G^+$ (pure extension).
\item If $G_i\in {\mathfrak K}^{rtf}$ ($i<\alpha$) is $\subseteq_{pr}$-increasing
{\em then} $G_i\subseteq_{pr}\bigcup\limits_{j<\alpha}G_j\in{\mathfrak K}^{rtf}$
for every $i<\alpha$.
\end{enumerate}
\end{Proposition}
The proof of the following lemma introduces a method quite central to this
paper.

\begin{Lemma}
\label{3.3}
Assume that
\begin{description}
\item[$(*)^1_\lambda$] $2^{\aleph_0}+\mu^+<\lambda=\cf(\lambda)<
\mu^{\aleph_0}$,
\item[$(*)^2_\lambda$] for every $\chi<\lambda$, there is $S\subseteq
[\chi]^{\le\aleph_0}$, such that:
\begin{description}
\item[(i)]  $|S|<\lambda$,
\item[(ii)] if $D$ is a non-principal ultrafilter on $\omega$ and $f:D
\longrightarrow\chi$ {\em then} for some $a\in S$ we have
\[\bigcap \{X\in D:f(X)\in a\}\notin D.\]
\end{description}
\end{description}
{\em Then}
\begin{description}
\item[($\alpha$)] in ${\mathfrak K}^{rtf}_\lambda$ there is no universal
member (under usual embeddings (i.e. not necessarily pure)),
\item[($\beta$)]  moreover, \underbar{for any} $G_i\in {\mathfrak
K}^{rtf}_\lambda$, for $i<i^*<\mu^{\aleph_0}$ \underbar{there is} $G\in
{\mathfrak K}^{rtf}_\lambda$ not embeddable into any one of $G_i$.
\end{description}
\end{Lemma}
Before we prove \ref{3.3} we consider the assumptions of \ref{3.3} in
\ref{3.4}, \ref{3.5}.

\begin{Claim}
\label{3.4}
\begin{enumerate}
\item In \ref{3.3} we can replace $(*)^1_\lambda$ by
\begin{description}
\item[$(**)^1_\lambda$ (i)]  $2^{\aleph_0}<\mu<\lambda=\cf(\lambda)<
\mu^{\aleph_0}$,
\item[\qquad(ii)]  there is $\bar C=\langle C_\delta:\delta\in S^*\rangle$
such that $S^*$ is a stationary subset of $\lambda$, each $C_\delta$ is
a subset of $\delta$ with $\otp(C_\delta)$ divisible by $\mu$,
$C_\delta$ closed in $\sup(C_\delta)$ (which is normally $\delta$, but
not necessarily so) and
\[(\forall\alpha)[\alpha\in \nacc(C_\delta)\quad \Rightarrow\quad \cf(\alpha)
>2^{\aleph_0}]\]
(where $\nacc$ stands for ``non-accumulation points''),
and such that $\bar C$ guesses clubs of $\lambda$ (i.e. for every club $E$ of
$\lambda$, for some $\delta\in S^*$ we have $C_\delta\subseteq E$) and
$[\delta\in S^*\quad \Rightarrow\quad \cf(\delta)=\aleph_0]$.
\end{description}
\item In $(*)^1_\lambda$ and in $(*)^2_\lambda$, without loss of generality
$(\forall\theta<\mu)[\theta^{\aleph_0}<\mu]$ and $\cf(\mu)=\aleph_0$.
\end{enumerate}
\end{Claim}

\proof  \ \ \ 1) This is what we actually use in the proof (see below).

\noindent 2) Replace $\mu$ by $\mu'=\min\{\mu_1:\mu^{\aleph_0}_1\ge\mu$
(equivalently $\mu^{\aleph_0}_1=\mu^{\aleph_0}$)$\}$.
\hfill$\square_{\ref{3.4}}$

Compare to, say, \cite{Sh:447}, \cite{Sh:455}; the new assumption is
$(*)^2_\lambda$, note that it is a very weak assumption, in fact it might be
that it is always true.

\begin{Claim}
\label{3.5}
Assume that $2^{\aleph_0}<\mu<\lambda<\mu^{\aleph_0}$ and $(\forall \theta<
\mu)[\theta^{\aleph_0}<\mu]$ (see \ref{3.4}(2)).
Then each of the following is a sufficient condition to $(*)^2_\lambda$:
\begin{description}
\item[($\alpha$)]  $\lambda<\mu^{+\omega_1}$,
\item[($\beta$)]   if ${\mathfrak a}\subseteq\Reg\cap\lambda\setminus\mu$ and
$|{\mathfrak a}|\le 2^{\aleph_0}$ then we can find $h:{\mathfrak a}\longrightarrow
\omega$ such that:
\[\lambda>\sup\{\max\pcf({\mathfrak b}):{\mathfrak b}\subseteq {\mathfrak a}\mbox{
countable, and $h\restriction {\mathfrak b}$ constant}\}.\]
\end{description}
\end{Claim}

\proof Clause $(\alpha)$ implies Clause $(\beta)$: just use any
one-to-one function $h:\Reg\cap\lambda\setminus\mu\longrightarrow\omega$.
\smallskip

Clause $(\beta)$ implies (by \cite[\S6]{Sh:410} + \cite[\S2]{Sh:430}) that
for $\chi<\lambda$ there is $S\subseteq [\chi]^{\aleph_0}$, $|S|<\lambda$
such that for every $Y\subseteq\chi$, $|Y|=2^{\aleph_0}$, we can find
$Y_n$ such that $Y=\bigcup\limits_{n<\omega} Y_n$ and $[Y_n]^{\aleph_0}
\subseteq S$. (Remember: $\mu>2^{\aleph_0}$.) Without loss of generality
(as $2^{\aleph_0} < \mu < \lambda$):
\begin{description}
\item[$(*)$]  $S$ is downward closed.
\end{description}
So if $D$ is a non-principal ultrafilter on $\omega$ and $f:D\longrightarrow
\chi$ then letting $Y=\Rang(f)$ we can find $\langle Y_n:n<\omega\rangle$ as
above. Let $h:D\longrightarrow\omega$ be defined by $h(A)=\min\{n:f(A)\in
Y_n\}$. So
\[X\subseteq D\ \ \&\ \ |X|\le\aleph_0\ \ \&\ \ h\restriction X
\mbox{ constant }\Rightarrow\ f''(X)\in S\quad\mbox{(remember
$(*)$)}.\]
Now for each $n$, for some countable $X_n\subseteq D$ (possibly finite or
even empty) we have:
\[h \restriction X_n\ \mbox{ is constantly } n,\]
\[\ell <\omega \ \&\ (\exists A\in D)(h(A)=n\ \&\ \ell\notin A)
\Rightarrow (\exists B\in X_n)(\ell\notin B).\]
Let $A_n=:\bigcap\{A:A\in X_n\}=\bigcap\{A:A\in D\mbox{ and } h(X)=n\}$.
If the desired conclusion fails, then $\bigwedge\limits_{n<\omega}A_n\in D$.
So
\[(\forall A)[A\in D\quad \Leftrightarrow\quad\bigvee_{n<\omega} A\supseteq
A_n].\]
So $D$ is generated by $\{A_n:n<\omega\}$ but then $D$ cannot be a
non-principal ultrafilter.
\hfill$\square_{\ref{3.5}}$
\medskip

\begin{proof} 
\underline{of Lemma \ref{3.3}}
\newline 
Let $\bar C=\langle C_\delta:\delta\in S^*\rangle$
be as in $(**)^1_{\bar \lambda}$ (ii) from \ref{3.4} (for \ref{3.4}(1) its
existence is obvious, for \ref{3.3} - use \cite[VI,old III 7.8]{Sh:e}).
Let us suppose that $\bar A=\langle A_\delta:\delta\in S^*\rangle$, $A_\delta
\subseteq\nacc(C_\delta)$ has order type $\omega$ ($A_\delta$ like this will
be chosen later) and let $\eta_\delta$ enumerate $A_\delta$ increasingly. Let
$G_0$ be generated by $\{x_i:i<\lambda\} \cup \{x_{i,p,\ell}:i=(i-p)+
p$ that is $(\exists j) (i=j+p), p \in {\bf P}^*$ and $i<\omega\}$ freely
except 
$$
x_{i,p,0}=x_i, p x_{i,p,\ell+1}= x_{i,p,\ell}.
$$
Let $R$ be
\[\begin{aligned}
\big\{\bar a: &\bar a=\langle a_n:n<\omega\rangle\mbox{ is a sequence of
pairwise disjoint subsets of } \bP^*,\\
  &\mbox{with } a_n \cap \{0,\ldots, n-1\}=\emptyset\mbox{ such that}\\
  &\mbox{for infinitely many }n,\ a_n\ne\emptyset\big\}.
\end{aligned}\]
Let $G$ be a group generated by
\[G_0 \cup \{y^{\alpha,n}_{\bar a},z^{\alpha,n}_{\bar a,p}:\ \alpha<\lambda,
\ \bar a\in R,\ n<\omega,\ p \mbox{ prime}\}\]
freely except for:
\begin{description}
\item[(a)]  the equations of $G_0$,
\item[(b)]  $pz^{\alpha,n+1}_{\bar a,p}=z^{\alpha,n}_{\bar a,p}$ when
$p\in a_n$, $\alpha<\lambda$,
\item[(c)]  $z^{\delta,0}_{\bar a,p}=y^{\delta,n}_{\bar a}-
x_{\eta_\delta(n)+n}$ when $p\in a_n$ and $\delta\in S^*$.
\end{description}
Now $G\in {\mathfrak K}^{rtf}_\lambda$ by inspection.

For $\alpha<\lambda$ let $G_\alpha$ be the subgroup of $G$ generated by
$\{x_i,x_{i,p,\ell}:i<\alpha, \ell<\omega$ and $i=(i-p)+p\} \cup 
\{y^{\delta,n}_{\bar{a}}: \delta \in S^* \cap \alpha, \bar{a}\in R,n<\omega\}
\cup \{z^{\delta,n}_{\bar{a},p}:\delta\in S^* \cap \alpha, \bar{a} \in R,
n<\omega, p\in a_n\}$
\end{proof}

\noindent Before continuing the proof of \ref{3.3} we present a definition
and some facts.

\begin{Definition}
\label{3.7}
For a representation $\bar H$ of $H\in {\mathfrak K}^{rtf}_\lambda$, and $x\in H$,
$\delta\in S^*$ let
\begin{enumerate}
\item $\inv(x,C_\delta;\bar H)=:\{\alpha\in C_\delta:$ for some $Q\subseteq
\bP^*$, there is $y\in H_{\min[C_\delta\setminus(\alpha+1)]}$ such that
$Q\subseteq\bP(x-y,H)$ but for no $y\in H_\alpha$ is $Q\subseteq\bP(x-y,H)\}$

(so $\inv(x,C_\delta;\bar H)$ is a subset of $C_\delta$ of cardinality $\le
2^{\aleph_0}$).
\item  $\Inv^0(C_\delta,\bar H)=:\{\inv(x,C_\delta;\bar H):x\in
\bigcup\limits_i H_i\}$.
\item $\Inv^1(C_\delta,\bar H)=:\{a:a\subseteq C_\delta$ countable and for
some $x\in H$, $a\subseteq\inv(x,C_\delta;\bar H)\}$.
\item $\INv^\ell(\bar H,\bar C)=:\Inv^\ell(H,\bar H,\bar C)=:\langle
\Inv^\ell(C_\delta;\bar H):\delta\in S^*\rangle$ for $\ell\in\{0,1\}$.
\item $\INV^\ell(H,\bar C)=:\INv^\ell(H,\bar H,\bar C)/\id^a(\bar C)$,
where
\[\id^a(\bar C)=:\{T\subseteq\lambda:\mbox{ for some club $E$ of
$\lambda$ for no $\delta\in T$ is }C_\delta\subseteq E\}.\]
\item If $\ell$ is omitted, $\ell = 0$ is understood.
\end{enumerate}
\end{Definition}

\begin{Fact}
\label{3.8}
\begin{enumerate}
\item $\INV^\ell(H,\bar C)$ is well defined.
\item The $\delta$-th component of $\INv^\ell(\bar H,\bar C)$ is a family
of $\le\lambda$ subsets of $C_\delta$ each of cardinality $\le 2^{\aleph_0}$
and if $\ell=1$ each member is countable and the family is closed under
subsets.
\item {\em If} $G^*_i\in{\mathfrak K}^{rtf}_\lambda$ for $i<i^*$, $i^*<
\mu^{\aleph_0}$, $\bar G^i=\langle\bar G_{i,\alpha}:\alpha<\lambda\rangle$
is a representation of $G^*_i$,

{\em then} we can find $A_\delta\subseteq\nacc(C_\delta)$ of order type
$\omega$ such that: $i<i^*$, $\delta\in S^*\qquad \Rightarrow$\qquad for
no $a$ in the $\delta$-th component of $\INv^\ell(G^*_i,\bar G^i,\bar C)$ do
we have $|a \cap A_\delta|\ge\aleph_0$.
\end{enumerate}
\end{Fact}

\proof Straightforward.  (For (3) note $\otp(C_\delta)\ge\mu$, so there
are $\mu^{\aleph_0}>\lambda$ pairwise almost disjoint subsets of
$C_\delta$ each
of cardinality $\aleph_0$ and every $A\in\Inv(C_\delta,\bar G^i)$
disqualifies at most $2^{\aleph_0}$ of them.)
\hfill$\square_{\ref{3.8}}$

\begin{Fact}
\label{3.9}
Let $G$ be as constructed above for $\langle A_\delta:\delta\in
S^*\rangle,A_\delta\subseteq\nacc(C_\delta)$, $\otp(A_\delta)=\omega$
(where $\langle A_\delta:\delta\in S^*\rangle$ are chosen as in \ref{3.8}(3)
for the sequence $\langle G^*_i:i<i^* \rangle$ given for proving
\ref{3.3}, see $(\beta)$ there).

\noindent Assume $G \subseteq H\in {\mathfrak K}^{rtf}_\lambda$ and $\bar H$ is a
filtration of $H$. {\em Then}
\[\begin{array}{rr}
B=:\big\{\delta:A_\delta\mbox{ has infinite intersection with some}&\ \\
a\in\Inv(C_\delta,\bar H)\big\}&=\ \lambda\ \mod\ \id^a(\bar C).
\end{array}\]
\end{Fact}

\proof We assume otherwise and derive a contradiction. Let for $\alpha
<\lambda$, $S_\alpha\subseteq [\alpha]^{\le \aleph_0}$, $|S_\alpha|<\lambda$
be as guaranteed by $(*)^2_\lambda$.

Let $\chi>2^\lambda$, ${\mathfrak A}_\alpha\prec ({\mathcal H}(\chi),\in,
<^*_\chi)$ for $\alpha<\lambda$ increasing continuous, $\|{\mathfrak A}_\alpha\|
<\lambda$, $\langle {\mathfrak A}_\beta:\beta\le\alpha\rangle\in {\mathfrak
A}_{\alpha+1}$, ${\mathfrak A}_\alpha\cap\lambda$ an ordinal and: 
\[\langle S_\alpha:\alpha<\lambda\rangle,\ G,\ H,\ \bar C,\ \langle
A_\delta:\delta\in S^* \rangle,\ \bar H,\ \langle x_i, x_{i,p,\ell},
y^{\delta,n}_{\bar a}, z^{\delta,n}_{\bar a,p}:\;i,\delta,\bar a,n,p \rangle\]
all belong to ${\mathfrak A}_0$ and $2^{\aleph_0}+1\subseteq {\mathfrak A}_0$.
Then $E=\{\delta<\lambda:{\mathfrak A}_\delta \cap\lambda=\delta\}$ is a club
of $\lambda$, note that $\delta\in E \Rightarrow H_\delta \cap G=G_\delta$.
 Choose $\delta\in S^* \cap E\setminus B$ such that
$C_\delta \subseteq E$. (Why can we? As $\id^a(\bar C)$ contains all non
stationary subsets of $\lambda$, in particular $\lambda\setminus E$,
and $\lambda\setminus S^*$ and $B$, but $\lambda\notin \id^a(\bar
C)$.) Remember that $\eta_\delta$ enumerates $A_\delta$ (in the
increasing order). For each $\alpha=\eta_\delta(n)\in A_\delta$ (so
$\alpha\in E$
hence ${\mathfrak A}_\alpha \cap \lambda=\alpha$ but $\bar H\in {\mathfrak
A}_\alpha$ hence $H\cap {\mathfrak A}_\alpha= H_\alpha$) and
$Q\subseteq\bP^*$ choose, if possible, $y_{\alpha,Q}\in H_\alpha$ such that:
\[Q\subseteq\bP(x_{\alpha+n}-y_{\alpha,Q},H).\]
Let $I_\alpha=:\{Q\subseteq\bP^*:y_{\alpha,Q}$ well defined$\}$. Note (see
\ref{3.4} $(**)^1_\lambda$ and remember $\eta_\delta(n)\in A_\delta\subseteq
\nacc(C_\delta)$) that $\cf(\alpha)>2^{\aleph_0}$ (by (ii) of
\ref{3.4} $(**)^1_\lambda$) and hence for some $\beta_\alpha<\alpha$,
\[\{ y_{\alpha,Q}:Q\in I_\alpha\}\subseteq H_{\beta_\alpha}.\]
Now:
\begin{description}
\item[$\otimes_1$]  $I_\alpha$ is a downward closed family of subsets of
$\bP^*$, $\bP^*\notin I_\alpha$ for $\alpha
\in A_\delta$.
\end{description}
[Why? See the definition for the first phrase and note also that $H$ is
reduced for the second phrase.]
\begin{description}
\item[$\otimes_2$]  $I_\alpha$ is closed under unions of two members (hence
is an ideal on $\bP^*$).
\end{description}
[Why? If $Q_1,Q_2\in I_\alpha$ then (as $x_{\alpha+n}\in G\subseteq H$
witnesses
this):
\[\begin{aligned}
({\Cal H}(\chi),\in,<^*_\chi)\models &(\exists x)(x\in H\ \&\
Q_1\subseteq\bP(x- y_{\alpha,Q_1},H)\ \&\\
  &Q_2\subseteq\bP(x-y_{\alpha,Q_2},H)).
\end{aligned}\]
All the parameters are in ${\mathfrak A}_\alpha$ so there is $y\in
{\mathfrak A}_\alpha\cap H$ such that
\[Q_1\subseteq\bP(y-y_{\alpha,Q_1},H)\quad\mbox{ and }\quad Q_2\subseteq
\bP(y-y_{\alpha,Q_2},H).\]
By algebraic manipulations,
\[Q_1\subseteq \bP(x_{\alpha+n}-y_{\alpha,Q_1},H),\ Q_1\subseteq\bP(y-y_{\alpha,
Q_1},H)\quad\Rightarrow\quad Q_1\subseteq\bP(x_{\alpha+n}-y,H);\]
similarly for $Q_2$. So $Q_1\cup Q_2\subseteq\bP(x_{\alpha+n}-y,H)$ and hence
$Q_1\cup Q_2\in I_\alpha$.]

\begin{description}
\item[$\otimes_3$] If $\bar Q=\langle Q_n:n\in\Gamma\rangle$ are pairwise
disjoint subsets of $\bP^*$ and $Q_n$ disjoint to 
$\{0,\ldots,n\}$, for some infinite
$\Gamma\subseteq\omega$, then for some $n\in\Gamma$ we have $Q_n\in
I_{\eta_\delta(n)}$.
\end{description}
[Why?  Otherwise let $a_n$ be $Q_n$ if $n\in \Gamma$, and $\emptyset$
if $n\in \omega\setminus \Gamma$, and let  $\bar a=\langle a_n: n<
\omega\rangle$. Now
$n\in\Gamma\quad\Rightarrow\quad\eta_\delta(n)\in
\inv(y^{\delta, 0}_{\bar a},C_\delta;\bar H)$ and hence
\[A_\delta\cap\inv(y^{\delta, 0}_{\bar a},C_\delta;\bar
H)\supseteq\{\eta_\delta(n):n\in \Gamma\},\]
which is infinite, contradicting the choice of $A_\delta$.]

\begin{description}
\item[$\otimes_4$] for all but finitely many $n$ the Boolean algebra
${\Cal P}(\bP^*)/I_{\eta_\delta(n)}$ is finite.
\end{description}
[Why? If not, then by $\otimes_1$ second phrase, for each $n$ there are
infinitely many non-principal ultrafilters $D$ on $\bP^*$ disjoint to
$I_{\eta_\delta(n)}$, so for $n<\omega$ we can find an ultrafilter $D_n$
on $\bP^*$ disjoint to $I_{\eta_\delta(n)}$, distinct from $D_m$ for
$m<n$. Thus we can find $\Gamma\in [\omega]^{\aleph_0}$ and $Q_n\in D_n$
for $n\in\Gamma$ such that $\langle Q_n:n\in\Gamma\rangle$ are pairwise
disjoint (as $Q_n\in D_n$ clearly $|Q_n|=\aleph_0$) and without loss of
generality $Q_n\cap\{0,\ldots,n\}=\emptyset$. Why? Look: if $B_n
\in D_0\setminus D_1$ for $n\in\omega$ then
\[(\exists^\infty n)(B_n \in D_n)\quad\mbox{ or }\quad(\exists^\infty n)
(\bP^*\setminus B_n \in D_n),\]
etc. Let $Q_n=\emptyset$ for $n\in\omega\setminus\Gamma$, now $\bar Q=
\langle Q_n:n<\omega\rangle$ contradicts $\otimes_3$.]

\begin{description}
\item[$\otimes_5$] If the conclusion (of \ref{3.9}) fails, then for no
$\alpha\in A_\delta$ do we have $I_\alpha\cap [\omega]^{\aleph_0}=
\emptyset \& {\Cal P}(\bP^*)/I_\alpha$ finite.
\end{description}
[Why? If not, choose such an $\alpha$ and $Q^*\subseteq\bP^*$, $Q^*
\notin I_\alpha$ such that $I=I_\alpha\restriction Q^*$ is a maximal
ideal on $Q^*$. So $D=:{\Cal P}(Q^*)\setminus I$ is a non-principal
ultrafilter. Remember $\beta=\beta_\alpha<\alpha$ is such that
$\{y_{\alpha,Q}:Q\in I_\alpha\}\subseteq H_\beta$. Now, $H_\beta\in
{\mathfrak A}_{\beta+1}$, $|H_\beta|<\lambda$. Hence $(*)^2_\lambda$ from
\ref{3.3} (note that it does not matter whether we consider an ordinal
$\chi<\lambda$ or a cardinal $\chi<\lambda$, or any other set of
cardinality $< \lambda$) implies that there is $S_{H_\beta}\in
{\mathfrak A}_{\beta+1}$, $S_{H_\beta}\subseteq [H_\beta]^{\le \aleph_0}$,
$|S_{H_\beta}|<\lambda$ as there. Now it does not matter if we deal with
functions from an ultrafilter on $\omega$ \underbar{or} an ultrafilter on
$Q^*$. We define $f:D\longrightarrow H_\beta$ as follows: for $U\in D$ we
let $f(U)=y_{\alpha,Q^* \setminus U}$.  (Note: $Q^*\setminus U\in I_\alpha$,
hence $y_{\alpha,Q^* \setminus U}$ is well defined.) So, by the choice of
$S_{H_\beta}$ (see (ii) of $(*)^2_\lambda$), for some countable $f'
\subseteq f$, $f'\in {\mathfrak A}_{\beta+1}$ and $\bigcap\{U:U\in\Dom(f')\}
\notin D$ (reflect for a minute). Let $\Dom(f')=\{U_0,U_1,\ldots\}$.
Then $\bigcup\limits_{n<\omega}(Q^*\setminus U_n)\notin I_\alpha$. But as
in the proof of $\otimes_2$, as
\[\langle y_\alpha,(Q^* \setminus U_n):n<\omega\rangle\in
{\mathfrak A}_{\beta+1}\subseteq {\mathfrak A}_\alpha,\]
we have $\bigcup\limits_{n<\omega}(Q^*\setminus U_n)\in I_\alpha$, an
easy contradiction.]

Now recalling $\{0,\ldots,n\} \in I_{\eta_\delta (n)}$ by the choice of $G_0$
and $I_{\eta_\delta(n)}$ clearly $\otimes_4$,$\otimes_5$ give a contradiction.
\hfill$\square_{\ref{3.3}}$

\begin{Remark}
\label{3.10}
We can deal similarly with $R$-modules, $|R|<\mu$ \underbar{if} $R$ has
infinitely many prime ideals $I$. Also the treatment of
${\mathfrak K}^{rs(p)}_\lambda$ is similar to the one for modules over
rings with one prime.

\noindent Note: if we replace ``reduced" by
\[x\in G \setminus\{0\}\quad \Rightarrow\quad (\exists p\in\bP^*)(x\notin
pG)\]
then here we could have defined
\[\bP(x,H)=:\{p\in \bP^*:x\in pH\}\]
and the proof would go through with no difference (e.g. choose a fixed
partition $\langle \bP^*_n: n< \omega\rangle$ of $\bP^*$ to infinite
sets, and let $\bP'(x, H)=\{n: x\in pH\mbox{ for every }p\in
\bP^*_n\}$). Now the groups are less divisible.
\end{Remark}

\begin{Remark}
\label{3.11}
We can get that the groups are slender, in fact, the construction gives it.
\end{Remark}

\section{Below the continuum there may be universal structures}
Both in \cite{Sh:456} (where we deal with universality for $(<\lambda)$-stable
(Abelian) groups, like ${\mathfrak K}^{rs(p)}_\lambda$) and in \S3, we restrict
ourselves to $\lambda>2^{\aleph_0}$, a restriction which does not appear
in \cite{Sh:447}, \cite{Sh:455}. Is this restriction necessary?  In this
section we shall show that at least to some extent, it is.

We first show under MA that for $\lambda<2^{\aleph_0}$, any $G\in
{\mathfrak K}^{rs(p)}_\lambda$ can be embedded into a ``nice" one; our aim
is to reduce the consistency of ``there is a universal in
${\mathfrak K}^{rs(p)}_\lambda$" to ``there is a universal in
${\mathfrak K}^{tr}_{\langle\aleph_0:n<\omega\rangle\char 94\langle \lambda
\rangle}$". Then we proceed to
prove the consistency of the latter. Actually a weak form of MA suffices.

\begin{Definition}
\label{4.2}
\begin{enumerate}
\item $G\in {\mathfrak K}^{rs(p)}_\lambda$ is {\em tree-like} if:
\begin{description}
\item[(a)]  we can find a basic subgroup $B=
\bigoplus\limits_{\scriptstyle i<\lambda_n\atop\scriptstyle n<\omega}
{\mathbb Z} x^n_i$, where
\[\lambda_n=\lambda_n(G)=:\dim\left((p^nG)[p]/p^{n+1}(G)[p]\right)\]
(see Fuchs \cite{Fu}) such that: ${\mathbb Z} x^n_i \cong {\mathbb Z}/p^{n+1}
{\mathbb Z}$ and
\begin{description}
\item[$\otimes_0$] every $x\in G$ has the form
\[\sum\limits_{n,i}\{a^n_i p^{n-k} x^n_i:n\in [k,\omega)\mbox{ and }
i<\lambda\}\]
where $a^n_i\in{\mathbb Z}$ and
\[n<\omega\quad \Rightarrow\quad w_n[x]=:\{i:a^n_i\, p^{n-k}x^n_i\neq 0\}
\mbox{ is finite}\]
\end{description}
(this applies to any $G\in {\mathfrak K}^{rs(p)}_\lambda$ we considered so far;
we write $w_n[x]=w_n[x,\bar Y]$ when $\bar Y=\langle x^n_i:n,i\rangle$).
Moreover
\item[(b)]  $\bar Y=\langle x^n_i:n,i\rangle$ is tree-like inside $G$,
which means
that we can find $F_n:\lambda_{n+1}\longrightarrow\lambda_n$ such that
letting $\bar F=\langle F_n:n<\omega\rangle$, $G$ is generated by some
subset of
$\Gamma(G,\bar Y,\bar F)$ where:
\[\hspace{-0.5cm}\begin{aligned}
\Gamma(G,\bar Y,\bar F)=\big\{x:&\mbox{for some }\eta\in\prod\limits_{n
<\omega}\lambda_n, \mbox{ for each } n<\omega \mbox{ we have}\\
  &F_n(\eta(n+1))=\eta(n)\mbox{ and }x=\sum\limits_{n\ge k}
p^{n-k}x^n_{\eta(n)}\big\}.
\end{aligned}\]
\end{description}
\item $G\in {\mathfrak K}^{rs(p)}_\lambda$ is {\em semi-tree-like} if above
we replace (b) by
\begin{description}
\item[(b)$'$]  we can find a set $\Gamma\subseteq\{\eta:\eta$ is a partial
function from $\omega$ to $\sup\limits_{n<\omega} \lambda_n$ with
$\eta(n)< \lambda_n\}$ such that:
\begin{description}
\item[($\alpha$)]  $\eta_1\in\Gamma,\ \eta_2\in\Gamma,\ \eta_1(n)=\eta_2(n)
\quad \Rightarrow\quad\eta_1\restriction n=\eta_2\restriction n$,
\item[($\beta$)]   for $\eta\in\Gamma$ and $n\in\Dom(\eta)$, there is
\[y_{\eta,n}=\sum \{p^{m-n}x^m_{\eta(m)}:m \in \Dom(\eta)\mbox{ and }
m \ge n\}\in G,\]
\item[($\gamma$)]  $G$ is generated by
\[\{x^n_i:n<\omega,i<\lambda_n\}\cup\{y_{\eta,n}:\eta\in\Gamma,n\in
\Dom(\eta)\}.\]
\end{description}
\end{description}
\item $G\in {\mathfrak K}^{rs(p)}_\lambda$ is {\em almost tree-like} if
in (b)$'$ we add
\begin{description}
\item[($\delta$)]  for some $A\subseteq\omega$ for every $\eta\in\Gamma$,
$\Dom(\eta)=A$.
\end{description}
\end{enumerate}
\end{Definition}

\begin{Proposition}
\label{4.3}
\begin{enumerate}
\item Suppose $G\in {\mathfrak K}^{rs(p)}_\lambda$ is almost tree-like, as
witnessed by $A\subseteq\omega$, $\lambda_n$ (for $n<\omega$), $x^n_i$ (for
$n\in A$, $i<\lambda_n$), and if $n_0<n_2$ are successive members of $A$,
$n_0<n<n_2$ then $\lambda_n\ge\lambda_{n_0}$ or just
\[\lambda_n\ge|\{\eta(n_0):\eta\in\Gamma\}|.\]
{\em Then} $G$ is tree-like (possibly with other witnesses).
\item  If in \ref{4.2}(3) we just demand $\eta\in\Gamma\quad\Rightarrow
\quad\bigvee\limits_{n<\omega}\Dom(\eta)\setminus n=A\setminus n$;
then changing the $\eta$'s and the $y_{\eta,n}$'s we can regain the
``almost tree-like".
\end{enumerate}
\end{Proposition}

\proof 1) For every successive members $n_0<n_2$ of $A$ for
\[\alpha\in S_{n_0}=:\{\alpha:(\exists\eta)[\eta\in\Gamma\ \&\ \eta(n_0)
=\alpha]\},\]
choose ordinals $\gamma(n_0,\alpha,\ell)$ for $\ell\in (n_0,n_2)$ such that
\[\gamma(n_0,\alpha_1,\ell)=\gamma(n_0,\alpha_2,\ell)\quad\Rightarrow
\quad\alpha_1=\alpha_2.\]
We change the basis by replacing for $\alpha\in S_{n_0}$, $\{x^n_\alpha\}\cup
\{x^\ell_{\gamma(n_0,\alpha,\ell)}:\ell\in (n_0,n_2)\}$ (note: $n_0<n_2$
but possibly $n_0+1=n_2$), by:
\[\begin{aligned}
\biggl\{ x^{n_0}_\alpha + px^{n_0+1}_{\gamma(n_0,\alpha,n_0+1)},
&x^{n_0+1}_{\gamma(n_0,\alpha,n_0+1)} +
px^{n_0+2}_{\gamma(n_0,\alpha,n_0+2)},\ldots, \\
  &x^{n_2-2}_{\gamma(n_0,\alpha,n_2-2)} +
px^{n_2-1}_{\gamma(n_0,\alpha,n_2-1)},
x^{n_2-1}_{\gamma(n_0,\gamma,n_2-1)} \biggr\}.
\end{aligned}\]

2) For $\eta\in \Gamma$ let $n(\eta)=\min\{ n: n\in A\cap\Dom(\eta)$
and $\Dom(\eta)\setminus n=A\setminus n\}$, and let
$\Gamma_n=\{\eta\in \Gamma: n(\eta)=n\}$ for $n\in A$. We choose by
induction on $n< \omega$ the objects $\nu_\eta$ for $\eta\in \Gamma_n$
and $\rho^n_\alpha$ for $\alpha< \lambda_n$ such that: $\nu_\eta$ is a
function with domain $A$, $\nu_\eta\restriction (A\setminus
n(\eta))=\eta\restriction (A\setminus n(\eta))$ and
$\nu_\eta\restriction (A\cap n(\eta))= \rho^n_{\eta(n)}$,
$\nu_\eta(n)< \lambda_n$ and $\rho^n_\alpha$ is a function with domain
$A\cap n$, $\rho^n_\alpha(\ell)< \lambda_\ell$ and $\rho^n_\alpha
\restriction (A\cap \ell) = \rho^\ell_{\rho^n_\alpha(\ell)}$ for
$\ell\in A\cap n$. There are no problems and $\{\nu_\eta: \eta\in
\Gamma_n\}$ is as required.
\hfill$\square_{\ref{4.3}}$

\begin{Theorem}[MA]
\label{4.1}
Let $\lambda<2^{\aleph_0}$. Any $G\in {\mathfrak K}^{rs(p)}_\lambda$ can be
embedded into some $G'\in {\mathfrak K}^{rs(p)}_\lambda$ with countable
density which is tree-like.
\end{Theorem}

\proof By \ref{4.3} it suffices to get $G'$ ``almost tree-like" and
$A\subseteq\omega$ which satisfies \ref{4.3}(1). The ability to make $A$
thin helps in proving Fact E below. By \ref{1.1} without loss of generality
$G$ has a base (i.e. a dense subgroup of the form)
$B=\bigoplus\limits_{\scriptstyle n<\omega\atop\scriptstyle i<\lambda_n}
{\mathbb Z} x^n_i$, where ${\mathbb Z} x^n_i\cong{\mathbb Z}/p^{n+1}{\mathbb Z}$ and
$\lambda_n=\aleph_0$ (in fact $\lambda_n$ can be $g(n)$ if $g\in
{}^\omega\omega$ is not bounded (by algebraic manipulations), this will be
useful if we consider the forcing from \cite[\S2]{Sh:326}).

Let $B^+$ be the extension of $B$ by $y^{n,k}_i$ ($k<\omega$, $n<\omega$,
$i<\lambda_n$) generated freely except for $py^{n,k+1}_i=y^{n,k}_i$ (for
$k<\omega$), $y^{n,\ell}_i=p^{n-\ell}x^n_i$ for $\ell\le n$, $n<\omega$,
$i<\lambda_n$. So $B^+$ is a divisible $p$-group, let $G^+ =:
B^+\bigoplus\limits_B
G$. Let $\{z^0_\alpha:\alpha<\lambda\}\subseteq G[p]$ be a basis of $G[p]$
over $\{p^n x^n_i:n,i<\omega\}$ (as a vector space over ${\mathbb Z}/p{\mathbb
Z}$ i.e. the two sets are disjoint, their union is a basis); remember
$G[p]=\{x\in G:px=0\}$. So we can find $z^k_\alpha\in G$ (for $\alpha<
\lambda$, $k<\omega$ and $k\ne 0$) such that
\[pz^{k+1}_\alpha-z^k_\alpha=\sum_{i\in w(\alpha,k)} a^{k,\alpha}_i
x^k_i,\]
where $w(\alpha,k)\subseteq\omega$ is finite (reflect on the Abelian group
theory).

We define a forcing notion $P$ as follows: a condition $p \in P$ consists
of (in brackets are explanations of intentions):
\begin{description}
\item[(a)]  $m<\omega$, $M\subseteq m$,
\end{description}
[$M$ is intended as $A\cap\{0,\ldots,m-1\}$]
\begin{description}
\item[(b)]  a finite $u\subseteq m\times\omega$ and $h:u\longrightarrow
\omega$ such that $h(n,i)\ge n$,
\end{description}
[our extensions will not be pure, but still we want that the group produced
will be reduced, now we add some $y^{n,k}_i$'s and $h$ tells us how many]
\begin{description}
\item[(c)]  a subgroup $K$ of $B^+$:
\[K=\langle y^{n,k}_i:(n,i)\in u,k<h(n,i)\rangle_{B^+},\]
\item[(d)]  a finite $w\subseteq\lambda$,
\end{description}
[$w$ is the set of $\alpha<\lambda$ on which we give information]
\begin{description}
\item[(e)]  $g:w\rightarrow m + 1$,
\end{description}
[$g(\alpha)$ is in what level $m'\le m$ we ``start to think" about $\alpha$]
\begin{description}
\item[(f)]  $\bar\eta=\langle\eta_\alpha:\alpha\in w\rangle$ (see (i)),
\end{description}
[of course, $\eta_\alpha$ is the intended $\eta_\alpha$ restricted to $m$ and
the set of all $\eta_\alpha$ forms the intended $\Gamma$]
\begin{description}
\item[(g)]  a finite $v\subseteq m\times\omega$,
\end{description}
[this approximates the set of indices of the new basis]
\begin{description}
\item[(h)]  $\bar t=\{t_{n,i}:(n,i)\in v\}$ (see (j)),
\end{description}
[approximates the new basis]
\begin{description}
\item[(i)]  $\eta_\alpha\in {}^M\omega$, $\bigwedge\limits_{\alpha\in w}
\bigwedge\limits_{n\in M} (n,\eta_\alpha(n))\in v$,
\end{description}
[toward guaranteeing clause $(\delta)$ of \ref{4.2}(3) (see \ref{4.3}(2))]
\begin{description}
\item[(j)]  $t_{n,i}\in K$ and ${\mathbb Z} t_{n,i} \cong {\mathbb Z}/p^n
{\mathbb Z}$,
\item[(k)]  $K=\bigoplus\limits_{(n,i)\in v} ({\mathbb Z} t_{n,i})$,
\end{description}
[so $K$ is an approximation to the new basic subgroup]
\begin{description}
\item[(l)]  if $\alpha\in w$, $g(\alpha)\le\ell\le m$ and $\ell\in M$ then
\[z^\ell_\alpha-\sum\{t^{n-\ell}_{n,\eta_\alpha(n)}:\ell\le n\in
\Dom(\eta_\alpha)\}\in p^{m-\ell}(K+G),\]
\end{description}
[this is a step toward guaranteeing that the full difference (when
$\Dom(\eta_\alpha)$ is possibly infinite) will be in the closure of
$\bigoplus\limits_{\scriptstyle n\in [i,\omega)\atop\scriptstyle
i<\omega} {\mathbb Z} x^n_i$].

We define the order by:

\noindent $p \le q$ \qquad if and only if
\begin{description}
\item[$(\alpha)$]      $m^p\le m^q$, $M^q \cap m^p = M^p$,
\item[$(\beta)$]       $u^p\subseteq u^q$, $h^p\subseteq h^q$,
\item[$(\gamma)$]      $K^p\subseteq_{pr} K^q$,
\item[$(\delta)$]      $w^p\subseteq w^q$,
\item[$(\varepsilon)$] $g^p\subseteq g^q$,
\item[$(\zeta)$]       $\eta^p_\alpha\trianglelefteq\eta^q_\alpha$,
(i.e. $\eta^p_\alpha$ is an initial segment of $\eta^q_\alpha$)
\item[$(\eta)$]        $v^p\subseteq v^q$,
\item[$(\theta)$]      $t^p_{n,i}=t^q_{n,i}$ for $(n,i)\in v^p$.
\end{description}
\medskip

\noindent{\bf A Fact}\hspace{0.15in}  $(P,\le)$ is a partial order.
\medskip

\noindent{\em Proof of the Fact:}\ \ \ Trivial.
\medskip

\noindent{\bf B Fact}\hspace{0.15in}  $P$ satisfies the c.c.c. (even is
$\sigma$-centered).
\medskip

\noindent{\em Proof of the Fact:}\ \ \ It suffices to observe the following.

Suppose that
\begin{description}
\item[$(*)$(i)]    $p,q \in P$,
\item[\quad(ii)]   $M^p=M^q$, $m^p=m^q$, $h^p=h^q$, $u^p=u^q$, $K^p=K^q$,
$v^p=v^q$, $t^p_{n,i}=t^q_{n,i}$,
\item[\quad(iii)]  $\langle\eta^p_\alpha:\alpha\in w^p\cap w^q\rangle =
\langle\eta^q_\alpha:\alpha\in w^p\cap w^q\rangle$,
\item[\quad(iv)]   $g^p\restriction (w^p \cap w^q)=g^q \restriction(w^p
\cap w^q)$.
\end{description}
Then the conditions $p,q$ are compatible (in fact have an upper bound with
the same common parts): take the common values (in (ii)) or the union (for
(iii)).
\medskip

\noindent{\bf C Fact}\hspace{0.15in}  For each $\alpha<\lambda$ the set
${\Cal I}_\alpha=:\{p\in P:\alpha\in w^p\}$ is dense (and open).
\medskip

\noindent{\em Proof of the Fact:}\ \ \ For $p\in P$ let $q$ be like $p$
except that:
\[w^q=w^p\cup\{\alpha\}\quad\mbox{ and }\quad g^q(\beta)=\left\{
\begin{array}{lll}
g^p(\beta) &\mbox{if}& \beta\in w^p \\
m^p        &\mbox{if}& \beta=\alpha,\ \beta\notin w^p.
\end{array}\right.\]
\medskip

\noindent{\bf D Fact}\hspace{0.15in}  For $n<\omega$, $i<\omega$ the
following set is a dense subset of $P$:
\[{\Cal J}^*_{(n,i)}=\{p\in P:x^n_i\in K^p\ \&\ (\forall n<m^p)(\{n\}\times
\omega)\cap u^p \mbox{ has }>m^p\mbox{ elements}\}.\]
\medskip

\noindent{\em Proof of the Fact:}\ \ \ Should be clear.
\medskip

\noindent{\bf E Fact}\hspace{0.15in} For each $m<\omega$ the set ${\Cal J}_m
=:\{p\in P:m^p\ge m\}$ is dense in $P$.
\medskip

\noindent{\em Proof of the Fact:}\ \ \ Let $p\in P$ be given such that $m^p
<m$. Let $w^p=\{\alpha_0,\ldots,\alpha_{r-1}\}$ be without repetitions;
we know that in $G$, $pz^0_{\alpha_\ell}=0$ and $\{z^0_{\alpha_\ell}:
\ell<r\}$ is independent $\mod\ B$, hence also in $K+G$ the set
$\{z^0_{\alpha_\ell}:\ell<r\}$ is independent $\mod\ K$. Clearly
\begin{description}
\item[(A)] $pz^{k+1}_{\alpha_\ell}=z^k_{\alpha_\ell}\mod\ K$
for $k\in [g(\alpha_\ell),m^p)$, hence
\item[(B)] $p^{m^p}z^{m^p}_{\alpha_\ell}=z^{g(\alpha_\ell)}_{\alpha_\ell}
\mod\ K$.
\end{description}
Remember
\begin{description}
\item[(C)]  $z^{m^p}_{\alpha_\ell}=\sum\{a^{k,\alpha_\ell}_i p^{k-m^p} x^k_i:
k \ge m^p,i\in w(\alpha_\ell,k)\}$,
\end{description}
and so, in particular, (from the choice of $z^0_{\alpha_\ell}$)
\[p^{m^p+1}z^{m^p}_{\alpha_\ell}=0\quad\mbox{ and }\quad
p^{m^p}z^{m^p}_{\alpha_\ell}\ne 0.\]
For $\ell<r$ and $n\in [m^p,\omega)$ let
\[s^n_\ell=:\sum\big\{a^{k,\alpha_\ell}_i p^{k-m^p} x^k_i:k \ge m^p
\mbox{ but }k<n\mbox{ and } i\in w(\alpha_\ell,k)\big\}.\]
But $p^{k-m^p}x^k_i = y^{k,m^p}_i$, so
\[s^n_\ell=\sum\big\{a^{k,\alpha_\ell}_i y^{k,m^p}_i:k\in [m^p,n)
\mbox{ and }i\in (\alpha_\ell,k)\big\}.\]
Hence, for some $m^*>m,m^p$ we have: $\{p^m\,s^{m^*}_\ell:\ell<r\}$ is
independent in $G[p]$ over $K[p]$ and therefore in $\langle x^k_i:k\in
[m^p,m^*],i<\omega\rangle$. Let
\[s^*_\ell=\sum\big\{a^{k,\alpha_\ell}_i:k\in [m^p,m^*)\mbox{ and }i\in
w(\alpha_\ell,k)\}.\]
Then $\{ s^*_\ell:\ell<r\}$ is independent in
\[B^+_{[m,m^*)}=\langle y^{l,m^*-1}_i:k\in [m^p,m^*)\mbox{ and }i<\omega
\rangle.\]
Let $i^*<\omega$ be such that: $w(\alpha_\ell,k)\subseteq\{0,\ldots,i^*-1\}$
for $k\in [m^p,m^*)$, $\ell=1,\ldots,r$. Let us start to define $q$:
\[\begin{array}{c}
m^q=m^*,\quad M^q=M^p\cup\{m^*-1\},\quad w^q=w^p,\quad g^q=g^p,\\
u^q=u^p\cup ([m^p,m^*)\times\{0,\ldots,i^*-1\}),\\
h^q\mbox{ is } h^p\mbox{ on } u^p\mbox{ and }h^q(k,i)=m^*-1\mbox{
otherwise},\\
K^q\mbox{ is defined appropriately, let } K'=\langle x^n_i:n\in [m^p,m^*),
i<i^*\rangle.
\end{array}\]
Complete $\{s^*_\ell:\ell<r\}$ to $\{s^*_\ell:\ell<r^*\}$, a basis of $K'[p]$,
and choose $\{t_{n,i}:(n,i)\in v^*\}$ such that: $[p^mt_{n,i}=0\ \
\Leftrightarrow\ \ m>n]$, and for $\ell<r$
\[p^{m^*-1-\ell}t_{m^*-1,\ell} = s^*_\ell.\]
The rest should be clear.
\medskip

The generic gives a variant of the desired result: almost tree-like basis; the
restriction to $M$ and $g$ but by \ref{4.3} we can finish.
\hfill$\square_{\ref{4.2}}$

\begin{Conclusion}
[MA$_\lambda$($\sigma$-centered)]
\label{4.4}
For $(*)_0$ to hold it suffices that $(*)_1$ holds where
\begin{description}
\item[$(*)_0$]  in ${\mathfrak K}^{rs(p)}_\lambda$, there is a universal
member,
\item[$(*)_1$]  in ${\mathfrak K}^{tr}_{\bar\lambda}$ there is a universal
member, where:
\begin{description}
\item[(a)] $\lambda_n=\aleph_0$, $\lambda_\omega=\lambda$, $\ell g(\bar\lambda)
=\omega+1$\qquad \underbar{or}
\item[(b)] $\lambda_\omega=\lambda$, $\lambda_n\in [n,\omega)$, $\ell g
(\bar \lambda)=\omega+1$.
\end{description}
\end{description}
\end{Conclusion}

\begin{Remark}
\label{4.4A}
Any $\langle\lambda_n:n<\omega\rangle$, $\lambda_n<\omega$ which is not
bounded suffices.
\end{Remark}

\proof For case (a) - by \ref{4.1}.

\noindent For case (b) - the same proof. \hfill$\square_{\ref{4.4}}$

\begin{Theorem}
\label{4.5}
Assume $\lambda<2^{\aleph_0}$ and
\begin{description}
\item[(a)]  there are $A_i\subseteq\lambda$, $|A_i|=\lambda$ for
$i<2^\lambda$ such that $i\ne j \Rightarrow |A_i \cap A_j| \le \aleph_0$.
\end{description}
Let $\bar\lambda=\langle \lambda_\alpha:\alpha\le\omega\rangle$, $\lambda_n
= \aleph_0$, $\lambda_\omega=\lambda$.

\noindent{\em Then} there is $P$ such that:
\medskip
\begin{description}
\item[$(\alpha)$]  $P$ is a c.c.c. forcing notion,
\item[$(\beta)$]   $|P|=2^\lambda$,
\item[$(\gamma)$]  in $V^P$, there is $T\in {\mathfrak K}^{tr}_{\bar\lambda}$
into which every $T' \in ({\mathfrak K}^{tr}_{\bar \lambda})^V$ can be embedded.
\end{description}
\end{Theorem}

\proof  Let $\bar T=\langle T_i:i<2^\lambda\rangle$ list the trees $T$ of
cardinality $\le\lambda$ satisfying
\[{}^{\omega >}\omega\subseteq T \subseteq {}^{\omega \ge} \omega\quad
\mbox{ and }\quad T\cap {}^\omega\omega\mbox{ has cardinality $\lambda$, for
simplicity.}\]
Let $T_i\cap {}^\omega\omega=\{\eta^i_\alpha:\alpha\in A_i \}$.

We shall force $\rho_{\alpha,\ell}\in {}^\omega\omega$ for $\alpha<
\lambda$, $\ell<\omega$, and for each $i<2^\lambda$ a function $g_i:A_i
\longrightarrow\omega$ such that: there is an automorphism $f_i$ of
$({}^{\omega>}\omega,\triangleleft)$ which induces an embedding of $T_i$
into $\left(({}^{\omega>}\omega)\cup \{\rho_{\alpha,g_i(\alpha)}:\alpha<
\lambda\},\triangleleft\right)$.  We shall define $p\in P$ as an approximation.

\noindent A condition $p\in P$ consists of:
\begin{description}
\item[(a)]    $m<\omega$ and a finite subset $u$ of ${}^{m \ge}\omega$, closed
under initial segments such that $\langle\rangle\in u$,
\item[(b)]    a finite $w\subseteq 2^\lambda$,
\item[(c)]    for each $i\in w$, a finite function $g_i$ from $A_i$ to
$\omega$,
\item[(d)]    for each $i\in w$, an automorphism $f_i$ of
$(u,\triangleleft)$,
\item[(e)]    a finite $v\subseteq\lambda\times\omega$,
\item[(f)]    for $(\alpha,n)\in v$, $\rho_{\alpha,n}\in u\cap
({}^m\omega)$,
\end{description}
such that
\begin{description}
\item[(g)]  if $i\in w$ and $\alpha\in\Dom(g_i)$ then:
\begin{description}
\item[$(\alpha)$]   $(\alpha,g_i(\alpha))\in v$,
\item[$(\beta)$]    $\eta^i_\alpha\restriction m\in u$,
\item[$(\gamma)$]   $f_i(\eta^i_\alpha\restriction m)=\rho_{\alpha,
g_i(\alpha)}$,
\end{description}
\item[(h)]   $\langle\rho_{\alpha,n}:(\alpha,n)\in v\rangle$ is with no
repetition (all of length $m$),
\item[(i)]   for $i\in w$, $\langle\eta^i_\alpha\restriction m:\alpha\in
\Dom(g_i)\rangle$ is with no repetition.
\end{description}

The order on $P$ is: $p \le q$ if and only if:
\begin{description}
\item[$(\alpha)$]   $u^p \subseteq u^q$, $m^p\le m^q$,
\item[$(\beta)$]    $w^p \subseteq w^q$,
\item[$(\gamma)$]   $f^p_i \subseteq f^q_i$ for $i\in w^p$,
\item[$(\delta)$]   $g^p_i \subseteq g^q_i$ for $i\in w^p$,
\item[$(\varepsilon)$]  $v^p \subseteq v^q$,
\item[$(\zeta)$]    $\rho^p_{\alpha,n}\trianglelefteq\rho^q_{\alpha,n}$,
when $(\alpha,n) \in v^p$,
\item[$(\eta)$]     if $i\ne j\in w^p$ then for every $\alpha\in A_i\cap
A_j\setminus (\Dom(g^p_i)\cap \Dom(g^p_j))$ we have $g^q_i(\alpha)\ne
g^q_j(\alpha)$.
\end{description}
\medskip

\noindent{\bf A Fact}\hspace{0.15in}  $(P,\le)$ is a partial order.
\medskip

\noindent{\em Proof of the Fact:}\ \ \ Trivial.
\medskip

\noindent{\bf B Fact}\hspace{0.15in}  For $i<2^\lambda$ the set $\{p:i\in
w^p\}$ is dense in $P$.
\medskip

\noindent{\em Proof of the Fact:}\ \ \ If $p\in P$, $i\in 2^\lambda
\setminus w^p$, define $q$ like $p$ except $w^q=w^p\cup\{i\}$,
$\Dom(g^q_i)=\emptyset$.
\medskip

\noindent{\bf C Fact}\hspace{0.15in}  If $p\in P$, $m_1\in
(m^p,\omega)$, $\eta^*\in u^p$, $m^*<\omega$, $i\in w^p$, $\alpha\in\lambda
\setminus\Dom(g^p_i)$ {\em then} we can find $q$ such that $p\le q\in
P$, $m^q>m_1$, $\eta^* \char 94\langle m^*\rangle\in u^q$ and $\alpha\in
\Dom(g_i)$ and $\langle\eta^j_\beta\restriction m^q:j\in w^q$ and $\beta\in
\Dom(g^q_j)\rangle$ is with no repetition, more exactly
$\eta^{j(1)}_{\beta_1}\setminus m^q= \eta^{j(2)}_{\beta_2}\restriction
m^q \Rightarrow \eta^{j(1)}_{\beta_1}=\eta^{j(2)}_{\beta_2}$.
\medskip

\noindent{\em Proof of the Fact:}\ \ \ Let $n_0\le m^p$ be maximal such that
$\eta^i_\alpha\restriction n_0 \in u^p$. Let $n_1<\omega$ be minimal such
that $\eta^i_\alpha\restriction n_1\notin\{\eta^i_\beta\restriction n_1:
\beta\in\Dom(g^p_i)\}$ and moreover the sequence
\[\langle\eta^j_\beta\restriction n_1:j\in w^p\ \&\ \beta\in\Dom(g^p_j)\
\ \mbox{ or }\ \ j=i\ \&\ \beta=\alpha\rangle\]
is with no repetition. Choose a natural number  $m^q>m^p+1,n_0+1,n_1+2$ and
let $k^*=:3+\sum\limits_{i\in w^p}|\Dom(g^p_i)|$. Choose $u^q\subseteq
{}^{m^q\ge}\omega$ such that:
\begin{description}
\item[(i)]    $u^p\subseteq u^q\subseteq {}^{m^q\ge}\omega$, $u^q$ is downward
closed,
\item[(ii)]   for every $\eta\in u^q$ such that $\ell g(\eta)<m^q$, for
exactly $k^*$ numbers $k$, $\eta\char 94\langle k\rangle\in u^q\setminus
u^p$,
\item[(iii)]  $\eta^j_\beta\restriction\ell\in u^q$ when $\ell\le m^q$
and $j\in w^p$, $\beta\in\Dom(g^p_j)$,
\item[(iv)]   $\eta^i_\alpha\restriction\ell\in u^q$ for $\ell\le m^q$,
\item[(v)]    $\eta^*\char 94\langle m^* \rangle\in u^q$.
\end{description}
Next choose $\rho^q_{\beta,n}$ (for pairs $(\beta,n)\in v^p)$ such that:
\[\rho^p_{\beta,n}\trianglelefteq\rho^q_{\beta,n}\in u^q\cap {}^{m^q}
\omega.\]
For each $j\in w^p$ separately extend $f^p_j$ to an automorphism $f^q_j$
of $(u^q,\triangleleft)$ such that for each $\beta\in\Dom(g^p_j)$ we have:
\[f^q_j(\eta^j_\beta\restriction m^q)=\rho^q_{\beta,g_j}(\beta).\]
This is possible, as for each $\nu\in u^p$, and $j\in w^p$, we can
separately define
\[f^q_j\restriction\{\nu':\nu\triangleleft\nu'\in u^q\ \mbox{ and }\ \nu'
\restriction (\ell g(\nu)+1)\notin u^p\}\]
--its range is
\[\{\nu':f^p_j(\nu)\triangleleft \nu'\in u^q\ \mbox{ and }\ \nu'
\restriction (\ell g(\nu)+1)\notin u^p\}.\]
The point is: by Clause (ii) above those two sets are isomorphic and for
each $\nu$ at most one $\rho^p_{\beta,n}$ is involved (see Clause (h) in the
definition of $p \in P$). Next let $w^q=w^p$, $g^q_j=g^p_j$ for $j\in w
\setminus\{i\}$, $g^q_i\restriction\Dom(g^p_i)=g^p_i$, $g^q_i(\alpha)=
\min(\{n:(\alpha,n)\notin v^p\})$, $\Dom(g^q_i)=\Dom(g^p_i)\cup\{\alpha\}$,
and $\rho^q_{\alpha,g^q_i(\alpha)}=f^g_i(\eta^i_\alpha\restriction m^q)$ and
$v^q=v^p\cup\{(\alpha,g^q_i(\alpha))\}$.
\medskip

\noindent{\bf D Fact}\hspace{0.15in}  $P$ satisfies the c.c.c.
\medskip

\noindent{\em Proof of the Fact:}\ \ \ Assume $p_\varepsilon\in P$ for
$\varepsilon<\omega_1$. By Fact C, without loss of generality each
\[\langle\eta^j_\beta\restriction m^{p_\varepsilon}:j\in w^{p_\varepsilon}
\mbox{ and }\beta\in\Dom(g^{p_\varepsilon}_j)\rangle\]
is with no repetition. Without loss of generality, for all $\varepsilon
<\omega_1$
\[U_\varepsilon=:\big\{\alpha<2^\lambda:\alpha\in w^{p_\varepsilon}\mbox{ or }
\bigvee_{i\in w^p}[\alpha\in\Dom(g_i)]\mbox{ or }\bigvee_k(k,\alpha)\in
v^{p_\varepsilon}\big\}\]
has the same number of elements and for $\varepsilon\ne\zeta<\omega_1$,
there is a unique one-to-one order preserving function from $U_\varepsilon$
onto $U_\zeta$ which we call $\OP_{\zeta,\varepsilon}$, which also maps
$p_\varepsilon$ to $p_\zeta$ (so $m^{p_\zeta}=m^{p_\varepsilon}$; $u^{p_\zeta}
=u^{p_\varepsilon}$; $\OP_{\zeta,\varepsilon}(w^{p_\varepsilon})=
w^{p_\zeta}$; if $i\in w^{p_\varepsilon}$, $j=\OP_{\zeta,\varepsilon}(i)$,
then $f_i\circ\OP_{\varepsilon,\zeta}\equiv f_j$; and {\em if\/} $\beta=
\OP_{\zeta,\varepsilon}(\alpha)$ and $\ell<\omega$ {\em then}
\[(\alpha,\ell)\in v^{p_\varepsilon}\quad\Leftrightarrow\quad (\beta,\ell)
\in v^{p_\zeta}\quad\Rightarrow\quad\rho^{p_\varepsilon}_{\alpha,\ell}=
\rho^{p_\zeta}_{\beta,\ell}).\]
Also this mapping is the identity on $U_\zeta\cap U_\varepsilon$ and
$\langle U_\zeta:\zeta<\omega_1\rangle$ is a $\triangle$-system.

Let $w=:w^{p_0}\cap w^{p_1}$.  As $i\ne j\ \Rightarrow\ |A_i\cap A_j|\le
\aleph_0$, without loss of generality
\begin{description}
\item[$(*)$]  if $i\ne j\in w$ then
\[U_{\varepsilon}\cap (A_i\cap A_j)\subseteq w.\]
\end{description}
We now start to define $q\ge p_0,p_1$. Choose $m^q$ such that $m^q\in
(m^{p_\varepsilon},\omega)$ and
\[\begin{array}{ll}
m^q>\max\big\{\ell g(\eta^{i_0}_{\alpha_0}\cap\eta^{i_1}_{\alpha_1})+1:&
i_0\in w^{p_0},\ i_1 \in w^{p_1},\ \OP_{1,0}(i_0)=i_1,\\
\ &\alpha_0\in\Dom(g^{p_0}_{i_0}),\ \alpha_1\in\Dom(g^{p_1}_{i_1}),\\
\ &\OP_{1,0}(\alpha_0)=\alpha_1\big\}.
\end{array}\]
Let $u^q\subseteq {}^{m^q\ge}\omega$ be such that:
\begin{description}
\item[(A)]  $u^q\cap\left({}^{m^{p_0}\ge}\omega\right)=u^q\cap\left(
{}^{m^{p_1}\ge}\omega\right)=u^{p_0}=u^{p_1}$,
\item[(B)]  for each $\nu\in u^q$, $m^{p_0}\le\ell g(\nu)<m^q$, for exactly
two numbers $k<\omega$, $\nu\char 94 \langle k\rangle\in u^q$,
\item[(C)]  $\eta^i_\alpha\restriction\ell\in u^q$ for $\ell\le m^q$
\underbar{when}: $i\in w^{p_0}$, $\alpha\in\Dom(g^{p_0}_i)$ \underbar{or}
$i\in w^{p_1}$, $\alpha\in\Dom(g^{p_1}_i)$.
\end{description}
[Possible as $\{\eta^i_\alpha\restriction m^{p_\varepsilon}:i\in
w^{p_\varepsilon},\alpha\in\Dom(g^{p_\varepsilon}_i)\}$ is with no
repetitions (the first line of the proof).]

Let $w^q=:w^{p_0}\cup w^{p_1}$ and $v^q=:v^{p_0}\cup v^{p_1}$ and for
$i \in w^q$
\[g^q_i=\left\{\begin{array}{lll}
g^{p_0}_i &\mbox{\underbar{if}}& i\in w^{p_0}\setminus w^{p_1},\\
g^{p_1}_i &\mbox{\underbar{if}}& i\in w^{p_1}\setminus w^{p_0},\\
g^{p_0}_i \cup g^{p_1}_i &\mbox{\underbar{if}}& i\in w^{p_0}\cap w^{p_1}.
\end{array}\right.\]
Next choose $\rho^q_{\alpha,\ell}$ for $(\alpha,\ell)\in v^q$ as follows.
Let $\nu_{\alpha,\ell}$ be $\rho^{p_0}_{\alpha,\ell}$ if defined,
$\rho^{p_1}_{\alpha,\ell}$ if defined (no contradiction). If $(\alpha,\ell)
\in v^q$ choose $\rho^q_{\alpha,\ell}$ as any $\rho$ such that:
\begin{description}
\item[$\otimes_0$]  $\nu_{\alpha,\ell}\triangleleft\rho\in u^q\cap
{}^{(m^q)}\omega$.
\end{description}
But not all choices are O.K., as we need to be able to define $f^q_i$ for
$i\in w^q$. A possible problem will arise only when $i\in w^{p_0}\cap
w^{p_1}$.  Specifically we need just (remember that $\langle
\rho^{p_\varepsilon}_{\alpha,\ell}:(\alpha,\ell)\in v^{p_\varepsilon}
\rangle$ are pairwise distinct by clause (b) of the Definition of $p\in P$):
\begin{description}
\item[$\otimes_1$]  if $i_0\in w^{p_0},(\alpha_0,\ell)=(\alpha_0,
g_{i_0}(\alpha_0)),\alpha_0\in\Dom(g^{p_0}_{i_0})$, $i_1=\OP_{1,0}(i_0)$
and $\alpha_1=\OP_{1,0}(\alpha_0)$ and $i_0=i_1$

{\em then} $\ell g(\eta^{i_0}_{\alpha_0}\cap\eta^{i_1}_{\alpha_1})=\ell g
(\rho^q_{\alpha_0,\ell}\cap\rho^q_{\alpha_1,\ell})$.
\end{description}
We can, of course, demand $\alpha_0\neq \alpha_1$ (otherwise the
conclusion of $\otimes_1$ is trivial).
Our problem is expressible for each pair $(\alpha_0,\ell),(\alpha_1,\ell)$
separately as: first the problem is in defining the
$\rho^q_{(\alpha,\ell)}$'s and second, if $(\alpha'_1,\ell')$,
$(\alpha'_2,\ell)$ is another such pair then $\{(\alpha_1,\ell),
(\alpha_2,\ell)\}$, $\{(\alpha'_1,\ell'),(\alpha'_2,\ell')\}$ are either
disjoint or equal. Now for a given pair $(\alpha_0,\ell),(\alpha_1,\ell)$ how
many $i_0=i_1$ do we have?  Necessarily $i_0\in w^{p_0}\cap w^{p_1}=w$. But
if $i'_0\ne i''_0$ are like that then $\alpha_0\in A_{i'_0}\cap
A_{i''_0}$, contradicting $(*)$ above because $\alpha_0\neq
\alpha_1=\OP_{1,0}(\alpha_0)$. So there is at most one candidate
$i_0=i_1$, so there is no problem to satisfy $\otimes_1$. Now we can define
$f^q_i$ (i$\in w^q$) as in the proof of Fact C.
\medskip

The rest should be clear. \hfill$\square_{\ref{4.4}}$

\begin{Conclusion}
\label{4.6}
Suppose $V\models GCH$, $\aleph_0<\lambda<\chi$ and $\chi^\lambda=\chi$.
Then for some c.c.c. forcing notion $P$ of cardinality $\chi$, not collapsing
cardinals nor changing cofinalities, in $V^P$:
\begin{description}
\item[(i)]   $2^{\aleph_0}=2^\lambda=\chi$,
\item[(ii)]  ${\mathfrak K}^{tr}_\lambda$ has a universal family of cardinality
$\lambda^+$,
\item[(iii)] ${\mathfrak K}^{rs(p)}_\lambda$ has a universal family of cardinality
$\lambda^+$.
\end{description}
\end{Conclusion}

\proof First use a preliminary forcing $Q^0$ of Baumgartner \cite{B},
adding $\langle A_\alpha:\alpha<\chi\rangle$, $A_\alpha\in
[\lambda]^\lambda$, $\alpha\neq\beta\quad \Rightarrow\quad |A_\alpha\cap
A_\beta|\le\aleph_0$ (we can have $2^{\aleph_0}=\aleph_1$ here, or
$[\alpha\ne \beta\quad \Rightarrow\quad A_\alpha\cap A_\beta$ finite], but
not both).  Next use an FS iteration $\langle P_i,\dot{Q}_i:i<\chi\times
\lambda^+\rangle$ such that each forcing from \ref{4.4} appears and each
forcing as in \ref{4.5} appears. \hfill$\square_{\ref{4.6}}$

\begin{Remark}
\label{4.7}
We would like to have that there is a universal member in
${\mathfrak K}^{rs(p)}_\lambda$; this sounds very reasonable but we did not try.

In our framework, the present result shows limitations to ZFC results which
the methods applied in the previous sections can give.
\end{Remark}

\section{Back to ${\mathfrak K}^{rs(p)}$, real non-existence results}
By \S1 we know that if $G$ is an Abelian group with set of elements
$\lambda$, $C\subseteq\lambda$, then for an element $x\in G$ the distance
from $\{y:y<\alpha\}$ for $\alpha\in C$ does not code an appropriate invariant.
If we have infinitely many such distance functions, e.g. have infinitely many
primes, we can use more complicated invariants related to $x$ as in \S3.
But if we have one prime, this approach does not help.

If one element fails, can we use infinitely many?  A countable subset $X$ of
$G$ can code a countable subset of $C$:
\[\{\alpha\in C:\mbox{ closure}(\langle X\rangle_G)\cap\alpha\nsubseteq
\sup(C\cap\alpha)\},\]
but this seems silly - we use heavily the fact that $C$ has many countable
subsets (in particular $>\lambda$) and $\lambda$ has at least as many.
However, what if $C$ has a small family (say of cardinality $\le\lambda$ or
$<\mu^{\aleph_0}$) of countable subsets such that every subset of cardinality,
say continuum, contains one? Well, we need more: we catch a countable subset
for which the invariant defined above is infinite (necessarily it is at most
of cardinality $2^{\aleph_0}$, and because of \S4 we are not trying any more
to deal with $\lambda\le 2^{\aleph_0}$). The set theory needed is expressed
by $T_J$ below, and various ideals also defined below, and the result itself
is \ref{5.7}.

Of course, we can deal with other classes like torsion free reduced groups,
as they  have the characteristic non-structure property of unsuperstable
first order theories; but the relevant ideals will vary: the parallel to
$I^0_{\bar\mu}$ for $\bigwedge\limits_n \mu_n=\mu$, $J^2_{\bar\mu}$ seems
to be always O.K.

\begin{Definition}
\label{5.1}
\begin{enumerate}
\item For $\bar\mu=\langle\mu_n:n<\omega\rangle$ let $B_{\bar\mu}$ be
\[\bigoplus\{K^n_\alpha:n<\omega,\alpha<\mu_n\},\qquad K^n_\alpha=
\langle{}^* t^n_\alpha\rangle_{K^n_\alpha}\cong {\mathbb Z}/p^{n+1} {\mathbb Z}.\]
Let $B_{\bar\mu\restriction n}=\bigoplus\{K^m_\alpha:\alpha<\mu_m,m<n\}
\subseteq B_{\bar\mu}$ (they are in ${\mathfrak K}^{rs(p)}_{\le\sum\limits_{n}
\mu_n}$). Let $\hat B$ be the $p$-torsion completion of $B$ (i.e. completion
under the norm $\|x\|=\min\{2^{-n}:p^n\mbox{ divides }x\}$).
\item  Let $I^1_{\bar \mu}$ be the ideal on $\hat B_{\bar \mu}$
generated by $I^0_{\bar \mu}$, where
\[\begin{array}{ll}
I^0_{\bar\mu}=\big\{A\subseteq\hat B_{\bar\mu}:&\mbox{for every large
enough }n,\\
 &\mbox{for no }y\in\bigoplus\{K^m_\alpha:m\le n\mbox{ and }\alpha<\mu_m\}\\
 &\mbox{but }y\notin\bigoplus\{K^m_\alpha:m<n\mbox{ and }\alpha<\mu_m\}
  \mbox{ we have}:\\
 &\mbox{for every }m\mbox{ for some }z\in\langle A\rangle\mbox{ we have:}\\
 &p^m\mbox{ divides }z-y\big\}.
\end{array}\]
(We may write $I^0_{\hat B_{\bar\mu}}$, but the ideal depends also on
$\langle\bigoplus\limits_{\alpha<\mu_n} K^n_\alpha:n<\omega\rangle$ not
just on $\hat B_{\bar\mu}$ itself).

\item  For $X,A\subseteq\hat B_{\bar\mu}$,
\[\mbox{ recall }\ \ \langle A\rangle_{\bar
B_{\bar\mu}}=\big\{\sum\limits_{n<n^*} a_ny_n:y_n
\in A,\ a_n\in{\mathbb Z}\mbox{ and } n^*\in{\mathbb N}\big\},\]
\[\mbox{ and let }\ \ c\ell_{\hat B_{\bar\mu}}(X)=\{x:(\forall
n)(\exists y\in X)(x-y\in p^n \hat B_{\bar \mu})\}.\]

\item Let $J^1_{\bar \mu}$ be the ideal which $J^{0.5}_{\bar \mu}$
generates, where
\[\begin{array}{ll}
J^{0.5}_{\bar\mu}=\big\{A\subseteq\prod\limits_{n<\omega}\mu_n:
 &\mbox{for some }n<\omega\mbox{ for no }m\in [n,\omega)\\
 &\mbox{and }\beta<\gamma<\mu_m\mbox{ do we have}:\\
 &\mbox{for every }k\in [m,\omega)\mbox{ there are }\eta,\nu\in A\mbox{
 such}\\
 &\mbox{that:}\ \ \eta(m)=\beta,\,\nu(m)=\gamma,\,\eta\restriction m=
\nu\restriction m\\
 &\mbox{and }\eta\restriction (m,k)=\nu\restriction (m,k)\big\}.
\end{array}\]

\item
\[\begin{array}{rr}
J^0_{\bar\mu}=\{A\subseteq\prod\limits_{n<\omega}\mu_n:
 &\mbox{for some }n<\omega\mbox{ and } k,\mbox{ the mapping }\eta\mapsto
\eta \restriction n\\
 &\mbox{is }(\le k)\mbox{-to-one }\}.
\end{array}\]

\item $J^2_{\bar\mu}$ is the ideal of nowhere dense subsets of
$\prod\limits_n\mu_n$ (under the following natural topology: a
neighbourhood of $\eta$ is $U_{\eta,n}=\{\nu:\nu\restriction n=\eta
\restriction n\}$ for some $n$).

\item  $J^3_{\bar \mu}$ is the ideal of meagre subsets of $\prod\limits_n
\mu_n$, i.e. subsets which are included in countable union of members of
$J^2_{\bar \mu}$.
\end{enumerate}
\end{Definition}

\begin{Observation}
\label{5.2}
\begin{enumerate}
\item  $I^0_{\bar\mu}$, $J^0_{\bar\mu}$, $J^{0.5}_{\bar\mu}$ are
$(<\aleph_1)$-based, i.e. for $I^0_{\bar \mu}$: if $A\subseteq\hat
B_{\bar\mu}$, $A\notin I^0_{\bar\mu}$ then there is a countable $A_0
\subseteq A$ such that $A_0\notin I^0_{\bar\mu}$.
\item  $I^1_{\bar\mu}$, $J^0_{\bar\mu}$, $J^1_{\bar\mu}$,
$J^2_{\bar\mu}$, $J^3_{\bar\mu}$ are ideals, $J^3_{\bar\mu}$ is
$\aleph_1$-complete.
\item  $J^0_{\bar\mu}\subseteq J^1_{\bar\mu}\subseteq J^2_{\bar\mu}
\subseteq J^3_{\bar\mu}$.
\item  There is a function $g$ from $\prod\limits_{n<\omega}\mu_n$ into
$\hat B_{\bar\mu}$ such that for every $X\subseteq\prod\limits_{n<\omega}
\mu_n$:
\[X\notin J^1_{\bar\mu}\quad \Rightarrow\quad g''(X)\notin
I^1_{\bar\mu}.\]
\end{enumerate}
\end{Observation}

\proof E.g. 4)\ \ Let $g(\eta)=\sum\limits_{n<\omega}p^n({}^*t^n_{\eta(n)})$.

Let $X \subseteq\prod\limits_{n<\omega}\mu_n$, $X\notin J^1_{\bar\mu}$.
Assume $g''(X)\in\bar I^1_{\bar\mu}$, so for some $\ell^*$ and $A_\ell
\subseteq\hat B_{\bar\mu}$, ($\ell<\ell^*$) we have $A_\ell\in I^0_{\bar\mu}$,
and $g''(X)\subseteq\bigcup\limits_{\ell<\ell^*} A_\ell$, so $X=
\bigcup\limits_{\ell<\ell^*} X_\ell$, where
\[X_\ell=:\{\eta\in X:g(\eta)\in A_\ell\}.\]
As $J^1_{\bar \mu}$ is an ideal, for some $\ell<\ell^*$, $X_\ell\notin
J^1_{\bar\mu}$. So by the definition of $J^1_{\bar\mu}$, for some infinite
$\Gamma\subseteq\omega$ for each $m\in\Gamma$ we have $\beta_m<\gamma_m<
\mu_m$ and for every $k\in [m,\omega)$ we have $\eta_{m,k},\nu_{m,k}$, as
required in the definition of $J^1_{\bar \mu}$.  So $g(\eta_{m,k}),
g(\nu_{m,k}) \in A_\ell$ (for $m\in \Gamma$, $k\in (m,\omega)$).  Now
\[{}^* t^m_{\gamma_m} - {}^* t^m_{\beta_m}=g(\eta_{m,k})-g(\nu_{m,k})
\mod\ p^k \hat B_{\bar \mu},\]
but $g(\eta_{m,k})-g(\nu_{m,k})\in\langle A_\ell\rangle_{\hat B_{\bar\mu}}$.
Hence
\[(\exists z\in\langle A_\ell\rangle_{\hat B_{\bar\mu}})[{}^* t^m_{\gamma_m}
-{}^* t^m_{\beta_m}=z\mod\ p^k \hat B_{\bar\mu}],\]
as this holds for each $k$, ${}^* t^m_{\gamma_m}-{}^* t^m_{\beta_m}\in
c \ell(\langle A_\ell\rangle_{\hat B_{\bar\mu}})$.

This contradicts $A_\ell
\in I^0_{\bar \mu}$. \hfill$\square_{\ref{5.2}}$

\begin{Definition}
\label{5.3}
Let $I\subseteq{\Cal P}(X)$ be downward closed (and for simplicity $\{\{x\}:
x\in X\}\subseteq I$). Let $I^+={\Cal P}(X)\setminus I$. Let
\[\begin{array}{ll}
\bU^{<\kappa}_I(\mu)=:\min\big\{|{\Cal P}|:&{\Cal P}\subseteq [\mu]^{<\kappa},
\mbox{ and for every } f:X\longrightarrow\mu\mbox{ for some}\\
  &Y\in {\Cal P},\mbox{ we have }\{x\in X:f(x)\in Y\}\in
I^+\big\}.
\end{array}\]
Instead of $<\kappa^+$ in the superscript of $\bU$ we write $\kappa$. If
$\kappa>|\Dom(I)|^+$, we omit it (since then its value does not matter).
\end{Definition}

\begin{Remark}
\label{5.4}
\begin{enumerate}
\item If $2^{<\kappa}+|\Dom(I)|^{<\kappa}\le\mu$ we can find $F\subseteq$
partial functions from $\Dom(I)$ to $\mu$ such that:
\begin{description}
\item[(a)]  $|F|=\bU^{<\kappa}_I(\mu)$,
\item[(b)]  $(\forall f:X\longrightarrow\mu)(\exists Y\in I^+)[f\restriction
Y \in F]$.
\end{description}
\item Such functions (as $\bU^{<\kappa}_I(\mu)$) are investigated in {\bf
pcf} theory (\cite{Sh:g}, \cite[\S6]{Sh:410}, \cite[\S2]{Sh:430},
\cite{Sh:513}).
\item If $I\subseteq J\subseteq {\Cal P}(X)$, then $\bU^{<\kappa}_I(\mu)\le
\bU^{<\kappa}_J(\mu)$, hence by \ref{5.2}(3), and the above
\[\bU^{<\kappa}_{J^0_{\bar\mu}}(\mu)\le \bU^{<\kappa}_{J^1_{\bar\mu}}(\mu)
\le \bU^{<\kappa}_{J^2_{\bar\mu}}(\mu)\le
\bU^{<\kappa}_{J^3_{\bar\mu}}(\mu)\]
and by \ref{5.2}(4) we have $\bU^{<\kappa}_{I^1_{\bar \mu}}\leq
\bU^{<\kappa}_{J^1_{\bar \mu}}(\mu).$
\item On $\IND_\theta(\bar\kappa)$ (see \ref{5.5A} below) see \cite{Sh:513}.
\end{enumerate}
\end{Remark}

\begin{Definition}
\label{5.5A}
$\IND'_\theta(\langle\kappa_n:n<\omega\rangle)$ means that for every model
$M$ with universe $\bigcup\limits_{n<\omega}\kappa_n$ and $\le\theta$
functions, for some $\Gamma\in [\omega]^{\aleph_0}$ and $\eta\in
\prod\limits_{n<\omega}\kappa_n$ we have:
\[n\in\Gamma\quad\Rightarrow\quad\eta(n)\notin c\ell_M\{\eta(\ell):\ell
\ne n\}.\]
\end{Definition}

\begin{Remark}
Actually if $\theta\geq \aleph_0$, this implies that we can fix
$\Gamma$, hence replacing $\langle \kappa_n: n< \omega\rangle$ by an
infinite subsequence we can have $\Gamma=\omega$.
\end{Remark}

\begin{Theorem}
\label{5.5}
\begin{enumerate}
\item If $\mu_n\rightarrow (\kappa_n)^2_{2^\theta}$ and $\IND'_\theta(\langle
\kappa_n:n<\omega\rangle)$ {\em then} $\prod\limits_{n<\omega}\mu_n$ is not
the union of $\le\theta$ sets from $J^1_{\bar \mu}$.
\item If $\theta=\theta^{\aleph_0}$ and $\neg\IND'_\theta(\langle\mu_n:
n<\omega\rangle$) then $\prod\limits_{n<\omega}\mu_n$ is the union of
$\le\theta$ members of $J^1_{\bar\mu}$.
\item If $\lim\sup\limits_n \mu_n$ is $\ge 2$, then $\prod\limits_{n<\omega}
\mu_n\notin J^3_{\bar\mu}$ (so also the other ideals defined above are not
trivial by \ref{5.2}(3), (4)).
\end{enumerate}
\end{Theorem}

\proof 1)\ \ Suppose $\prod\limits_{n<\omega}\mu_n$ is
$\bigcup\limits_{i<\theta} X_i$, and each $X_i\in J^1_{\bar\mu}$. We define
for each $i<\theta$ and $n<k<\omega$ a two-place relation $R^{n,k}_i$ on
$\mu_n$:

\qquad $\beta R^{n,k}_i \gamma$ if and only if

\qquad there are $\eta,\nu\in X_i\subseteq\prod\limits_{\ell<k}\mu_\ell$
such that
\[\eta\restriction [0,n)=\nu\restriction [0,n)\quad\mbox{and }\
\eta\restriction (n,k)=\nu\restriction (n,k)\quad\mbox{and }\ \eta(n)
=\beta,\ \nu(n)=\gamma.\]
Note that $R^{n,k}_i$ is symmetric and
\[n<k_1<k_2\ \&\ \beta R^{n,k_2}_i \gamma\quad \Rightarrow\quad \beta
R^{n,k_1}_i \gamma.\]
As $\mu_n\rightarrow (\kappa_n)^2_{2^\theta}$, we can find $A_n\in
[\mu_n]^{\kappa_n}$ and a truth value $\bt^{n,k}_i$ such that for all
$\beta<\gamma$ from $A_n$, the truth value of $\beta R^{n,k}_i\gamma$ is
$\bt^{n,k}_i$. If for some $i$ the set
\[\Gamma_i=:\{n<\omega:\mbox{ for every }k\in (n,\omega)\mbox{ we have }
\bt^{n,k}_i=\mbox{ true}\}\]
is infinite, we get a contradiction to ``$X_i\in J^1_{\bar \mu}$", so for
some $n(i)<\omega$ we have $n(i)=\sup(\Gamma_i)$.

For each $n<k<\omega$ and $i<\theta$ we define a partial function
$F^{n,k}_i$ from $\prod\limits_{\scriptstyle \ell<k,\atop\scriptstyle\ell
\ne n} A_\ell$ into $A_n$:
\begin{quotation}
\noindent $F(\alpha_0\ldots\alpha_{n-1},\alpha_{n+1},\ldots,\alpha_k)$ is
the first $\beta\in A_n$ such that for some $\eta\in X_i$ we have
\[\begin{array}{c}
\eta\restriction [0,n)=\langle\alpha_0,\ldots,\alpha_{n-1}\rangle,\quad
\eta(n)=\beta,\\
\eta\restriction (n,k)=\langle\alpha_{n+1},\ldots,\alpha_{k-1}\rangle.
\end{array}\]
\end{quotation}
So as $\IND'_\theta(\langle\kappa_n:n<\omega\rangle)$ there is $\eta=\langle
\beta_n:n<\omega\rangle\in\prod\limits_{n<\omega} A_n$ such that for
infinitely many $n$, $\beta_n$ is not in the closure of $\{\beta_\ell:\ell
<\omega,\,\ell\ne n\}$ by the $F^{n,k}_i$'s. As $\eta\in\prod\limits_{n<
\omega} A_n\subseteq\prod\limits_{n<\omega}\mu_n=\bigcup\limits_{i<\theta}
X_i$, necessarily for some $i<\theta$, $\eta\in X_i$.  Let $n\in(n(i),\omega)$
be such that $\beta_n$ is not in the closure of
$\{\beta_\ell:\ell<\omega\mbox{ and }\ell\neq n\}$
and let $k>n$ be such that $\bt^{n,k}_i=\mbox{ false}$. Now $\gamma=:
F^{n,k}_i(\beta_0,\ldots,\beta_{n-1},\beta_{n+1},\ldots,\beta_{k-1})$ is well
defined $\le\beta_n$ (as $\beta_n$ exemplifies that there is such $\beta$) and
is $\ne \beta_n$ (by the choice of $\langle\beta_\ell:\ell<\omega\rangle$), so
by the choice of $n(i)$ (so of $n$, $k$ and earlier of $\bt^{n, k}_i$
and of $A_n$) we get
contradiction to ``$\gamma<\beta_n$ are from $A_n$".

\noindent 2)\ \ Let $M$ be an algebra with universe $\sum\limits_{n<\omega}
\mu_n$ and $\le\theta$ functions (say $F^n_i$ for $i<\theta$, $n<\omega$,
$F^n_i$ is $n$-place) exemplifying $\neg\IND'_\theta(\langle\mu_n:n<\omega
\rangle)$. Let
\[\Gamma=:\{\langle(k_n,i_n):n^*\le n<\omega\rangle:n^*<\omega\mbox{ and }
\bigwedge_n n<k_n<\omega\mbox{ and }i_n<\theta\}.\]
For $\rho=\langle(k_n,i_n):n^*\le n<\omega\rangle\in\Gamma$ let
\[\begin{aligned}
A_\rho=:&\big\{\eta\in\prod\limits_{n<\omega}\mu_n:\mbox{for every }n
\in [n^*,\omega)\mbox{ we have}\\
  &\qquad\eta(n)=F^{k_n-1}_{i_n}\left(\eta(0),\ldots,\eta(n-1),\eta(n+1),
  \ldots,\eta(k_n)\right)\big\}.
\end{aligned}\]
So, by the choice of $M$, $\prod\limits_{n<\omega}\mu_n=\bigcup\limits_{\rho
\in\Gamma} A_\rho$. On the other hand, it is easy to check that $A_\rho\in
J^1_{\bar \mu}$. \hfill$\square_{\ref{5.5}}$

\begin{Theorem}
\label{5.6}
If $\mu=\sum\limits_{n<\omega}\lambda_n$, $\lambda^{\aleph_0}_n<\lambda_{n+1}$
and $\mu<\lambda=\cf(\lambda)<\mu^{+\omega}$\\
then $\bU^{\aleph_0}_{I^0_{\langle\lambda_n:n<\omega\rangle}}(\lambda)=
\lambda$ and even $\bU^{\aleph_0}_{J^3_{\langle \lambda_n:n<\omega\rangle}}
(\lambda)=\lambda$.
\end{Theorem}

\proof See \cite[\S6]{Sh:410}, \cite[\S2]{Sh:430}, and \cite{Sh:513}
for considerably more.

\begin{Lemma}
\label{5.7}
Assume $\lambda>2^{\aleph_0}$ and
\begin{description}
\item[$(*)$(a)]  $\prod\limits_{n<\omega}\mu_n<\mu$ and $\mu^+<\lambda=
\cf(\lambda)<\mu^{\aleph_0}$,
\item[\ \ (b)]   $\hat B_{\bar\mu}\notin I^0_{\bar\mu}$ and $\lim_n\sup\mu_n$
is infinite,
\item[\ \ (c)]   $\bU^{\aleph_0}_{I^0_{\bar\mu}}(\lambda)=\lambda$
(note $I^0_{\bar \mu}$ is not required to be an ideal).
\end{description}
\underline{Then} there is no universal member in ${\mathfrak K}^{rs(p)}_\lambda$.
\end{Lemma}

\proof  Let $S\subseteq\lambda$, $\bar C=\langle C_\delta:\delta\in S\rangle$
guesses clubs of $\lambda$, chosen as in the proof of \ref{3.3} (so $\alpha
\in\nacc(C_\delta)\ \Rightarrow\ \cf(\alpha)>2^{\aleph_0}$). Instead of
defining the relevant invariant we prove the theorem directly, but we could
define it, somewhat cumbersomely (like \cite[III,\S3]{Sh:e}).

Assume $H\in {\mathfrak K}^{rs(p)}_\lambda$ is a pretender to universality;
without loss of generality with the set of elements of $H$ equal to
$\lambda$.

Let $\chi=\beth_7(\lambda)^+$, ${\bar{\mathfrak A}}=\langle {\mathfrak A}_\alpha:
\alpha<\lambda\rangle$ be an increasing continuous sequence of elementary
submodels of $({\Cal H}(\chi),\in,<^*_\chi)$, ${\bar {\mathfrak A}}\restriction
(\alpha+1)\in {\mathfrak A}_{\alpha+1}$, $\|{\mathfrak A}_\alpha\|<\lambda$, ${\mathfrak
A}_\alpha\cap\lambda$ an ordinal, ${\mathfrak A}=\bigcup\limits_{\alpha<\lambda}
{\mathfrak A}_\alpha$ and $\{H,\langle\mu_n:n<\omega\rangle,\mu,\lambda\}\in
{\mathfrak A}_0$, so $B_{\bar \mu},\hat B_{\bar \mu} \in {\mathfrak A}_0$ (where $\bar
\mu=\langle\mu_n:n<\omega\rangle$, of course).

For each $\delta\in S$, let ${\Cal P}_\delta=:[C_\delta]^{\aleph_0}\cap
{\mathfrak A}$. Choose $A_\delta\subseteq C_\delta$ of order type $\omega$
almost disjoint from each $a\in {\Cal P}_\delta$, and from $A_{\delta_1}$ for
$\delta_1\in\delta\cap S$; its existence should be clear as
$2^{\aleph_0}<\lambda< \mu^{\aleph_0}$. So
\begin{description}
\item[$(*)_0$]  every countable $A\in {\mathfrak A}$ is almost disjoint to
$A_\delta$.
\end{description}
By \ref{5.2}(2), $I^0_{\bar\mu}$ is $(<\aleph_1)$-based so by \ref{5.4}(1) and
the assumption (c) we have
\begin{description}
\item[$(*)_1$]  for every $f:\hat B_{\bar\mu}\longrightarrow\lambda$ for some
countable $Y\subseteq \hat B_{\bar\mu}$, $Y\notin I^0_{\bar\mu}$, we have
$f\restriction Y\in {\mathfrak A}$
\end{description}
(remember $(\prod\limits_{n<\omega}\mu_n)^{\aleph_0}=\prod\limits_{n<\omega}
\mu_n$).

\noindent Let $B$ be $\bigoplus\{G^n_{\alpha,i}:n<\omega,\alpha<\lambda,\,i<
\sum\limits_{k<\omega}\mu_k\}$, where
\[G^n_{\alpha,i}=\langle x^n_{\alpha,i}\rangle_{G^n_{\alpha,i}}\cong{\mathbb
Z}/p^{n+1}{\mathbb Z}.\]
\noindent
So $B$, $\hat B$, $\langle(n,\alpha,i,x^n_{\alpha,i}):n<\omega,\alpha<\lambda,
i<\sum\limits_{k<\omega}\mu_k\rangle$ are well defined. Let $G$ be the
subgroup of $\hat B$ generated by:
\[\begin{aligned}
B\cup\big\{x\in\hat B: &\mbox{for some }\delta\in S,\, x\mbox{ is in the
  closure of }\\
  &\bigoplus\{G^n_{\alpha,i}:n<\omega,i<\mu_n,\alpha\mbox{ is the }n\mbox{th
  element of } A_\delta\}\big\}.
\end{aligned}\]
As $\prod\limits_{n<\omega}\mu_n<\mu<\lambda$, clearly $G\in {\mathfrak
K}^{rs(p)}_\lambda$, without loss of generality the set of elements of $G$ is
$\lambda$ and let $h:G\longrightarrow H$ be an embedding. Let
\[E_0=:\{\delta<\lambda:({\mathfrak A}_\delta,h \restriction \delta,\;G
\restriction\delta)\prec({\mathfrak A},h,G)\},\]
\[E=:\{\delta<\lambda:\otp(E_0\cap\delta)=\delta\}.\]
They are clubs of $\lambda$, so for some $\delta\in S$, $C_\delta\subseteq E$
(and $\delta\in E$ for simplicity). Let $\eta_\delta$ enumerate $A_\delta$
increasingly.

There is a natural embedding $g = g_\delta$ of $B_{\bar \mu}$ into $G$:
\[g({}^* t^n_i) = x^n_{\eta_\delta(n),i}.\]
Let $\hat g_\delta$ be the unique extension of $g_\delta$ to an embedding of
$\hat B_{\bar\mu}$ into $G$; those embeddings are pure, (in fact $g''_\delta
(\hat B_{\bar\mu})\setminus g''_\delta(B_\mu)\subseteq G\setminus G\cap
{\mathfrak A}_\delta$). So $h\circ\hat g_\delta$ is an embedding of $\hat B_{\bar
\mu}$ into $H$, not necessarily pure but still an embedding, so the distance
function can become smaller but not zero and
\[h\circ\hat g_\delta(\hat B_{\bar\mu})\setminus h\circ g_\delta(B_\mu)
\subseteq H\setminus {\mathfrak A}_\delta.\]
Remember $\hat B_{\bar\mu}\subseteq {\mathfrak A}_0$ (as it belongs to ${\mathfrak
A}_0$ and has cardinality $\prod\limits_{n<\omega}\mu_n<\lambda$ and
$\lambda\cap {\mathfrak A}_0$ is an ordinal). By $(*)_1$ applied to
$f=h\circ\hat g_\delta$ there is a countable  $Y \subseteq \hat B_{\bar \mu}$
such that $Y \notin I^0_{\bar\mu}$ and $f \restriction Y \in {\mathfrak
A}$. But, from $f \restriction Y$ we shall below reconstruct
some countable sets not almost disjoint to $A_\delta$, reconstruct meaning in
${\mathfrak A}$, in contradiction to $(*)_0$ above.

As $Y\notin I^0_{\bar \mu}$ we can find an infinite
$S^*\subseteq\omega\setminus m^*$ and for $n\in
S^*$, $z_n\in\bigoplus\limits_{\alpha<\mu_n} K^n_\alpha\setminus\{0\}$ and
$y^\ell_n\in\hat B_{\bar\mu}$ (for $\ell<\omega$) such that:
\begin{description}
\item[$(*)_2$]  $z_n+y_{n,\ell}\in\langle Y\rangle_{\hat B_{\bar\mu}}$,\qquad
and
\item[$(*)_3$]  $y_{n,\ell}\in p^\ell\,\hat B_{\bar\mu}$.
\end{description}
Without loss of generality $pz_n=0\ne z_n$ hence $p\,y^\ell_n=0$. Let
\[\nu_\delta(n)=:\min(C_\delta\setminus (\eta_\delta(n)+1)),\quad z^*_n=
(h\circ\hat g_\delta)(z_n)\quad\mbox{ and }\quad y^*_{n,\ell}=(h\circ\hat
g_\delta)(y_{n,\ell}).\]
Now clearly $\hat g_\delta(z_n)=g_\delta(z_n)=x^n_{\eta_\delta(n),i}\in
G\restriction\nu_\delta(n)$, hence $(h\circ\hat g_\delta)(z_n)\notin H
\restriction\eta_\delta(n)$, that is $z^*_n\notin
H\restriction\eta_\delta(n)$.

So $z^*_n\in H_{\nu_\delta(n)}\setminus H_{\eta_\delta(n)}$ belongs to
the $p$-adic closure of $\Rang(f\restriction Y)$. As $H$, $G$, $h$ and
$f\restriction Y$ belongs to ${\mathfrak A}$, also $K$, the closure of
$\Rang(f\restriction Y)$ in $H$ by the $p$-adic topology belongs to
${\mathfrak A}$, and clearly $|K|\leq 2^{\aleph_0}$, hence
\[A^*=\{\alpha\in C_\delta: K\cap H_{\min(C_\delta\setminus
(\alpha+1))}\setminus H_\alpha \mbox{ is not empty}\}\]
is a subset of $C_\delta$ of cardinality $\leq 2^{\aleph_0}$ which
belongs to ${\mathfrak A}$, hence $[A^*]^{\aleph_0}\subseteq {\mathfrak A}$
but $A_\delta\subseteq A^*$ so $A_\delta\in {\mathfrak A}$, a contradiction.
 \hfill$\square_{\ref{5.7}}$

\section{Implications between the existence of universals}

\begin{Theorem}
\label{6.1}
Let $\bar n=\langle n_i:i<\omega\rangle$, $n_i\in [1,\omega)$. Remember
\[J^2_{\bar n}=\{A\subseteq\prod_{i<\omega} n_i:A \mbox{ is nowhere
dense}\}.\]
Assume $\lambda\ge 2^{\aleph_0}$, $T^{\aleph_0}_{J^3_{\bar n}}(\lambda)=
\lambda$ or just $T^{\aleph_0}_{J^2_{\bar n}}(\lambda)=\lambda$ for every
such $\bar n$, and
\[n<\omega\quad\Rightarrow\quad\lambda_n\le\lambda_{n+1}\le\lambda_\omega
=\lambda\quad\mbox{ and}\]
\[\lambda\le\prod_{n<\omega}\lambda_n\quad\mbox{ and }\quad\bar\lambda=
\langle\lambda_i:i\le\omega\rangle.\]
\begin{enumerate}
\item  If in ${\mathfrak K}^{fc}_{\bar \lambda}$ there is a universal member

{\em then} in ${\mathfrak K}^{rs(p)}_\lambda$ there is a universal member.
\item  If in ${\mathfrak K}^{fc}_\lambda$ there is a universal member for
${\mathfrak K}^{fc}_{\bar \lambda}$

{\em then} in
\[{\mathfrak K}^{rs(p)}_{\bar\lambda}=:\{G\in {\mathfrak K}^{rs(p)}_\lambda:\lambda_n
(G)\le\lambda\}\]
there is a universal member (even for ${\mathfrak K}^{rs(p)}_\lambda$).
\end{enumerate}
($\lambda_n(G)$ were defined in \ref{1.1}).
\end{Theorem}

\begin{Remark}
\begin{enumerate}
\item Similarly for ``there are $M_i\in {\mathfrak K}_{\lambda_1}$ ($i<\theta$)
with $\langle M_i:i<\theta\rangle$ being universal for ${\mathfrak K}_\lambda$''.
\item The parallel of \ref{1.1} holds for ${\mathfrak K}^{fc}_\lambda$.
\item By \S5 only the case $\lambda$ singular or $\lambda=\mu^+\ \&\
\cf(\mu)=\aleph_0\ \& \ (\forall \alpha<
\mu)(|\alpha|^{\aleph_0}<\mu)$
is of interest for \ref{6.1}.
\end{enumerate}
\end{Remark}

\proof 1)\ \ By \ref{1.1}, (2) $\Rightarrow$ (1).

More elaborately, by part (2) of \ref{6.1} below there is $H\in {\mathfrak
K}^{rs(p)}_{\bar \lambda}$ which is universal in ${\mathfrak
K}^{rs(p)}_{\bar \lambda}$. Clearly $|G|=\lambda$ so $H\in {\mathfrak
K}^{rs(p)}_\lambda$, hence for proving part (1) of \ref{6.1} it
suffices to prove that $H$ is a universal member of ${\mathfrak
K}^{rs(p)}_\lambda$. So let $G\in {\mathfrak K}^{rs(p)}_\lambda$, and we
shall prove that it is embeddable into $H$. By \ref{1.1} there is $G'$
such that $G\subseteq G'\in {\mathfrak K}^{rs(p)}_{\bar \lambda}$. By the
choice of $H$ there is an embedding $h$ of $G'$ into $H$. So
$h\restriction G$ is an embedding of $G$ into $H$, as required.

\noindent 2)\ \ Let $T^*$ be a universal member of ${\mathfrak K}^{fc}_{\bar
\lambda}$ (see \S2) and let $P_\alpha = P^{T^*}_\alpha$.

Let $\chi>2^\lambda$. Without loss of generality $P_n=\{n\}\times
\lambda_n$, $P_\omega=\lambda$. Let
\[B_0=\bigoplus\{G^n_t:n<\omega,t\in P_n \},\]
\[B_1=\bigoplus \{G^n_t: n< \omega\mbox{ and }t\in P_n\},\]
where $G^n_t\cong {\mathbb Z}/p^{n+1}{\mathbb Z}$, $G^n_t$ is generated by
$x^n_t$. Let ${\mathfrak B}\prec ({\Cal H}(\chi),\in,<^*_\chi)$, $\|{\mathfrak B}\|=
\lambda$, $\lambda+1\subseteq {\mathfrak B}$, $T^*\in {\mathfrak B}$, hence
$B_0$, $B_1\in {\mathfrak B}$ and $\hat B_0, \hat B_1\in {\mathfrak B}$ (the
torsion completion of $B$). Let $G^* =\hat B_1\cap {\mathfrak B}$.

Let us prove that $G^*$ is universal for ${\mathfrak K}^{rs(p)}_{\bar
\lambda}$ (by \ref{1.1} this suffices).
Let $G \in {\mathfrak K}^{rs(p)}_{\lambda}$, so by \ref{1.1} without loss of
generality $B_0 \subseteq G\subseteq\hat B_0$. We define $R$:
\[\begin{aligned}
R=\big\{\eta:&\eta\in\prod\limits_{n<\omega}\lambda_n\mbox{ and for
some }
x\in G\mbox{ letting }\\
	&x=\sum\{a^n_i\,p^{n-k}\,x^n_i:n<\omega,i\in w_n(x)\}\mbox{ where }\\
	&w_n(x)\in [\lambda_n]^{<\aleph_0},a^n_i\,p^{n-k}\,x^n_i\ne 0\mbox{
we have }\\
	&\bigwedge\limits_n\eta(n)\in w_n(x)\cup \{\ell:
\ell+|w_n(x)|\leq n\}\big\}.
\end{aligned}\]
Lastly let $M =:(R\cup\bigcup\limits_{n<\omega}\{n\}\times\lambda_n,\,P_n,\,
F_n)_{n<\omega}$ where $P_n=\{n\}\times\lambda_n$ and $F_n(\eta)=(n,\eta(n))$,
so clearly $M\in {\mathfrak K}^{fc}_{\bar\lambda}$. Consequently, there is an
embedding $g:M\longrightarrow T^*$, so $g$ maps $\{n\}\times\lambda_n$ into
$P^{T^*}_n$ and $g$ maps $R$ into $P^{T^*}_\omega$. Let $g(n,\alpha)=
(n,g_n(\alpha))$ (i.e. this defines $g_n$). Clearly $g\restriction (\cup
P^M_n)=g\restriction (\bigcup\limits_n\{n\}\times\lambda_n)$ induces an
embedding $g^*$ of $B_0$ to $B_1$ (by mapping the generators into the
generators).

\noindent The problem is why:
\begin{description}
\item[$(*)$]  if $x=\sum\{a^n_i\,p^{n-k}\,x^n_i:n<\omega,i\in w_n(x)\}\in G$

{\em then} $g^*(x)=\sum\{a^n_i\,p^{n-k}\,g^*(x^n_i):n<\omega,i\in w_n(x)\}\in
G^*$.
\end{description}
As $G^*=\hat B_1\cap {\mathfrak B}$, and $2^{\aleph_0}+1\subseteq {\mathfrak B}$, it is
enough to prove $\langle g^{\prime\prime}(w_n(x)):n<\omega\rangle\in
{\mathfrak B}$. Now for
notational simplicity $\bigwedge\limits_n [|w_n(x)|\ge n+1]$ (we can add an
element of $G^*\cap {\mathfrak B}$ or just repeat the arguments). For each $\eta
\in\prod\limits_{n<\omega} w_n(x)$ we know that
\[g(\eta)=\langle g(\eta(n)):n<\omega\rangle\in T^*\quad\mbox{ hence is in }
{\mathfrak B}\]
(as $T^*\in {\mathfrak B}$, $|T^*|\le\lambda$).  Now by assumption there is
$A\subseteq\prod\limits_{n<\omega} w_n(x)$ which is not nowhere dense
such that $g
\restriction A\in {\mathfrak B}$, hence for some $n^*$ and $\eta^* \in
\prod\limits_{\ell<n^*}w_\ell(x)$, $A$ is dense above $\eta^*$ (in
$\prod\limits_{n<\omega} w_n(x)$). Hence
\[\langle\{\eta(n):\eta\in A\}:n^* \le n<\omega\rangle=\langle w_n[x]:n^*\le
n<\omega\rangle,\]
but the former is in ${\mathfrak B}$ as $A\in {\mathfrak B}$, and from the latter the
desired conclusion follows.  \hfill$\square_{\ref{6.1}}$

\section{Non-existence of universals for trees with small density}
For simplicity we deal below with the case $\delta=\omega$, but the proof
works in general (as for ${\mathfrak K}^{fr}_{\bar\lambda}$ in \S2) with minor
changes.
Section 1 hinted we should look at ${\mathfrak K}^{tr}_{\bar\lambda}$ not only
for the case
$\bar\lambda=\langle\lambda:\alpha\le\omega\rangle$ (i.e. ${\mathfrak
K}^{tr}_\lambda$), but in particular for
\[\bar\lambda=\langle\lambda_n:n<\omega\rangle\char 94\langle\lambda\rangle,
\qquad \lambda^{\aleph_0}_n<\lambda_{n+1}<\mu<\lambda=\cf(\lambda)<
\mu^{\aleph_0}.\]
Here we get for this class (embeddings are required to preserve levels),
results stronger than the ones we got for the classes of Abelian groups we
have considered.

\begin{Theorem}
\label{7.1}
Assume that
\begin{description}
\item[(a)] $\bar\lambda=\langle\lambda_\alpha:\alpha\le\omega\rangle$,
$\lambda_n<\lambda_{n+1}<\lambda_\omega$, $\lambda=\lambda_\omega$,
all are regulars,
\item[(b)] $D$ is a filter on $\omega$ containing cobounded sets,
\item[(c)] $\tcf(\prod \lambda_n /D)=\lambda$ (indeed, we mean $=$, we could
just use $\lambda\in\pcf_D(\{\lambda_n:n<\omega\})$),
\item[(d)] $(\sum\limits_{n<\omega}\lambda_n)^+<\lambda<\prod\limits_{n<
\omega}\lambda_n$.
\item[(e)] ${\Cal P}({\mathfrak a})/ J_{<\lambda} 
(\{\lambda_n:n<\omega\})$ is infinite.
\end{description}
\underline{Then} there is no universal member in ${\mathfrak
K}^{tr}_{\bar\lambda}$.
\end{Theorem}

\proof Clearly if $n_i<n_{i+1}<\omega$ for $i<\omega$ and 
$\bar{\lambda}'=\langle
\lambda_{n_i}:i<\omega\rangle \smallfrown \langle \lambda \rangle$ and is ${\mathfrak
K}^{tr}_{\bar{\lambda}}$ there is a universal member \underline{then} in
${\mathfrak
K}^{tr}_{\bar{\lambda}}$, there is a universal member. So without loss of
generality $\lambda=\max\pcf\{\lambda_n:n<\omega\}$ without loss of
generality $D-\{A:$ max proof $\{\lambda_n:n \in \omega \setminus A\}
<\lambda\}$. Let $\mu=\sum\limits_{n<\omega} \lambda_n$. 

We now notice that there is a sequence $\bar{\mathcal P}=\langle {\mathcal 
P}_\alpha: \mu<\alpha<\lambda\rangle$ such that:
\begin{enumerate}
\item  $|{\Cal P}_\alpha|<\lambda$,
\item  $a\in {\Cal P}_\alpha\quad\Rightarrow\quad a$ is a closed subset of
$\alpha$
of order type $\leq\mu$,
\item  $a\in\bigcup\limits_{\alpha<\lambda} {\Cal P}_\alpha\ 
\&\ \beta\in\nacc(a)
\quad \Rightarrow\quad a\cap\beta\in {\Cal P}_\beta$,
\item  For all club subsets $E$ of $\lambda$, there are stationarily many
$\delta$ for which there is an $a\in\bigcup\limits_{\alpha<\lambda} P_\alpha$
such that
\[\cf(\delta)=\aleph_0\ \&\ a\in P_\delta\ \&\ \otp(a)=\mu
\ \&\ a\subseteq E.\]
\end{enumerate}
[Why? If $\lambda=\mu^{++}$, then it is the
successor of a regular, so we use \cite[\S4]{Sh:351}, i.e.
\[\{\alpha<\lambda:\cf(\alpha)\le\mu\}\]
is the union of $\leq \mu^+$ sets with squares.\\
If $\lambda>\mu^{++}$, then we can use
\cite[\S1]{Sh:420}, which guarantees that there is a stationary $S\in
I[\lambda]$.and note that we do not require $\sup (a)=\alpha$ even if
$a \in {\Cal P}_\alpha$ and otp$({\rm a})=\mu$]

We can now find a sequence
\[\langle f_\alpha,g_{\alpha,a}:\alpha<\lambda \quad \hbox{and} 
\quad a\in {\Cal P}_\alpha\rangle\]
such that:
\begin{description}
\item[(a)] $\bar f=\langle f_\alpha:\alpha<\lambda\rangle$ is a
$<_D$-increasing cofinal sequence in $\prod\limits_{n<\omega}\lambda_n$,
\item[(b)] $g_{\alpha,a}\in\prod\limits_{n<\omega}\lambda_n$,
\item[(c)] $\bigwedge\limits_{\beta<\alpha} f_\beta<_D g_{\alpha,a}<_D
f_{\alpha+1}$,
\item[(d)] $\lambda_n>|a|\ \&\ \beta\in\nacc(a)\quad \Rightarrow\quad
g_{\beta,a\cap\beta}(n)<g_{\alpha,a}(n)$.
\item[(e)] for every $f\in \prod \lambda_n$ for some $\alpha$ we have
$f<f_\alpha$
\end{description}
[How?  Choose $\bar f$ by $\tcf(\prod\limits_{n<\omega}\lambda_n/D)=\lambda
+ \max\pcf \{\lambda_n:n<\omega\}=\lambda$.
Then choose $g$'s by induction, possibly throwing out some of the
$f$'s; this is from \cite[II, \S1]{Sh:g}.]

Let $T\in {\mathfrak K}^{tr}_{\bar \lambda}$.

We introduce for $x\in\lev_\omega(T)$ and $\ell<\omega$ the notation
$F^T_\ell(x) = F_\ell(x)$ to denote the unique member of $\lev_\ell(T)$ which
is below $x$ in the tree order of $T$.

\noindent For $a\in\bigcup\limits_{\alpha<\lambda} {\Cal P}_\alpha$, let $a=\{
\alpha_{a,\xi}:\xi<\otp(a)\}$ be an increasing enumeration and let
$a(\zeta)=\{\alpha_{a,\xi}:\xi<\zeta\}$. We shall consider
two cases.  In the first one, we assume that the following statement $(*)$
holds.  In this case, the proof is easier, and maybe $(*)$ always holds for
some $D$, but we do not know this at present.
\begin{description}
\item[{{$(*)$}}]  There is a partition $\langle A_n:n < \omega \rangle$ of
$\omega$ into sets not disjoint to any member of $D$.
\end{description}
Let $A_n \in D^+$ be pairwise disjoint, for $n<\omega$, let $D_n$ be
the filter
generated by
$D$ and
$A_n$.  Let for $a\in\bigcup\limits_{\alpha<\lambda} {\Cal P}_\alpha$ with
$\otp(a)=
\mu$, and for $x\in\lev_\omega(T)$,
\[\inv(x,a,T)=:\langle\xi_n(x,a,T):n<\omega\rangle,\]
where
\[\begin{aligned}
\xi_n(x,a,T)=:\min\big\{\xi\leq \otp(a)\!:&\mbox{ if } 
\xi<\otp (a) \mbox{ then } \omega \xi+\omega \leq \otp (a) \\ 
        & \mbox{ and for some } m<\omega \mbox{ we have }\\
	& \langle F^T_\ell(x):
\ell<\omega\rangle<_{D_n} g_{\alpha', a'}\mbox{ where }\\
	&\alpha'=\alpha_{a,\omega\xi+m}\mbox{ and } a'=
a\cap\alpha'\big\}.
\end{aligned}\]
Let
\[\INv(a,T)=:\{\inv(x,a,T):x\in T\ \&\ \lev_T(x)=\omega\},\]
\[\begin{aligned}
\INV(T)=:\big\{c:&\mbox{for every club } E\subseteq\lambda,\mbox{ for some }
\delta\mbox{ and } a\in P\\
  &\mbox{we have }\otp(a)=\mu\ \&\ a\subseteq E\ \&\ a\in P_\delta\\
  &\mbox{and for some } x\in T\mbox{ of } \lev_T(x)=\omega,\ c=\inv(x,a,T)
  \big\}.
\end{aligned}\]
(Alternatively, we could have looked at the function giving each $a$ the value
$\INv(a,T)$, and then divide by a suitable club guessing ideal as in
the proof in \S3, see Definition \ref{3.7}.)

\noindent Clearly
\medskip

\noindent{\bf Fact}:\hspace{0.15in}  $\INV(T)$ has cardinality $\le\lambda$.
\medskip

The main point is the following
\medskip

\noindent{\bf Main Fact}:\hspace{0.15in}  If $\bh:T^1\longrightarrow T^2$ is
an embedding, {\em then\/}
\[\INV(T^1)\subseteq\INV(T^2).\]
\medskip

\noindent{\em Proof of the {\bf Main Fact} under $(*)$}\ \ \ We define for $n
\in\omega$
\[E_n=:\big\{\delta<\lambda_n:\,\delta>\bigcup_{\ell<n}\lambda_\ell\mbox{ and
}\left(\forall x\in\lev_{n+1}(T^1)\right)(F_n(\bh(x))<\delta
\Leftrightarrow F_n(x)<\delta)
\big\}.\]
We similarly define $E_\omega$, so $E_n$ ($n\in\omega$) and $E_\omega$ are
clubs (of $\lambda_n$ and $\lambda$ respectively). Now suppose $c\in\INV(T_1)
\setminus\INV(T_2)$. Without loss of generality $E_\omega$ is (also) a club of
$\lambda$ which exemplifies that $c\notin\INV(T_2)$. For $h\in
\prod\limits_{n<\omega}\lambda_n$, let
\[h^+(n)=:\min(E_n\setminus (h(n)+1),\quad\mbox{ and }\quad\beta[h]=\min\{\beta
<\lambda:h<f_\beta\}.\]
(Note that $h<f_{\beta[h]}$, not just $h<_D f_{\beta[h]}$ check?.) For a
sequence
$\langle h_i:i<i^*\rangle$ of functions from $\prod\limits_{n<\omega}
\lambda_n$, we use $\langle h_i:i<i^* \rangle^+$ for $\langle h^+_i:i<i^*
\rangle$. Now let
\[E^*=:\big\{\delta<\lambda:\mbox{if }\alpha<\delta\mbox{ then }
\beta[f^+_\alpha]<\delta\mbox{ and }\delta\in\acc(E_\omega)\big\}.\]
Thus $E^*$ is a club of $\lambda$.  Since $c\in\INV(T_1)$, there is $\delta<
\lambda$ and $a\in P_\delta$ such that for some $x\in\lev_\omega(T_1)$ we have
\[a\subseteq E^*\ \&\ \otp(a)=\mu\ \&\ c=\inv(x,a,T_1).\]
Let for $n\in\omega$, $\xi_n=:\xi_n(x,a,T_1)$, so $c=\langle\xi_n:n<\omega
\rangle$. Also let for $\xi<\mu$, $\alpha_\xi=:
\alpha_{a,\xi}$, so $a=\{\alpha_\xi:\xi<\mu\}$ is an increasing enumeration.
Now fix an $n<\omega$ and consider
$\bh(x)$. Then we know that for some $m$
\begin{description}
\item[($\alpha$)]  $\langle F^{T_1}_\ell(x):\ell<\omega\rangle<_{D_n}
g_{\alpha'}$ where $\alpha'=\alpha_{\omega\xi_n+m}$\qquad and
\item[($\beta$)]   for no $\xi<\xi_n$ is there such an $m$.
\end{description}
Now let us look at $F^{T_1}_\ell(x)$ and $F^{T_2}_\ell({\mathbf h}(x))$. They are
not necessarily equal, but
\begin{description}
\item[($\gamma$)]  $\min(E_\ell\setminus (F^{T_1}_\ell(x)+1))=\min(E_\ell
\setminus (F^{T_2}_\ell(\bh(x)+1))$
\end{description}
(by the definition of $E_\ell$).  Hence
\begin{description}
\item[($\delta$)]  $\langle F^{T_1}_\ell(x):\ell<\omega\rangle^+=\langle
F^{T_2}_\ell(\bh(x)):\ell<\omega\rangle^+$.
\end{description}
Now note that by the choice of $g$'s
\begin{description}
\item[($\varepsilon$)] $(g_{\alpha_\varepsilon,a\cap\alpha_\varepsilon})^+
<_{D_n} g_{\alpha_{\varepsilon+1},a\cap\alpha_{\varepsilon+1}}$.
\end{description}
\relax From $(\delta)$ and $(\varepsilon)$ it follows that $\xi_n(\bh(x),a,
T^2)=\xi_n(x,a,T^1)$.  Hence $c\in\INV(T^2)$.
\hfill$\square_{\mbox{Main Fact}}$
\medskip

Now it clearly suffices to prove:
\medskip

\noindent{\bf Fact A:}\hspace{0.15in}  For each $c=\langle\xi_n:n<\omega
\rangle\in {}^\omega \mu$ 
we can find a $T\in {\mathfrak K}^{tr}_{\bar\lambda}$ such that $c\in\INV(T)$.
\medskip

\noindent{\em Proof of the Fact A in case $(*)$ holds}\ \ \ For each $a\in
\bigcup\limits_{\delta<\lambda} {\Cal P}_\delta$ with $\otp(a)=
\mu$ we define $x_{c,a}=:\langle x_{c,a}(\ell):\ell<\omega
\rangle$ by:
\begin{quotation}
if $\ell\in A_n$, then $x_{c,a}(\ell)=g_{\alpha_{a,\omega\xi_n+1}} (\ell)$.
\end{quotation}
Let
\[T=\bigcup_{n<\omega}\prod_{\ell<n}\lambda_\ell\cup\big\{x_{c,a}:a\in
\bigcup_{\delta<\lambda} {\Cal P}_\delta\ \&\ \otp(a)=\mu
\big\}.\]
We order $T$ by $\triangleleft$.

It is easy to check that $T$ is as required. \hfill$\square_A$
\medskip

Now we are left to deal with the case that $(*)$ does not hold.  Let
\[\pcf(\{\lambda_n:n<\omega\})=\{\kappa_\alpha:\alpha\le\alpha^*\}\]
be an enumeration in increasing order so in particular
\[\kappa_{\alpha^*}=\max\pcf(\{\lambda_n:n<\omega\}).\]
Without loss of generality $\kappa_{\alpha^*}=\lambda$ (by throwing out some
elements if necessary; see beginning of the proof) and
$\lambda\cap\pcf(\{\lambda_n:n<\omega\})$ has no last element (this appears
explicitly in
\cite{Sh:g}, but is also straightforward from the pcf theorem).  In
particular, $\alpha^*$
is a limit ordinal.  Hence, without loss of generality
\[D=\big\{A\subseteq\omega:\lambda>\max\pcf\{\lambda_n:n\in\omega\setminus A\}
\big\}.\]
Let $\langle {\mathfrak a}_{\kappa_\alpha}:\alpha\le\alpha^*\rangle$ be a
generating sequence for $\pcf(\{\lambda_n:n<\omega\})$, i.e.
\[\max\pcf({\mathfrak a}_{\kappa_\alpha})=\kappa_\alpha\quad\mbox{ and }\quad
\kappa_\alpha\notin\pcf(\{\lambda_n:n<\omega\}\setminus {\mathfrak
a}_{\kappa_\alpha}).\]
(The existence of such a sequence is part of the pcf theorem).  Without
loss of generality,
\[{\mathfrak a}_{\alpha^*} = \{ \lambda_n:n < \omega\}.\]
Now note

\begin{Observation}
If $\cf(\alpha^*)=\aleph_0$, then $(*)$ holds.
\end{Observation}
Why? Let $\langle\alpha(n):n<\omega\rangle$ be a strictly increasing cofinal
sequence in $\alpha^*$. Let $\langle B_n:n<\omega\rangle$ partition $\omega$
into infinite pairwise disjoint sets and let
\[A_\ell=:\big\{k<\omega:\bigvee_{n\in B_\ell}[\lambda_k\in {\mathfrak
a}_{\kappa_{\alpha(n)}}\setminus\bigcup_{m<n} {\mathfrak a}_{\kappa_{\alpha(m)}}]
\big\}.\]
To check that this choice of $\langle A_\ell:\ell<\omega\rangle$ works, recall
that for all $\alpha$ we know that ${\mathfrak a}_{\kappa_\alpha}$ does not
belong to
the ideal generated by $\{{\mathfrak a}_{\kappa_\beta}: \beta<\alpha\}$ and use
the pcf calculus (i.e. $n \in B_\ell \Rightarrow \lambda_n \in \pcf({\mathfrak
a}_{\kappa_{\alpha(n)}}) \subseteq \pcf \{\lambda_k:k\in A_\ell\}$ so
$\max\pcf\{\lambda_k:k\in A_\ell\}) \geq \max\pcf \{\lambda_n:n\in
B_\ell\} \geq \sum\limits_{n\in B_\ell} \kappa_{\alpha(n)}$ but
$\pcf(\{\lambda_n:n<\omega\})
\setminus \sum\limits_{n\in B_\ell} \kappa_{\alpha(n)}=\{\lambda\}$ so we
are done).
\hfill$\square$

Now let us go back to the general case, assuming $\cf(\alpha^*)>\aleph_0$.
Our problem is the possibility that
\[{\Cal P}(\{\lambda_n:n<\omega\})/J_{<\lambda}[\{\lambda_n:n<\omega\}].\]
is finite. Let now $A_\alpha=:\{n:\lambda_n\in {\mathfrak a}_\alpha\}$, and
\[\begin{array}{lll}
J_\alpha  &=: &\big\{A\subseteq\omega:\max(\pcf\{\lambda_\ell:\ell\in A\})<
\kappa_\alpha\big\}\\
J'_\alpha &=: &\big\{A\subseteq\omega:\max\pcf(\{\lambda_\ell:\ell\in A\}\cap
{\mathfrak a}_{\kappa_\alpha})<\kappa_\alpha\big\}.
\end{array}\]

For $T\in {\mathfrak K}^{tr}_{\bar\lambda}$, $x\in\lev_\omega(T)$,
$\alpha<\alpha^*$ and $a\in\bigcup\limits_{\delta<\lambda} {\Cal P}_\delta$:

\[\begin{aligned}
\xi^*_\alpha(x,a,T)=:\min\big\{\xi:&\bigvee_m [\langle F^T_\ell(x):\ell<
\omega\rangle <_{J'_\alpha} g_{\alpha', a'}]\mbox{ where }\\
	&\qquad \alpha'=\alpha_{a,\omega \xi+m}\mbox{ and }a'=a\cap
\alpha'\big\}.
  \end{aligned}
\]
Let
\[\inv_\alpha(x,a,T)=:\langle\xi^*_{\alpha+n}(x,a,T):n<\omega\rangle,\]
\[\INv(a,T)=:\big\{\inv_\alpha(x,a,T):x\in T\ \&\ \alpha<\alpha^* \ \&\
\lev_T(x)=\omega\big\},\]
and

$\INV(T)=\big\{c\in\prod\limits_{n<\omega}\lambda_n:$ each  
 $\xi\in {\rm Rang}(i)$  is a succession ordinal and for every club 
$E^*$ of $\lambda$  for some $a\in\bigcup\limits_{\delta<\lambda} 
{\Cal P}_\delta$
with $\otp(a)=\mu$  and  a $\subseteq E^*$  for arbitrarily
 large $\alpha<\alpha^*$, there is $x\in\lev_\omega(T)$ 
 such that $\inv_\alpha(x,a,T)= c\big\}$.

As before, the point is to prove the Main Fact.
\medskip

\noindent{\em Proof of the {\bf Main Fact} in general}\ \ \  Suppose $\bh:
T^1\longrightarrow T^2$ and $c\in\INV(T^1)\setminus\INV(T^2)$. Let $E'$ be
a club of $\lambda$ which witnesses that $c\notin\INV(T^2)$. We define
$E_n,E_\omega$ as before, as well as $E^*$ ($\subseteq E_\omega $).
Now let us choose $a\in
\bigcup\limits_{\delta<\lambda} {\Cal P}_\delta$ with $a\subseteq E^*$ and
$\otp(a)=
\mu$ a witnessing $c\in \INV (T')$ for the club $E^*$. So
$a=\{\alpha_{a,\xi}:\xi<
\mu\}$, which we shorten as $a=\{\alpha_\xi:\xi<
\mu\}$. For each $\xi<\mu$, as
before, we know that
\[(g_{\alpha_\xi,a\cap\alpha_\xi})^+<_{J_{\alpha^*}} g_{\alpha_{\xi+1},a\cap
\alpha_{\xi+1}}.\]
By the definition of $<_{J_{\alpha^*}}$, there are $\beta_{\xi,\ell}<\alpha^*$
($\ell<\ell_\xi$) such that
\[\{\ell:g^+_{\alpha_\xi,a\cap\alpha_\xi}(\ell)\ge g_{\alpha_{\xi+1},a\cap
\alpha_{\xi+1}}(\ell)\}\subseteq\bigcup_{\ell<\ell_\xi} A_{\beta_{\xi,\ell}}.\]
Let $c=\langle\xi_n:n<\omega\rangle$ and let
\[\Upsilon=\big\{\beta_{\xi,\ell}:\mbox{for some } n\mbox{ and } m\mbox{ we
have }\xi=\omega\xi_n+m\ \&\ \ell<\ell_\xi \mbox{ or }\]
$ \xi=\omega(\xi_n -1)+m \& \ell<\ell_\xi\big\}.$
Thus $\Upsilon\subseteq\alpha^*$ is countable. Since $\cf(\alpha^*)>\aleph_0$,
the set $\Upsilon$ is bounded in $\alpha^*$. Now we know that $c$ appears as
an invariant for $a$ and arbitrarily large $\delta<\alpha^*$, for some
$x_{a,\delta}\in\lev_\omega(T_1)$ we have $c={\rm inv}_\alpha
(x_{a,\delta},a,T_1)$. If
$\delta>\sup(\Upsilon)$,
$c\in\INV(T^2)$ is exemplified by $a,\delta,\bh(x_{\alpha,\delta})$, just
as before.
\hfill$\square$
\medskip

We still have to prove that every $c=\langle\xi_n:n<\omega\rangle$ appears as
an invariant; i.e. the parallel of Fact A.
\medskip

\noindent{\em Proof of Fact A in the general case:}\ \ \ Define for each $a\in
\bigcup\limits_{\delta<\lambda} {\Cal P}_\delta$ with $\otp(a)=
\mu$ and $\beta<\alpha^*$
\[x_{c,a,\beta}=\langle x_{c,a,\beta}(\ell):\ell<\omega\rangle,\]
where
\[x_{c,a,\beta}(\ell)=\left\{
\begin{array}{lll}
  \alpha_{a,\omega\xi_n+\delta} &\mbox{\underbar{if}} &\lambda_\ell\in {\mathfrak
a}_{\beta+k}\setminus\bigcup\limits_{k'<k} {\mathfrak a}_{\beta+k'}\\
  0 &\mbox{\underbar{if}} &\lambda_\ell\notin {\mathfrak a}_{\beta+k}\mbox{ for
any } k<\omega.
\end{array}\right.\]
Form the tree as before. Now for any club $E$ of $\lambda$, we can find $a\in
\bigcup\limits_{\delta<\lambda} P_\delta$ with $\otp(a)=
\mu$, $a\subseteq E$ such that $\langle x_{c,a,\beta}:\beta<\alpha^*
\rangle$ shows that $c\in\INV(T)$. \hfill$\square_{\ref{7.1}}$

\begin{Remark}
\begin{enumerate}
\item Clearly, this proof shows not only that there is no one $T$ which is
universal for ${\mathfrak K}^{tr}_{\bar\lambda}$, but that any sequence of
$<\prod\limits_{n<\omega}\lambda_n$ trees will fail. This occurs generally in
this paper, as we have tried to mention in each particular case.
\item The case ``$\lambda<2^{\aleph_0}$" is included in the theorem, though
for the Abelian group application the $\bigwedge\limits_{n<\omega}
\lambda^{\aleph_0}_n<\lambda_{n+1}$ is necessary.
\end{enumerate}
\end{Remark}

\begin{Remark}
\label{7.1A}
\begin{enumerate}
\item If $\mu^+< \lambda=\cf(\lambda)<\chi<\mu^{\aleph_0}$ and
$\chi^{[\lambda]}<\mu^{\aleph_0}$
(or at least $T_{\id^a(\bar C)}(\chi)<\mu^{\aleph_0}$) we can get the results
for ``no $M \in {\mathfrak K}^x_\chi$ is universal for ${\mathfrak K}^x_\lambda$", see
\S8 (and \cite{Sh:456}).
\end{enumerate}
\end{Remark}

\noindent We can below (and subsequently in \S8) use $J^3_{\bar m}$ as in \S6.

\begin{Theorem}
\label{7.2}
\begin{enumerate}
\item Assume that $2^{\aleph_0}<\lambda_0$,
$\bar\lambda=\langle\lambda_n:n<\omega
\rangle\char 94\langle\lambda\rangle$, $\mu=\sum\limits_{n<\omega}\lambda_n$,
$\lambda_n<\lambda_{n+1}$, $\mu^+<\lambda=\cf(\lambda)<\mu^{\aleph_0}$.\\
\underbar{If}, for simplicity, $\bar m=\langle m_i:i<\omega\rangle$ 
and $\prod m_i<\lambda_0$ with 
$m_i\in [2,\lambda_0)$ and $\bU^{<\mu}_{J^2_{\bar m}}
(\lambda)=\lambda$ (remember
\[J^2_{\bar m}=\{A\subseteq\prod_{i<\omega} m_i:A \mbox{ is nowhere dense}\}\]
and definition \ref{5.3}),\\
\underbar{then} in ${\mathfrak K}^{tr}_{\bar\lambda}$ there is no universal
member (hence also in ${\mathfrak K}^{tr}_{\langle \lambda:i\leq \omega\rangle})$.
\item If we replace ${\bf U}^{<\mu}_{J^2_{\bar{m}}} (\lambda)=\lambda$ by
${\bf U}^{<\mu}_{J^{\gamma}_{1}}(\lambda) < \cov (\mu,\mu.\aleph_1,2)$ whee
$\Gamma\in(J^2_{\bar{m}})^+, J'_{\bar{m}}=J^2_m \restriction \Gamma$ and
$|\Gamma| <{\bf
U}^{<\Gamma}_{J'_{\bar{m}}}(\lambda)$ and $\mu^*\leq \mu$, 
$(\forall \mu_1) [\mu_1<\mu\Rightarrow \mu_1^{\aleph_0}<\mu^*]$ 
\underline{then} the conclusion still holds.
\end{enumerate}
\end{Theorem}

\proof  1)\ \ Let $\lambda^*={\bf U}^{<\mu}_{J'_{\bar{m}}}(\lambda)$. Let
$S\subseteq\lambda$, $\bar C=\langle C_\delta:\delta\in S
\rangle$ be a club guessing sequence on $\lambda$ with $\otp(C_\delta)\ge
\mu$. We can choose ${\bar{\mathfrak A}}=\langle {\mathfrak A}_\alpha:
\alpha<\lambda\rangle$, satisfying $\{J^2_{\bar m}$, 
$,\bar{C},\langle\lambda_n:n<\omega\rangle$, 
$\mu,\lambda,\bar{m}$, $T^*\}$ 
belong to 
${\mathfrak A}_0$ ($T^*$ is a candidate for the universal), 
and ${\mathfrak A}_\alpha\prec({\Cal H}(\chi),\in,<^*_\chi)$,
$\chi= \beth_7(\lambda)^+$, $\|{\mathfrak A}_\alpha\|<\lambda$, ${\mathfrak A}_\alpha$
increasingly continuous, $\langle {\mathfrak A}_\beta:\beta\le\alpha\rangle\in
{\mathfrak A}_{\alpha+1}$, ${\mathfrak A}_\alpha\cap\lambda$ is an ordinal, and
\[E=:\{\alpha:{\mathfrak A}_\alpha\cap\lambda=\alpha\}.\]
Choose ${\mathfrak A}\prec  ({\Cal H}(\chi),\in,<^*_\chi)$, satisfying $\lambda^* +1
\subseteq {\mathfrak A}, \lambda^*= \parallel {\mathfrak A} \parallel, \,\mbox{ and }
\quad \bigcup\limits_{\alpha<\lambda} {\mathfrak A}_\alpha \prec {\mathfrak A}$

Note that $\Gamma \subseteq {\mathfrak A}$. Without loss of generality $T^*$
satisfies
\begin{description}
\item[$(*)_1$] $\lev_\alpha (T^*) \subseteq {}^\alpha\lambda$ for $\alpha
\le \omega$ and $x\le_T y$ means $x\le y$ and for $\eta\in T^*$ 
such that 
$\lg (\eta)=n$ let ${\rm Suc}_{T^*} (\eta)=\{\eta'\langle 
\alpha\rangle:\alpha< \lambda_n\}$
\item[$(*)_2$] for every $\eta\in\lev_\omega(T^*)$ and $k<\omega$, there are
$\lambda$ members $\nu$ of $\lev_\omega(T^*)$ such that $\nu\restriction k=
\eta\restriction k$
\end{description}

\noindent NOTE: By $\bU^{<\mu}_{J^2_{\bar m}}(\lambda)=\lambda^*$,
\begin{description}
\item[$(*)$] \underbar{if} $x_\eta\in\lev_\omega(T^*)$ for $\eta\in\Gamma
\subseteq
\prod\bar m$

\underbar{then} for some $A\in (J'_{\bar m})^+$ the set $\langle(\eta,
x_\eta):\eta\in A\rangle$ belongs to ${\mathfrak A}$, hence $\{x_\eta(n):\eta\in A$
and for some $\nu\in A$ we have $\lg(\eta\cap\nu)=n\}$ belongs to ${\mathfrak A}$.
Hence: if the mapping $\eta\mapsto x_\eta$ for $\eta\in\Gamma$
is continuous then $\langle x_\rho:\rho\in cl_\Gamma (t) \rangle\in {\mathfrak A}$.
\end{description}
We let
\[\begin{aligned}
P^0=P^0({\mathfrak A})=\biggl\{\bar x:&\bar x=\langle x_\rho:\rho
  \in t\rangle\in {\mathfrak A}\mbox{ and } x_\rho\in\lev_{\ell g(\rho)}(T^*),
t\in(J'_{\bar{m}})^* \mbox{ and }\\
  &\mbox{the mapping }\rho\mapsto x_\rho\mbox{ preserves all of the
  relations}\\
  &\ell g(\rho_1\cap\rho_2)=n \biggr\}.
\end{aligned}\]
Assume $\bar x=\langle x_\rho:\rho\in t\rangle\in P^0$.
For $\delta\in S$ such that $C_\delta \subseteq E$ let
\[\inv(\bar x,C_\delta,T^*,{\bar {\mathfrak A}})=:\big\{\alpha\in C_\delta:
(\exists\rho\in\Dom(\bar x))(x_\rho\in {\mathfrak A}_{\min(C_\delta
\setminus (\alpha+1))}\setminus {\mathfrak A}_\alpha\big)\}.\]
Let
\[\begin{array}{l}
\Inv(C_\delta,T^*,{\mathfrak A})=:\\
\big\{a:\mbox{for some }\bar x\in P^0,\  a =\inv(\bar x,C_\delta,T^*,{\mathfrak A})
\mbox{ is a countable set}\big\}.
\end{array}\]
Note: $\inv(\bar x,C_\delta,T,{\mathfrak A})$ has cardinality at most continuum,
so $\Inv(C_\delta,T^*,{\mathfrak A})$ is a family of $\le 2^{\aleph_0}\times
|{\mathfrak A}|=\lambda^*$ countable subsets of $C_\delta$.

We continue as before. Let $\alpha_{\delta,\varepsilon}$ be the
$\varepsilon$-th member of $C_\delta$ for $\varepsilon<\mu$. So as
$\lambda^*<\cov(\mu,\mu,\aleph_1,2)$  we can find $\gamma_n<\mu$ for
$n<\omega$ such that: $a\in{\mathfrak A} \cap [\mu]^{<\mu} \Rightarrow
\aleph_0>|\{\gamma_n:n <\omega\}\cap a|$ finite; how?  using coding 
of finite sequences.
Hence we can find
$\gamma_n\in (\bigcup\limits_{\ell<n}\lambda_\ell,\lambda_n)$ limit such
that
for each $\delta\in S$ and $a\in\Inv(C_\delta,T^*,\bar{\mathfrak A})$ we have
$\{\gamma_n+\ell:n<\omega\mbox{ and }\ell<m_i\}\cap \{\otp(\alpha\cap
C_\delta):\alpha\in a\}$ is bounded in $\mu$.

Let $\Xi_n=\{\rho\in\prod_{\bar{m}}: n=\sup\{\ell:\eta(\ell)\not=0\}.$

For each $\delta\in S^*=\{\delta\in S:C_\delta \subseteq
E\}$ we choose by induction on $n<\omega$ for each $\rho\in \Xi_n$, a member
$y_{\delta,\rho}\in \lev_\omega (T^*) \cap
({\mathfrak A}_{\alpha_\delta,\gamma_{\ell+1}}\setminus
{\mathfrak A}_{\alpha_\delta,\gamma_\ell})$
such that for $\rho_1,\rho_2\in
\bigcup\limits_{m\leq n} \Xi_n$ and $k$ we have $[\rho_1\restriction k=\rho_2
\restriction k \Leftrightarrow y_{\delta,\rho_1\restriction k}=
y_{\delta,\rho_2\restriction k}]$. Let for $\rho
\in \prod \bar{m},x^\delta_\rho$ be the unique $x\in
\prod\limits_{n<\omega} \lambda_n$ such that 
$\nu\in\bigcup\limits_n \Xi_n \& \rho \restriction k= \nu \restriction k
\Rightarrow x \restriction k= y_{\delta,\nu}$, so 
$x^\delta_\rho=y_{\delta,\rho}$ for $\rho\in\bigcup\limits_n \Xi_n$.

Now we can find $T$ such that $\lev_n(T)=\lev_n(T^*)$ for
$n<\omega$ and
\[\lev_\omega(T)=\big\{x^\delta_\rho: \delta\in S^* \mbox{ and } \rho\in\prod
\bar{m}\big\}\] 
So, if $T^*$ is universal then there is an embedding $f:T\longrightarrow T^*$,
define ${\mathfrak A}'_\alpha$ (for $\alpha<\lambda$) as before such that 
${\mathfrak A}_\delta$,
$\langle {\mathfrak A}_\alpha:\alpha<\lambda\rangle \in {\mathfrak A}'_0$ and $T,f\in
{\mathfrak A}'_0$, so ${\mathfrak A}'_\alpha\mbox{ is
closed under }f\mbox{ and } f^{-1}$ and hence
\[E'=\{\alpha\in E:{\mathfrak A}'_\alpha \cap \lambda=\alpha\}\]
is a club of $\lambda$. By the choice of $\bar C$ for some $\delta\in S$ we
have $C_\delta\subseteq E'$. Now use $(*)$ with $x_\eta=f(x^\delta_\eta
)$. Thus we get $A\in (J'_{\bar m})^+$ such that
$\{(\eta,x_\eta):\eta\in A\}\subseteq b\in{\mathfrak A}$, 
$|b|\leq \mu^*$. We can assume $\{(\eta,x_\eta):\eta\in A\}\in {\mathfrak A}$
so there is
$\nu\in\bigcup\limits_k
\prod\limits_{i<k} m_i$ such that $A$ is dense above $\nu$, hence as $f$ is
continuous, $\langle(\eta,x_\eta):\nu\triangleleft\eta\in\Gamma\rangle
\in {\mathfrak A}$. So $\langle x_\eta:\eta\in\Gamma,\nu\trianglelefteq\eta
\rangle\in P^0({\mathfrak A})$, and hence the set
\[\{\alpha_{\delta,\gamma_\ell}:\nu \triangleleft \rho_\ell\}\]
is included in $\inv(\bar x,C_\delta,T^*,{\mathfrak A})$. Hence
\[a=\{\alpha_{\delta,\gamma_\ell}\!:\nu\triangleleft \rho_\ell\}\in\Inv(
C_\delta,T^*,{\mathfrak A}),\]
contradicting

``$\{\alpha_{\delta,\gamma_\ell}:\ell<\omega\}$ has finite intersection with
any $a\in\Inv(C_\delta,T^*,{\mathfrak A})$''.

\begin{Remark}
\label{7.3}
We can a priori fix a set of $\aleph_0$ candidates and say more on their order
of appearance, so that $\Inv(\bar x,C_\delta,T^*,{\bar{\mathfrak A}})$ has order
type $\omega$. This makes it easier to phrase a true invariant, i.e. $\langle
(\eta_n,t_n):n<\omega\rangle$ is as above, $\langle\eta_n:n<\omega\rangle$
lists ${}^{\omega >}\omega$ with no repetition, $\langle t_n\cap {}^\omega
\omega:n<\omega\rangle$ are pairwise disjoint. If $x_\rho\in\lev_\omega(T^*)$
for $\rho\in {}^\omega\omega$, $\bar T^*=\langle\bar T^*_\zeta:\zeta<\lambda
\rangle$ representation
\[\begin{array}{l}
\inv(\langle x_\rho:\rho\in {}^\omega\omega\rangle,C_\delta,\bar T^*)=\\
 \big\{\alpha\in C_\delta:\mbox{for some } n,\ (\forall\rho)[\rho\in t_n\cap
 {}^\omega\omega\quad\Rightarrow\quad x_\rho\in T^*_{\min(C_\delta\setminus
 (\alpha+1))}\setminus T^*_\alpha]\big\}.
\end{array}\]
\hfill$\square_{\ref{7.2}}$
\end{Remark}

By \ref{7.2}  and \cite{Sh:460} we get

\begin{nConclusion}
\label{7.8}
If $\beth_\omega<\lambda_{n+1}<\mu=\sum\limits_n
\lambda_{n}$ and  $\mu^+<\lambda=\cf(\lambda)<\mu^{\aleph_0}$ and
$\bar{\lambda}=\langle \lambda_n:n<\omega\rangle \wedge \lambda$
\underline{then} in ${\mathfrak K}^{tr}_{\bar{\lambda}}$ there is no
universal member.
\end{nConclusion}
+-+
\begin{Remark}
\label{7.9}
\begin{enumerate}
\item We can in \ref{7.8} and earlies deal with simpler
cardinals, but it is more cumbersom.
\item We can give sufficient conditions for having $T_\alpha\in
{\mathfrak K}^{tr}_{\bar{\lambda}}$ for $\alpha<\lambda^{**}$, such that
into no $T^*\in{\mathfrak k}^{tr}_{\bar{\lambda}}$ can we embed
$>{\bf U}_{J^2_{m\restriction \Gamma}}(\lambda)$ of then
[see, more, see  \cite{Sh:622} end of section 2]
\item If we like to use e.g. $|\Gamma|<2^{\aleph_0}$, 
we need to deal with
a closure of a copy of $\Gamma$.
\end{enumerate}
\end{Remark}

\section{Universals in singular cardinals}
In \S3, \S5, \ref{7.2}, we can in fact deal with ``many'' singular cardinals
$\lambda$. This is done by proving a stronger assertion on some regular
$\lambda$. Here ${\mathfrak K}$ is a class of models.

\begin{Lemma}
\label{8.1}
\begin{enumerate}
\item There is no universal member in ${\mathfrak K}_{\mu^*}$ if for some $\lambda
<\mu^*$, $\theta\ge 1$ we have:
\begin{description}
\item[\mbox{$\otimes_{\lambda,\mu^*,\theta}[{\mathfrak K}]$}] not only there is no
universal member in ${\mathfrak K}_\lambda$ but if we assume:
\[\langle M_i:i<\theta\rangle\ \mbox{ is given, }\ \|M_i\|\le\mu^*<\prod_n
\lambda_n,\ M_i\in {\mathfrak K},\]
then there is a structure $M$ from ${\mathfrak K}_\lambda$ (in some cases
of a simple form) not embeddable in any $M_i$.
\end{description}
\item Assume
\begin{description}
\item[$\otimes^\sigma_1$]  $\langle\lambda_n:n<\omega\rangle$ is given,
$\lambda^{\aleph_0}_n<\lambda_{n+1}$,
\[\mu=\sum_{n<\omega}\lambda_n<\lambda=\cf(\lambda)\leq
\mu^*<\prod_{n<\omega}\lambda_n\]
and $\mu^+<\lambda$ or at least there is a club guessing $\bar C$ as
in $(**)^1_\lambda$ (ii) of \ref{3.4} for
$(\lambda,\mu)$.
\end{description}
\underbar{Then} there is no universal member in ${\mathfrak K}_{\mu^*}$ (and
moreover $\otimes_{\lambda,\mu^*,\theta}[{\mathfrak K}]$ holds) in the following
cases
\begin{description}
\item[$\otimes_2$(a)] for torsion free groups, i.e. ${\mathfrak K}={\mathfrak
K}^{rtf}_{\bar\lambda}$ if $\cov(\mu^*,\lambda^+,\lambda^+,\lambda)<
\prod\limits_{n<\omega}\lambda_n$, see notation \ref{0.4} on $\cov$)
\item[\quad(b)] for ${\mathfrak K}={\mathfrak K}^{tcf}_{\bar\lambda}$,
\item[\quad(c)] for ${\mathfrak K}={\mathfrak K}^{tr}_{\bar\lambda}$ as in \ref{7.2} -
$\cov(\bU_{J^3_{\bar m}}(\mu^*),\lambda^+,\lambda^+,\lambda)<\prod\limits_{n<
\omega}\lambda_n$,
\item[(d)] for ${\mathfrak K}^{rs(p)}_{\bar\lambda}$: like case (c) (for
appropriate ideals), replacing $tr$ by $rs(p)$.
\end{description}
\end{enumerate}
\end{Lemma}

\begin{Remark}
\label{8.1A}
\begin{enumerate}
\item For \ref{7.2} as $\bar m=\langle\omega:i<\omega\rangle$ it is clear that
the subtrees $t_n$ are isomorphic. We can use $m_i\in [2,\omega)$, and use
coding; anyhow it is immaterial since ${}^\omega \omega,{}^\omega 2$ are
similar.
\item We can also vary $\bar\lambda$ in \ref{8.1} $\otimes_2$, case (c).
\item We can replace $\cov$ in $\otimes_2$(a),(c) by
\[\sup\pp_{\Gamma(\lambda)}(\chi):\cf(\chi)=\lambda,\lambda<\chi\le
\bU_{J^3_{\bar m}}(\mu^*)\}\]
(see \cite[5.4]{Sh:355}, \ref{2.3}).
\end{enumerate}
\end{Remark}

\proof Should be clear, e.g.\\
{\em Proof of Part 2), Case (c)}\ \ \  Let $\langle T_i:i<i^* \rangle$ be
given, $i^*<\prod\limits_{n<\omega}\lambda_n$ such that
\[\|T_i\|\le\mu^*\quad\mbox{ and }\quad\mu^\otimes=:\cov(\bU_{J^3_{\bar
m}}(\mu^*),\lambda^+,\lambda^+,\lambda)<\prod_{n<\omega}\lambda_n.\]
By \cite[5.4]{Sh:355} and $\pp$ calculus (\cite[2.3]{Sh:355}), $\mu^\otimes=
\cov(\mu^\otimes,\lambda^+,\lambda^+,\lambda)$. Let
$\chi=\beth_7(\lambda)^+$. For $i<i^*$ choose ${\mathfrak B}_i\prec (H(\chi)\in
<^*_\chi)$, $\|{\mathfrak B}_i\|=\mu^\otimes$, $T_i\in {\mathfrak B}_i$, $\mu^\otimes
+1\subseteq {\mathfrak B}_i$. Let $\langle Y_\alpha: \alpha<\mu^\otimes\rangle$ be
a family of subsets of $T_i$ exemplifying the Definition of $\mu^\otimes=
\cov(\mu^\otimes,\lambda^+,\lambda^+,\lambda)$.\\
Given $\bar x=\langle x_\eta:\eta\in {}^\omega\omega\rangle$, $x_\eta\in
\lev_\omega(T_i)$, $\eta\mapsto x_\eta$ continuous (in our case this means
$\ell g(\eta_1\cap\eta_2)=\ell g(x_{\eta_1}\cap x_{\eta_2})=:\ell g(\max
\{\rho:\rho\triangleleft\eta_1\ \&\ \rho\triangleleft\eta_2\})$. Then for some
$\eta\in {}^{\omega>}\omega$,
\[\langle x_\rho:\eta\triangleleft\rho\in {}^\omega\omega\rangle\in{\mathfrak
B}.\]
So given $\left<\langle x^\zeta_\eta:\eta\in {}^\omega\omega\rangle:\zeta<
\lambda\right>$, $x^\zeta_\eta\in\lev_\omega(T_i)$ we can find $\langle
(\alpha_j,\eta_j):j<j^*<\lambda\rangle$ such that:
\[\bigwedge_{\zeta<\lambda}\bigvee_j\langle x^\zeta_\eta:\eta_j\triangleleft
\eta\in {}^\omega\omega\rangle\in Y_\alpha.\]
Closing $Y_\alpha$ enough we can continue as usual. 
\hfill$\square_{\ref{8.1}}$

\section{Metric spaces and implications}

\begin{Definition}
\label{9.1}
\begin{enumerate}
\item ${\mathfrak K}^{mt}$ is the class of metric spaces $M$ (i.e. $M=(|M|,d)$,
$|M|$ is the set of elements, $d$ is the metric, i.e. a two-place function
from $|M|$ to ${\mathbb R}^{\geq 0}$ such that
$d(x,y)=0\quad\Leftrightarrow\quad x=0$ and
$d(x,z)\le d(x,y)+d(y,z)$ and $d(x,y) = d(y,x)$).

An embedding $f$ of $M$ into $N$ is a one-to-one function from $|M|$ into
$|N|$ which is continuous, i.e. such that:
\begin{quotation}
\noindent if in $M$, $\langle x_n:n<\omega\rangle$ converges to $x$

\noindent then in $N$, $\langle f(x_n):n<\omega\rangle$ converges to $f(x)$.
\end{quotation}
\item ${\mathfrak K}^{ms}$ is defined similarly but $\Rang(d)\subseteq\{2^{-n}:n
<\omega\}\cup\{0\}$ and instead of the triangular inequality we require
\[d(x,y)=2^{-i},\qquad d(y,z)=2^{-j}\quad \Rightarrow\quad d(x,z) \le
2^{-\min\{i-1,j-1\}}.\]
\item ${\mathfrak K}^{tr[\omega]}$ is like ${\mathfrak K}^{tr}$ but $P^M_\omega=|M|$
and embeddings preserve $x\;E_n\;y$ (not necessarily its negation) are
one-to-one, and remember $\bigwedge\limits_n
x\;E_n\;y\quad\Rightarrow\quad x \restriction n = y \restriction n$).
\item ${\mathfrak K}^{mt(c)}$ is the class of semi-metric spaces $M=(|M|,d)$,
which means that for the constant $c\in\mathbb R^+$ the triangular inequality is
weakened to $d(x,z)\le cd(x,y)+cd(y,z)$ with embedding as in \ref{9.1}(1)
(so for $c=1$ we get ${\mathfrak K}^{mt}$).
\item ${\mathfrak K}^{mt[c]}$ is the class of pairs $(A,d)$ such that $A$ is a
non-empty set, $d$ a two-place symmetric function from $A$ to ${\mathbb
R}^{\ge 0}$ such that $[d(x,y)=0\ \Leftrightarrow\ x=y]$ and
\[d(x_0,x_n)\le c\sum\limits_{\ell<n} d(x_\ell,x_{\ell+1})\ \ \mbox{ for any
$n<\omega$ and $x_0,\ldots,x_n\in A$.}\]
\item ${\mathfrak K}^{ms(c)}$, ${\mathfrak K}^{ms[c]}$ are defined parallely.
\item ${\mathfrak K}^{rs(p),\mbox{pure}}$ is defined like ${\mathfrak K}^{rs(p)}$ but
the embeddings are pure.
\end{enumerate}
\end{Definition}

\begin{Remark}
There are, of course, other notions of embeddings; isometric embeddings if $d$
is preserved, co-embeddings if the image of an open set is open, bi-continuous
means an embedding which is a co-embedding.  The isometric embedding is the
weakest, its case is essentially equivalent to the ${\mathfrak
K}^{tr}_{\lambda}$ case (as in \ref{9.3}(3)); for the open case there
is a universal:
discrete space. The universal for ${\mathfrak K}^{mt}_\lambda$ under bicontinuous
case exist in cardinality $\lambda^{\aleph_0}$, see \cite{Ko57}.
\end{Remark}

\begin{Definition}
\label{9.1A}
\begin{enumerate}
\item $\Univ^0({\mathfrak K}^1,{\mathfrak K}^2)=\{(\lambda,\kappa,\theta):$ there are
$M_i\in {\mathfrak K}^2_\kappa$ for $i<\theta$ such that any $M\in {\mathfrak
K}^1_\lambda$ can be embedded into some $M_i\}$. We may omit $\theta$
if it is 1. We may omit the superscript 0.
\item $\Univ^1({\mathfrak K}^1,{\mathfrak K}^2)=\{(\lambda,\kappa,\theta):$ there are
$M_i\in {\mathfrak K}^2_\kappa$ for $i<\theta$ such that any $M\in {\mathfrak
K}^1_\lambda$ can be represented as the union of $<\lambda$ sets $A_\zeta$
($\zeta<\zeta^*<\lambda)$ such that each $M\restriction A_\zeta$ can be
embedded into some $M_i\}$ and is a $\leq_{{\mathfrak K}^1}$-submodel of $M$.
\item If above ${\mathfrak K}^1={\mathfrak K}^2$ we write it just once; (naturally we
usually assume ${\mathfrak K}^1 \subseteq {\mathfrak K}^2$).
\end{enumerate}
\end{Definition}

\begin{Remark}
\begin{enumerate}
\item We prove our theorems for $Univ^0$, we can get parallel things for
$\Univ^1$.
\item Many previous results of this paper can be rephrased using a pair of
classes.
\item We can make \ref{9.2} below deal with pairs and/or function $H$ changing
cardinality.
\item $\Univ^\ell$ has the obvious monotonicity properties.
\end{enumerate}
\end{Remark}

\begin{Proposition}
\label{9.2}
\begin{enumerate}
\item Assume ${\mathfrak K}^1,{\mathfrak K}^2$ has the same models as their members
and every embedding for ${\mathfrak K}^2$ is an embedding for ${\mathfrak K}^1$.\\
Then $\Univ({\mathfrak K}^2)\subseteq\Univ({\mathfrak K}^1)$.
\item Assume there is for $\ell=1,2$ a function $H_\ell$ from ${\mathfrak K}^\ell$
into ${\mathfrak K}^{3-\ell}$ such that:
\begin{description}
\item[(a)]  $\|H_1(M_1)\|=\|M_1\|$ for $M_1\in {\mathfrak K}^1$,
\item[(b)]  $\|H_2(M_2)\|=\|M_2\|$ for $M_2\in {\mathfrak K}^2$,
\item[(c)]  if $M_1\in {\mathfrak K}^1$, $M_2\in {\mathfrak K}^2$, $H_1(M_1)\in {\mathfrak
K}^2$ is embeddable into $M_2$ \underbar{then} $M_1$ is embeddable into
$H_2(M_2)\in {\mathfrak K}^1$.
\end{description}
\underbar{Then}\quad $\Univ({\mathfrak K}^2)\subseteq\Univ({\mathfrak K}^1)$.
\end{enumerate}
\end{Proposition}

\begin{Definition}
\label{9.2A}
We say ${\mathfrak K}^1\le {\mathfrak K}^2$ if the assumptions of \ref{9.2}(2)
hold. We say ${\mathfrak K}^1\equiv {\mathfrak K}^2$ if ${\mathfrak K}^1\le {\mathfrak K}^2
\le {\mathfrak K}^1$ (so larger means with fewer cases of universality).
\end{Definition}

\begin{Theorem}
\label{9.3}
\begin{enumerate}
\item The relation ``${\mathfrak K}^1\le {\mathfrak K}^2$'' is a quasi-order
(i.e. transitive and reflexive).
\item If $({\mathfrak K}^1,{\mathfrak K}^2)$ are as in \ref{9.2}(1) then ${\mathfrak K}^1
\le {\mathfrak K}^2$ (use $H_1 = H_2 =$ the identity).
\item For $c_1>1$ we have ${\mathfrak K}^{mt(c_1)}\equiv {\mathfrak K}^{mt[c_1]}\equiv
{\mathfrak K}^{ms[c_1]}\equiv {\mathfrak K}^{ms(c_1)]}$.
\item ${\mathfrak K}^{tr[\omega]} \le {\mathfrak K}^{rs(p)}$.
\item ${\mathfrak K}^{tr[\omega]} \le {\mathfrak K}^{tr(\omega)}$.
\item ${\mathfrak K}^{tr(\omega)} \le {\mathfrak K}^{rs(p),\mbox{pure}}$.
\end{enumerate}
\end{Theorem}

\proof 1)\ \ Check.\\
2)\ \  Check. \\
3)\ \  Choose $n(*)<\omega$ large enough and ${\mathfrak K}^1,{\mathfrak K}^2$ any two
of the four. We define $H_1,H_2$ as follows. $H_1$ is the identity. For $(A,d)
\in{\mathfrak K}^\ell$ let $H_\ell((A,d))=(A,d^{[\ell]})$ where $d^{[\ell]}(x,y)=
\inf\{1/(n+n(*)):2^{-n}\ge d(x,y)\}$ (the result is not necessarily a metric
space, $n(*)$ is chosen so that the semi-metric inequality holds). The point
is to check clause (c) of \ref{9.2}(2); so assume $f$ is a function which
${\mathfrak K}^2$-embeds $H_1((A_1,d_1))$ into $(A_2,d_2)$; but
\[H_1((A_1,d_1))=(A_1,d_1),\quad H_2((A_2,d_2))=(A_2,d^{[2]}_2),\]
so it is enough to check that $f$ is a function which ${\mathfrak K}^1$-embeds
$(A_1,d^{[1]}_1)$ into $(A_2,d^{[2]}_2)$ i.e. it is one-to-one (obvious) and
preserves limit (check).\\
4)\ \ For $M=(A,E_n)_{n<\omega}\in {\mathfrak K}^{tr[\omega]}$, without loss of
generality $A\subseteq {}^\omega\lambda$ and
\[\eta E_n\nu\qquad\Leftrightarrow\qquad\eta\in A\ \&\ \nu\in A\ \&\
\eta\restriction n=\nu\restriction n.\]
Let $B^+=\{\eta\restriction n:\eta\in A\mbox{ and } n<\omega\}$. We define
$H_1(M)$ as the (Abelian) group generated by
\[\{x_\eta:\eta\in A\cup B\}\cup\{y_{\eta,n}:\eta\in A,n<\omega\}\]
freely except
\[\begin{array}{rcl}
p^{n+1}x_\eta=0 &\quad\mbox{\underbar{if}}\quad &\eta\in B, \ell g(\eta)=n\\
y_{\eta,0}=x_\eta &\quad\mbox{\underbar{if}}\quad &\eta\in A\\
py_{\eta,n+1}-y_\eta=x_{\eta\restriction n} &\quad\mbox{\underbar{if}}\quad
&\eta\in A, n<\omega\\
p^{n+1}y_{\eta, n}=0 &\quad\mbox{\underbar{if}}\quad &\eta\in B, n<\omega.
\end{array}\]
For $G\in {\mathfrak K}^{rs(p)}$ let $H_2(G)$ be $(A,E_n)_{n<\omega}$ with:
\[A = G,\quad xE_ny\quad\underbar{iff}\quad G\models\mbox{``}p^n\mbox{ divides
}(x-y)\mbox{''.}\]
$H_2(G)\in {\mathfrak K}^{tr[\omega]}$ as ``$G$ is separable" implies $(\forall
x)(x \ne 0\ \Rightarrow\ (\exists n)[x\notin p^nG])$. Clearly clauses
(a), (b) of
Definition \ref{9.1}(2) hold. As for clause (c), assume $(A,E_n)_{n<\omega}
\in {\mathfrak K}^{tr[\omega]}$. As only the isomorphism type counts without loss
of generality $A\subseteq {}^\omega\lambda$. Let $B=\{\eta\restriction n:n<
\omega:\eta\in A\}$ and $G=H_1((A,E_n)_{n<\omega})$ be as above. Suppose that
$f$ embeds $G$ into some $G^*\in {\mathfrak K}^{rs(p)}$, and let $(A^*,E^*_n)_{n<
\omega}$ be $H_2(G^*)$. We should prove that $(A,E_n)_{n<\omega}$ is
embeddable into $(A^*,E^*_n)$.\\
Let $f^*:A\longrightarrow A^*$ be $f^*(\eta)=x_\eta\in A^*$. Clearly $f^*$ is
one to one from $A$ to $A^*$; if $\eta E_n \nu$ then $\eta\restriction n=\nu
\restriction n$ hence $G \models p^n \restriction (x_\eta-x_\nu)$ hence
$(A^*,A^*_n)_{n<\omega}\models\eta E^*_n\nu$. \hfill$\square_{\ref{9.3}}$

\begin{Remark}
\label{9.3A}
In \ref{9.3}(4) we can prove ${\mathfrak K}^{tr[\omega]}_{\bar\lambda}\le{\mathfrak
K}^{rs(p)}_{\bar\lambda}$.
\end{Remark}

\begin{Theorem}
\label{9.4}
\begin{enumerate}
\item ${\mathfrak K}^{mt} \equiv {\mathfrak K}^{mt(c)}$ for $c \ge 1$.
\item ${\mathfrak K}^{mt} \equiv {\mathfrak K}^{ms[c]}$ for $c > 1$.
\end{enumerate}
\end{Theorem}

\proof 1)\ \ Let $H_1:{\mathfrak K}^{mt}\longrightarrow {\mathfrak K}^{mt(c)}$ be the
identity. Let $H_2:{\mathfrak K}^{mt(c)}\longrightarrow {\mathfrak K}^{mt}$ be defined
as follows:\\
$H_2((A,d))=(A,d^{mt})$, where
\[\begin{array}{l}
d^{mt}(y,z)=\\
\inf\big\{\sum\limits^n_{\ell=0} d(x_\ell,x_{\ell,n}):n<\omega\ \&\
x_\ell\in A\mbox{ (for $\ell\le n$)}\ \&\ x_0=y\ \&\ x_n=z\big\}.
\end{array}\]
Now
\begin{description}
\item[$(*)_1$]  $d^{mt}$ is a two-place function from $A$ to ${\mathbb
R}^{\ge 0}$,
is symmetric, $d^{mt}(x,x)=0$ and it satisfies the triangular inequality.
\end{description}
This is true even on ${\mathfrak K}^{mt(c)}$, but here also
\begin{description}
\item[$(*)_2$]  $d^{mt}(x,y) = 0 \Leftrightarrow x=0$.
\end{description}
[Why? As by the Definition of ${\mathfrak K}^{mt[c]},d^{mt}(x,y)\ge{\frac 1c}
d(x,y)$. Clearly clauses (a), (b) of \ref{9.2}(2) hold.]\\
Next,
\begin{description}
\item[$(*)_3$]  If $M_1,N\in {\mathfrak K}^{mt}$, $f$ is an embedding (for ${\mathfrak
K}^{mt}$) of $M_1$ into $N$ then $f$ is an embedding (for ${\mathfrak K}^{mt[c]}$)
of $H_1(M)$ into $H_1(N)$
\end{description}
[why? as $H_1(M)=M$ and $H_2(N)=N$],
\begin{description}
\item[$(*)_4$]  If $M,N\in {\mathfrak K}^{mt[c]}$, $f$ is an embedding (for
${\mathfrak K}^{mt[c]}$) of $M$ into $N$ then $f$ is an embedding (for ${\mathfrak
K}^{mt}$) of $H_2(M)$ into $H_1(M)$
\end{description}
[why? as $H^*_\ell$ preserves $\lim\limits_{n\to\infty} x_n=x$ and
$\lim\limits_{n\to\infty} x_n\ne x$].

So two applications of \ref{9.2} give the equivalence. \\
2)\ \ We combine $H_2$ from the proof of (1) and the proof of \ref{9.3}(3).
\hfill$\square_{\ref{9.4}}$

\begin{Definition}
\label{9.6}
\begin{enumerate}
\item If $\bigwedge\limits_n\mu_n=\aleph_0$ let
\[\hspace{-0.5cm}\begin{array}{ll}
 J^{mt}=J^{mt}_{\bar\mu}=\big\{A\subseteq\prod\limits_{n<\omega}\mu_n:&
 \mbox{for every } n \mbox{ large enough, } \\
\ & \mbox{ for every }\eta\in\prod\limits_{\ell
 <n}\mu_\ell\\
\ &\mbox{the set }\{\eta'(n):\eta\triangleleft\eta'\in A\}\mbox{ is finite}
 \big\}.
\end{array}\]
\item Let $T=\bigcup\limits_{\alpha\le\omega}\prod\limits_{n<\alpha}\mu_n$,
 $(T,d^*)$ be a metric space such that
\[\prod_{\ell<n}\mu_\ell\cap\mbox{closure}\left(\bigcup_{m<n}\prod_{\ell<m}
\mu_\ell\right)=\emptyset;\]
now
\[\begin{aligned}
I^{mt}_{(T,d^*)}=:\big\{A\subseteq\prod\limits_{n<\omega}\mu_n:&\mbox{ for
  some }n,\mbox{ the closure of } A\mbox{ (in $(T,d^*)$)}\\
  &\mbox{ is disjoint to }\bigcup\limits_{m\in [n,\omega)}\prod\limits_{\ell
  <m} \mu_\ell\big\}.
\end{aligned}\]
\item Let $H\in {\mathfrak K}^{rs(p)}$, $\bar H=\langle H_n:n<\omega\rangle$, $H_n
\subseteq H$ pure and closed, $n<m\ \Rightarrow\ H_n\subseteq H_m$ and
$\bigcup\limits_{n<\omega} H_n$ is dense in $H$.  Let
\[\begin{array}{ll}
I^{rs(p)}_{H,\bar H}=:\big\{A\subseteq H:&\mbox{for some } n\mbox{ the closure
of }\langle A\rangle_H\mbox{ intersected with}\\
 &\bigcup\limits_{\ell<\omega}H_\ell\mbox{ is included in }H_n\big\}.
\end{array}\]
\end{enumerate}
\end{Definition}

\begin{Proposition}
\label{9.5}
Suppose that $2^{\aleph_0}<\mu$ and $\mu^+<\lambda=\cf(\lambda)<
\mu^{\aleph_0}$ and
\begin{description}
\item[$(*)_\lambda$] $\bU_{J^{mt}_{\bar\mu}}(\lambda)=\lambda$ or at least
$\bU_{J^{mt}_{\bar\mu}}(\lambda)<\lambda^{\aleph_0}$ for some $\bar\mu=\langle
\mu_n:n<\omega\rangle$ such that $\prod\limits_{n<\omega}\mu_n<\lambda$.
\end{description}
Then ${\mathfrak K}^{mt}_\lambda$ has no universal member.
\end{Proposition}

\begin{Proposition}
\label{9.7}
\begin{enumerate}
\item $J^{mt}$ is $\aleph_1$-based.
\item The minimal cardinality of a set which is not in the
$\sigma$-ideal generated by $J^{mt}$ is ${\mathfrak b}$.
\item $I^{mt}_{(T,d^*)},I^{rs(p)}_{H,\bar H}$ are $\aleph_1$-based.
\item $J^{mt}$ is a particular case of $I^{mt}_{(T,d^*)}$ (i.e. for some
choice of $(T,d^*)$).
\item $I^0_{\bar \mu}$ is a particular case of $I^{rs(p)}_{H,\bar H}$.
\end{enumerate}
\end{Proposition}

\proof of \ref{9.5}. Let
$$
\begin{array}{ll}
T_\alpha=\{(\eta, \nu)\in{}^\alpha\lambda\times {}^\alpha(\omega+
1):& \mbox{ for every }n\mbox{ such that }n+1< \alpha\\
 \ & \mbox{ we have
}\nu(n)< \omega\}
\end{array}
$$
and for $\alpha\le\omega$ let $T=\bigcup\limits_{\alpha\le
\omega}T_\alpha$. We define on $T$ the relation $<_T$:
\[(\eta_1,\nu_1)\le(\eta_1,\nu_2)\quad\mbox{ iff }\quad\eta_1\trianglelefteq
\eta_2\ \&\ \nu_1\triangleleft\nu_2.\]
We define a metric:\\
if $(\eta_1,\nu_1)\ne(\eta_2,\nu_2)\in T$ and $(\eta,\nu)$ is their maximal
common initial segment and $(\eta,\nu)\in T$ then necessarily $\alpha=
\ell g((\eta,\nu))<\omega$ and we let:
\begin{quotation}
\noindent if $\eta_1(\alpha)\ne\eta_2(\alpha)$ then
\[d\left((\eta_1,\nu_1),(\eta_2,\nu_2)\right)=2^{-\sum\{\nu(\ell):\ell<
\alpha\}},\]
if $\eta_1(\alpha)=\eta_2(\alpha)$ (so $\nu_1(\alpha)\ne\nu_2(\alpha)$ then
\[d\left((\eta_1,\nu_1),(\eta_2,\nu_2)\right)=2^{-\sum\{\nu(\ell):\ell<\alpha
\}}\times 2^{-\min\{\nu_1(\alpha),\nu_2(\alpha)\}}.\]
\end{quotation}
Now, for every $S\subseteq\{\delta<\lambda:\cf(\delta)=\aleph_0\}$, and $\bar
\eta=\langle\eta_\delta:\delta\in S\rangle$, $\eta_\delta\in {}^\omega
\delta$, $\eta_\delta$ increasing let $M_\eta$ be $(T,d)\restriction A_{\bar
\eta}$, where
\[A_{\bar\eta}=\bigcup_{n<\omega} T_n\cup\{(\eta_\delta,\nu):\delta\in S,\;\nu
\in {}^\omega\omega\}.\]
The rest is as in previous cases (note that $\langle(\eta\char 94\langle
\alpha \rangle,\nu\char 94\langle n\rangle):n<\omega\rangle$ converges to
$(\eta\char 94\langle\alpha\rangle,\nu\char 94\langle\omega\rangle)$
and even if $(\eta\char 94\langle \alpha\rangle, \nu\char 94\langle
n\rangle)\leq (\eta_n, \nu_n)\in T_\omega$ then $\langle(\eta_n,
\nu_n): n<\omega\rangle$ converge to $(\eta\char 94 \langle
\alpha\rangle, \nu\char 94\langle \omega\rangle)$).
\hfill$\square_{\ref{9.7}}$

\begin{Proposition}
\label{9.8}
If $\IND_{\chi'}(\langle\mu_n:n<\omega\rangle)$, then $\prod\limits_{n<\omega}
\mu_n$ is not the union of $\le\chi$ members of $I^0_{\bar\mu}$ (see
Definition \ref{5.5A} and Theorem \ref{5.5}).
\end{Proposition}

\proof Suppose that $A_\zeta=\{\sum\limits_{n<\omega}
p^nx^n_{\alpha_n}:\langle \alpha_n:n<\omega\rangle\in X_\zeta\}$ and
$\alpha_n<\mu_n$ are such that if $\sum
p^nx^n_{\alpha_n}\in A_\zeta$ then for infinitely many $n$ for every
$k<\omega$ there is $\langle \beta_n:n<\omega\rangle$,
\[(\forall\ell<k)[\alpha_\ell=\beta_\ell\ \ \Leftrightarrow\ \ \ell=n]\qquad
\mbox{ and }\qquad\sum_{n<\omega}p^nx^n_{\beta_n}\in A_\zeta\ \ \mbox{ (see
\S5).}\]
This clearly follows. \hfill$\square_{\ref{9.8}}$

\section{On Modules}
Here we present the straight generalization of the one prime case like Abelian
reduced separable $p$-groups. This will be expanded in \cite{Sh:622}
(including the proof of \ref{10.4new}). (see ??)

\begin{Hypothesis}
\label{10.1}
\begin{description}
\item[(A)] $R$ is a ring, $\bar{\mathfrak e}=\langle{\mathfrak e}_n:n<\omega\rangle$,
${\mathfrak e}_n$ is a definition of an additive subgroup of $R$-modules by an
existential positive formula (finitary or infinitary) decreasing with $n$, we
write ${\mathfrak e}_n(M)$ for this additive subgroup, ${\mathfrak e}_\omega(M)=
\bigcap\limits_n {\mathfrak e}_n(M)$.
\item[(B)] ${\mathfrak K}$ is the class of $R$-modules.
\item[(C)] ${\mathfrak K}^*\subseteq {\mathfrak K}$ is a class of $R$-modules, which
is closed under direct summand, direct limit and for which there is $M^*$, $x^*
\in M^*$, $M^*=\bigoplus\limits_{\ell\le n} M^*_\ell\oplus M^{**}_n$, $M^*_n
\in {\mathfrak K}$, $x^*_n\in {\mathfrak e}_n(M^*_n)\setminus {\mathfrak e}_{n+1}(M^*)$,
$x^*-\sum\limits_{\ell<n} x^*_\ell\in {\mathfrak e}_n(M^*)$.
\end{description}
\end{Hypothesis}

\begin{Definition}
\label{10.2}
For $M_1,M_2\in {\mathfrak K}$, we say $h$ is a $({\mathfrak K},\bar {\mathfrak
e})$-homomorphism from $M_1$ to $M_2$ if it is a homomorphism and it maps $M_1
\setminus {\mathfrak e}_\omega(M_1)$ into $M_2\setminus {\mathfrak
e}_\omega(M_2)$;
we say $h$ is an $\bar {\mathfrak e}$-pure homomorphism if for each $n$ it
maps $M_1\setminus {\mathfrak e}_n(M_1)$ into $M_2\setminus {\mathfrak
e}_n(M_2)$.
\end{Definition}

\begin{Definition}
\label{10.3}
\begin{enumerate}
\item Let $H_n \subseteq H_{n+1} \subseteq H$, $\bar H=\langle H_n:n<\omega
\rangle$, $c\ell$ is a closure operation on $H$, $c\ell$ is a function from
${\Cal P}(H)$ to itself and
\[X \subseteq c \ell(X) = c \ell(c \ell(X)).\]
Define
\[I_{H,\bar H,c\ell}=\big\{A\subseteq H:\mbox{for some }k<\omega\mbox{ we have
} c\ell(A)\cap\bigcup_{n<\omega} H_n\subseteq H_k\big\}.\]
\item We can replace $\omega$ by any regular $\kappa$ (so $H=\langle H_i:i<
\kappa\rangle$).
\end{enumerate}
\end{Definition}

\begin{Claim}
\label{10.4new}
Assume $|R|+\mu^+< \lambda = \cf(\lambda)< \mu^{\aleph_0}$, then for
every $M\in {\mathfrak K}_\lambda$ there is $N\in {\mathfrak K}_\lambda$ with
no $\bar {\mathfrak e}$-pure homomorphism from $N$ into $M$.
\end{Claim}

\begin{Remark}
In the interesting cases $c\ell$ has infinitary character.\\
The applications here are for $\kappa=\omega$. For the theory, $pcf$
is nicer for higher $\kappa$.
\end{Remark}

\section{Open problems}

\begin{Problem}
\begin{enumerate}
\item If $\mu^{\aleph_0}\ge\lambda$ then any $(A,d)\in {\mathfrak K}^{mt}_\lambda$
can be embedded into some $M'\in {\mathfrak K}^{mt}_\lambda$ with density
$\le\mu$.
\item If $\mu^{\aleph_0}\ge\lambda$ then any $(A,d)\in {\mathfrak K}^{ms}_\lambda$
can be embedded into some $M'\in {\mathfrak K}^{ms}_\lambda$ with density
$\le\mu$.
\end{enumerate}
\end{Problem}

\begin{Problem}
\begin{enumerate}
\item Other inclusions on $\Univ({\mathfrak K}^x)$ or show consistency of non
inclusions (see \S9).
\item Is ${\mathfrak K}^1\le {\mathfrak K}^2$ the right partial order? (see \S9).
\item By forcing reduce consistency of $\bU_{J_1}(\lambda)>\lambda+
2^{\aleph_0}$ to that of $\bU_{J_2}(\lambda)>\lambda+2^{\aleph_0}$.
\end{enumerate}
\end{Problem}

\begin{Problem}
\begin{enumerate}
\item The cases with the weak $\pcf$ assumptions, can they be resolved in
ZFC? (the $pcf$ problems are another matter).
\item Use \cite{Sh:460}, \cite{Sh:513} to get ZFC results for large enough
cardinals.
\end{enumerate}
\end{Problem}

\begin{Problem}
If $\lambda^{\aleph_0}_n<\lambda_{n+1}$, $\mu=\sum\limits_{n<\omega}
\lambda_n$, $\lambda=\mu^+<\mu^{\aleph_0}$ can $(\lambda,\lambda,1)$ belong to
$\Univ({\mathfrak K})$?  For ${\mathfrak K}={\mathfrak K}^{tr},{\mathfrak K}^{rs(p)},{\mathfrak
K}^{trf}$?
\end{Problem}

\begin{Problem}
\begin{enumerate}
\item If $\lambda=\mu^+$, $2^{<\mu}=\lambda<2^\mu$ can $(\lambda,\lambda,1)
\in \Univ({\mathfrak K}^{\mbox{or}}=$ class of linear orders)?
\item Similarly for $\lambda=\mu^+$, $\mu$ singular, strong limit, $\cf(\mu)=
\aleph_0$, $\lambda<\mu^{\aleph_0}$.
\item Similarly for $\lambda=\mu^+$, $\mu=2^{<\mu}=\lambda^+ <2^\mu$.
\end{enumerate}
\end{Problem}

\begin{Problem}
\begin{enumerate}
\item Analyze the existence of universal member from ${\mathfrak
K}^{rs(p)}_\lambda$, $\lambda<2^{\aleph_0}$.
\item \S4 for many cardinals, i.e. is it consistent that:
$2^{\aleph_0}> \aleph_\omega$ and for every $\lambda< 2^{\aleph_0}$
there is a universal member of ${\mathfrak K}^{rs(p)}_\lambda$?
\end{enumerate}
\end{Problem}



\begin{Problem}
\begin{enumerate}
\item If there are $A_i\subseteq\mu$ for $i<2^{\aleph_0}$, $|A_i\cap A_j|<
\aleph_0$, $2^\mu=2^{\aleph_0}$ find forcing adding $S\subseteq [{}^\omega
\omega]^\mu$ universal for $\{(B, \vartriangleleft):
{}^{\omega>}\omega \subseteq B\subseteq {}^{\omega\geq}\omega, |B|\leq
\lambda\}$ under (level preserving) natural embedding.
\end{enumerate}
\end{Problem}

\begin{Problem}
For simple countable $T$, $\kappa=\kappa^{<\kappa}<\lambda\subseteq \kappa$
force existence of universal for $T$ in $\lambda$ still $\kappa=\kappa^{<
\kappa}$ but $2^\kappa=\chi$.
\end{Problem}

\begin{Problem}
Make \cite[\S4]{Sh:457}, \cite[\S1]{Sh:500} work for a larger class of
theories more than simple.
\end{Problem}

See on some of these problems \cite{Sh:614}, \cite{Sh:622}.

\bibliographystyle{amsalpha}
\bibliography{shlhetal}
\end{document}